\documentclass[12pt]{article}
\usepackage[margin=1in]{geometry}

\let\originalleft\left
\let\originalright\right
\renewcommand{\left}{\mathopen{}\mathclose\bgroup\originalleft}
\renewcommand{\right}{\aftergroup\egroup\originalright}

\usepackage{amssymb}
\usepackage{amsmath}
\usepackage{mathtools}
\usepackage{xfrac}

\usepackage{amsthm}
\usepackage{bookmark}
\usepackage{hyperref}

\numberwithin{equation}{section}
\newtheorem{theorem}[equation]{Theorem}
\newtheorem{lemma}[equation]{Lemma}
\newtheorem{corollary}[equation]{Corollary}
\newtheorem{proposition}[equation]{Proposition}
\newtheorem{definition}[equation]{Definition}
\newtheorem{remark}[equation]{Remark}
\newtheorem{example}[equation]{Example}

\DeclareMathOperator{\dom}{dom}
\DeclareMathOperator{\End}{End}
\DeclareMathOperator{\rank}{rank}
\DeclareMathOperator{\grad}{grad}
\DeclareMathOperator{\Hom}{Hom}

\begin{document}

\title{The Non-Abelian X-Ray Transform on Asymptotically Hyperbolic Spaces}
\author{Haim Grebnev}
\maketitle

\begin{abstract}

In this paper we formulate and prove a gauge equivalence for unitary connections and skew-Hermitian Higgs fields of suitable regularity that are mapped to the same function under the non-abelian X-ray transform on nontrapping asymptotically hyperbolic spaces with negative curvature and no nontrivial twisted conformal Killing tensor fields with certain regularity.
If one furthermore fixes such a connection with zero curvature, a corollary provides an injectivity result for the non-abelian X-ray transform over skew-Hermitian Higgs fields.
\end{abstract}

\tableofcontents

\section{Conventions/notations} \label{Section_1} 

In this paper we employ the following conventions/notations:

	\begin{enumerate} \item We employ the Einstein summation convention.

	\item The dimension of our manifold will be $n+1$, and we will denote all indices related to its dimension by $0,1,\ldots,n$.
When using the Einstein summation convention on indices related to the manifold's dimension, we employ the convention that Latin indices run from 0 to $n$ while Greek indices run from 1 to $n$.

	\item The notation $C^k$ denotes $k$ times continuously differentiable objects.
If the object is scalar valued, we always assume that it is complex valued (e.g. $C^k \left(\overline{M} \right)=C^k \left(\overline{M};\mathbb{C} \right)$).

	\item Whenever we say ``smooth," we mean ``$C^\infty$." All diffeomorphisms are smooth.

	\item If $\pi  : E\rightarrow N$ is a vector bundle over a manifold $N$ and $S\subseteq N$ is a subset of $N$, then we let $\left.E\right|_S$ denote the restriction of $E$ to the fibers over $S$ (i.e. more precise notation would be $\pi ^{-1}\left[S\right]$).

	\item Continuing off of point 5), we write $C^\infty \left(N;E \right)$ for smooth sections of $E$ (i.e. not simply smooth maps from $N$ to $E$).

	\item We denote the geodesic vector field over the tangent and cotangent bundles (i.e. infinitesimal generator of the geodesic flow) by $X$.
We will also let $X$ denote its restriction to the sphere and cosphere bundles $SM$ and $S^\ast M$ (described below), which makes sense because $X$ is tangent to it.\end{enumerate}
\section{Introduction} \label{Section_2} 

\subsection{Motivation} \label{Section_2.1} 

We begin by providing motivation for the non-abelian X-ray transform - the object of our interest, delaying precise definitions for a later section.
This transform is a generalization of the so-called ``scalar X-ray transform," the latter of which is used in reconstructing images of the internals of patients after irradiating them with X-rays at various angles.
The typical mathematical problem for the scalar X-ray transform is the following: suppose that we have a bounded subset $D\subseteq\mathbb{R}^n$ with smooth boundary and a continuous function $\phi : D\rightarrow \left(0,\infty \right)$ over it.
In our analogy, $D$ represents the shape of our patient and $\phi$ the body's X-ray absorption coefficient at various points.
Suppose we have a parametrized line $l \left(t \right)$ that enters $D$ at $t=0$ and leaves at $t=t_{\mathrm{exit}}$, which represents the motion of an X-ray moving through the body.
The ray's intensity $I \left(t \right)$ along this path decays according to the law
\begin{equation} \label{MathItem_2.1} \frac{dI}{dt}=-\phi I,\ \ \ \ \ \ \ \ \ \ I \left(0 \right)=I_0,\end{equation}
where $I_0$ represents the initial intensity of the ray.
We record $I \left(t_{\mathrm{exit}} \right)$ (i.e. the intensity of the ray when it exits) and then repeat this procedure for all possible lines $l$ that pass through $D$, using the same value for $I_0$ every time.
The inverse problem is then to recover the coefficient $\phi$ from knowing such data, which is equivalent to recovering a gray-scale image of the patient's internals.
Immediately we may note that a necessary condition for this to be possible is that two different $\phi$'s cannot generate the same observation data.
If this is the case, we say that the scalar X-ray transform that takes $\phi$ to the observed data is ``injective." It turns out that this is the case: one can indeed recover $\phi$ from the mentioned data - see for instance Theorem 1.1.6 in \cite{Bibitem_42}.

We mention that in the literature it is standard to use a slightly different formulation of the mentioned X-ray transform.
We have that (\ref{MathItem_2.1}) is a separable differential equation; solving for $I$ gives that
$$I \left(t_{\mathrm{exit}} \right)=I_0e^{-\int_{0}^{t_{\mathrm{exit}}}{\phi\circ l \left(t \right)dt}}.$$
From here we see that knowing $I \left(t_{\mathrm{exit}} \right)$ and $\int_{0}^{t_{\mathrm{exit}}}{\phi\circ l \left(t \right)dt}$ is equivalent, and so it is standard to instead study the map
\begin{equation} \label{MathItem_2.2} \phi\longmapsto\left\{\int_{0}^{t_{\mathrm{exit}}}{\phi\circ l \left(t \right)dt} : \mathrm{for\ all\ possible\ lines} \ l\right\},\end{equation}
which is called the \textit{scalar X-ray transform} (typically using some parametrization for the set of lines $l$).
This form of the transform has the advantage of being linear in $\phi$.
However, we will not make use of this formulation because it does not suit the way we generalize the above scalar problem to the non-abelian X-ray transform.

The non-abelian X-ray transform is defined by taking (\ref{MathItem_2.1}) and turning it into a system of equations by letting $I$ by a column vector and $\phi$ a square matrix.
The question is also whether one can recover the matrix $\phi$ from the collected data, or in other words if the transform involved is injective.
One application of this is in the recently introduced \textit{polarimetric neutron tomography}, which attempts to reconstruct the structure of magnetic fields inside materials after sending neutron beams through them and measure their spin changes - see for instance \cite{Bibitem_22} and \cite{Bibitem_9}.
We will mention a few more applications of this problem with references in Section \ref{Section_2.6} below.
We also remark that it is not possible to define an analogous linear transform as (\ref{MathItem_2.2}) once we set $I$ to be a column vector and $\phi$ to be a matrix because the equation (\ref{MathItem_2.1}) is not necessarily separable anymore.
We note that even defining $I$ to be a matrix will not force (\ref{MathItem_2.1}) to be separable due to the non-commutativity of matrices, which is the motivation behind calling this the \textit{non-abelian} X-ray transform.

We will actually be interested in a more sophisticated generalization of this problem, for instance by allowing the paths ``$l$" to be geodesics with respect to some Riemannian metric $g$ on $D$ where $D$ is now a smooth manifold.
We will also formulate $I$ and $\phi$ to be a section and endomorphism field respectively of a smooth vector bundle over $M$ and generalize the time derivative to be a connection (which may also be unknown) in the direction of the curve's velocity.
It turns out that in this more general setting it is not possible to recover all of the coefficients involved except in a special case when an additional assumption is made on the connection - see Theorem \ref{MathItem_2.8} and Corollary \ref{MathItem_2.11} below.
However, in the case when we cannot recover the coefficients, we do have a ``gauge equivalence" understanding of coefficients that produce the same data.

Such transforms have been well studied in the context of compact domains, and so a direction of research that has received significant attention in recent years is whether its possible to generalize X-ray results to noncompact domains.
We will be pursuing this direction of research, in particular we will generalize known results for the non-abelian X-ray transform to a certain class of noncompact manifolds called ``asymptotically hyperbolic spaces."
\subsection{Asymptotically Hyperbolic Spaces} \label{Section_2.2} 

In this section we introduce the geometry on which our transform will be defined.
Let $\overline{M}$ be a compact smooth manifold with smooth boundary of dimension $n+1$ with $n\geq1$, whose interior we denote by $M$.
A Riemannian metric $g$ on $M$ is called \textbf{asymptotically hyperbolic (AH)} if for any (and hence all) boundary defining function $\rho : \overline{M}\rightarrow\left[0,\infty \right)$ (i.e. $\rho$ is smooth, $\rho=0$ on and only on $\partial \overline{M}$, and $\left.d\rho\right|_{\partial \overline{M}}\neq0$) the tensor $\overline{g}=\rho^2g$ extends to a smooth Riemannian metric on all of $\overline{M}$ with $\left|d\rho\right|_{\rho^2g}^2\equiv1$ along $\partial \overline{M}$.
The boundary $\partial \overline{M}$ is thought of as the ``infinity" where the metric $g$ blows up.
Hence recalling that hyperbolic space has constant sectional curvature $-1$, the known fact that the sectional curvatures of $g$ tend to $-1$ as one approaches $\partial \overline{M}$ (Proposition 1.10 in \cite{Bibitem_31}) explains why such metrics are given the name ``AH.\footnote{ More generally, the sectional curvatures approaches $-\left|d\rho\right|^2$ restricted to the boundary.}"

In fact, the Poincar\'e ball model of hyperbolic space is the archetypical example of an AH space.
It is given by $\overline{M}=\left\{\left|x\right|\le1\right\}\subseteq\mathbb{R}^{n+1}$, where $\left|x\right|$ denotes the Euclidean length, and
$$g=4\frac{ \left(dx^1 \right)^2+\ldots+ \left(dx^{n+1} \right)^2}{ \left(1-\left|x\right|^2 \right)^2}.$$
Indeed if one takes the boundary defining function $\rho=1-\left|x\right|^2$, then an elementary exercise shows that $\left|d\rho\right|_{\rho^2g}^2\equiv1$ along $\left\{\left|x\right|=1\right\}$.

Given an AH metric on $M$, it is proven in \cite{Bibitem_16} that there exist infinitely many boundary defining functions $\rho$ such that $\left|d\rho\right|_{\rho^2g}^2\equiv1$ on the neighborhood $\left\{\rho<\varepsilon\right\}$ of $\partial \overline{M}$ for some $\varepsilon>0$ (i.e. not simply on $\partial \overline{M}$).
Fixing such a $\rho$, if $ \left(y^1,\ldots,y^n \right)$ are coordinates of $\partial \overline{M}$, then it follows from the theory of Fermi coordinates applied to $\rho^2g$ (see Corollary 6.42 in \cite{Bibitem_26}) that
\begin{equation} \label{MathItem_2.3} g=\frac{d\rho^2+h_{\mu\nu}dy^\mu d y^\nu}{\rho^2}.\end{equation}
on $\left[0,\varepsilon \right)\times\dom{ \left(y^\mu \right)}$.
We call such boundary coordinates $ \left(\rho,y^\mu \right)$ of $\overline{M}$ \textbf{asymptotic boundary normal coordinates}.
Such functions $\rho$ are called \textbf{geodesic boundary defining functions} (the reason is that the curves $t\mapsto \left(t,y^\mu \right)$ are geodesics of $\rho^2g$).
Throughout the paper we will often assume that our boundary coordinates are of this form because several results that we cite from \cite{Bibitem_15} are only stated in such coordinates.

By Proposition 1.8 in \cite{Bibitem_31}, AH spaces are complete.
Furthermore, in some cases we will assume that $g$ is also \textbf{nontrapping} which means that for any complete $g$-geodesic $\gamma :  \left(-\infty,\infty \right)\rightarrow M$, $\liminf_{t\rightarrow\pm\infty}{\rho \left(\gamma \left(t \right) \right)}=0$.
Intuitively speaking, this condition requires that $\gamma$ eventually ``escapes to infinity."
\subsection{Results} \label{Section_2.3} 

We now state our results.
Suppose that $ \left(M\subseteq\overline{M},g \right)$ is an asymptotically hyperbolic (AH) space and that $\rho$ is a boundary defining function.
Assume that we have a smooth complex Hermitian vector bundle $ \left(\mathcal{E},\langle \cdot,\cdot\rangle _\mathcal{E} \right)$ over $\overline{M}$, meaning that $\langle \cdot,\cdot\rangle _\mathcal{E}$ is an inner product on every fiber $\mathcal{E}_x$ that varies smoothly over $\overline{M}$.
Suppose also that we have a smooth section of the endomorphism bundle $\Phi : \overline{M}\rightarrow\End{\mathcal{E}}$.
We will often require that $\Phi$ is skew-Hermitian with respect to $\langle \cdot,\cdot\rangle _\mathcal{E}$, which means that it satisfies $\langle \Phi u,v\rangle _\mathcal{E}=-\langle u,\Phi v\rangle _\mathcal{E}$.
If this is the case, we will write $\Phi\in C^\infty \left(\overline{M};{\End}_{\mathrm{sk}}{\mathcal{E}} \right)$.
Lastly, assume that we have a smooth connection $\nabla^\mathcal{E}$ in $\mathcal{E}$ over $\overline{M}$.
We will often require that $\nabla^\mathcal{E}$ is unitary with respect to $\langle \cdot,\cdot\rangle _\mathcal{E}$, which means that
$${V\langle u,v\rangle }_\mathcal{E}=\langle \nabla_V^\mathcal{E}u,v\rangle _\mathcal{E}+\langle u,\nabla_V^\mathcal{E}v\rangle _\mathcal{E},$$
whenever $V$ is a smooth vector field over $\overline{M}$ and $u,v$ are smooth sections of $\mathcal{E}$.
It is easy to check that if $ \left(x^i \right)$ are coordinates of $\overline{M}$, $ \left(b_i \right)$ are an \textit{orthonormal} frame, and ${}^\mathcal{E}\Gamma_{ij}^k$ are the connection symbols of $\nabla^\mathcal{E}$ with respect to $ \left(\sfrac{\partial }{\partial  x^i} \right)$ and $ \left(b_i \right)$, then $\nabla^\mathcal{E}$ is unitary if and only if the connection symbols satisfy the skew-symmetry property ${}^\mathcal{E}\Gamma_{ij}^k=-{}^\mathcal{E}\Gamma_{ik}^j$.

Take any unit-speed complete geodesic $\gamma :  \left(-\infty,\infty \right)\rightarrow M$ such that $\liminf_{t\rightarrow\pm\infty}{\rho \left(\gamma \left(t \right) \right)}=0$.
It follows from Lemma 2.3 in \cite{Bibitem_15} that the limits of $\gamma \left(t \right)$ in $\overline{M}$ exist as $t\rightarrow\pm\infty$ and is equal to a point on the boundary $\partial \overline{M}$.
The analog of (\ref{MathItem_2.1}) that we will be considering is the following initial value problem for a section $u :  \left(-\infty,\infty \right)\rightarrow\mathcal{E}$:
\begin{equation} \label{MathItem_2.4} \nabla_{\dot{\gamma} \left(t \right)}^\mathcal{E}u \left(t \right)+\Phi \left(\gamma \left(t \right) \right)u \left(t \right)=0,\ \ \ \ \lim_{t\rightarrow-\infty}{u \left(t \right)}=e,\end{equation}
where $e$ is any element in $\mathcal{E}_{x_0}$ where $x_0\in\partial \overline{M}$ is the limit of $\gamma \left(t \right)$ as $t\rightarrow-\infty$.
The data point that we ``record" is
\begin{equation} \label{MathItem_2.5} \lim_{t\rightarrow\infty}{u \left(t \right)}.\end{equation}
The question that we are interested then becomes whether we can recover $\Phi$ and $\nabla^\mathcal{E}$ from the data recorded for all such possible pairs $\gamma$ and $e$.

A bit of vocabulary: (\ref{MathItem_2.4}) is a type of differential equation called a \textbf{transport equation}, and $\Phi$ is called a \textbf{Higgs field}.
Going from the pair $ \left(\nabla^\mathcal{E},\Phi \right)$ to the map that takes every $ \left(\gamma,e \right)$ as above to its associated data (\ref{MathItem_2.5}) is called the \textbf{non-abelian X-ray transform}, of which we give a more precise definition in Section \ref{Section_2.5} below.

To make rigorous sense of our problem however, we need to establish the well-definedness of the solution to (\ref{MathItem_2.4}) and the data (\ref{MathItem_2.5}).
Considering that we are making use of the values of the solution to (\ref{MathItem_2.4}) at plus or minus infinities, we accomplish this by imposing a decay condition on $\Phi$.
Part 1) in the following lemma establishes the well-definedness of our problem; part 2) is a result that we will need later in the paper.

\begin{lemma} \label{MathItem_2.6}  Suppose that $ \left(M\subseteq\overline{M},g \right)$ is an asymptotically hyperbolic space, $\rho$ is a boundary defining function, $\mathcal{E}$ is a smooth complex vector bundle over $\overline{M}$, $\Phi\in\rho C^\infty \left(\overline{M};\End{\mathcal{E}} \right)$, and that $\nabla^\mathcal{E}$ is a smooth connection in $\mathcal{E}$.
Suppose also that $\gamma :  \left(-\infty,\infty \right)\rightarrow M$ is a complete geodesic such that $\liminf_{t\rightarrow\pm\infty}{\rho \left(\gamma \left(t \right) \right)}=0$.
Let $x_0=\lim_{t\rightarrow-\infty}{\gamma \left(t \right)}\in\partial \overline{M}$, which exists by Lemma 2.3 in \cite{Bibitem_15}.

	\begin{enumerate} \item Then for any $e\in\mathcal{E}_{x_0}$, the unique solution $u$ to (\ref{MathItem_2.4}) exists and so does the limit (\ref{MathItem_2.5}).

	\item Furthermore, the maps $e\in\mathcal{E}_{x_0}\mapsto\lim_{t\rightarrow\infty}{u \left(t \right)}$ and $e\in\mathcal{E}_{x_0}\mapsto u \left(t_0 \right)$ for any fixed time $t_0\in \left(-\infty,\infty \right)$ are isomorphisms between fibers of $\mathcal{E}$.\end{enumerate}

\end{lemma}

For future use, we remark that the above lemma and its proof work equally well if one changes ``$t\rightarrow\pm\infty$" to ``$t\rightarrow\mp\infty$" in its statement and in (\ref{MathItem_2.4}) and (\ref{MathItem_2.5}).

Before we state our main result, we establish a way of talking about the decay regularity of the connection $\nabla^\mathcal{E}$:

\begin{definition} \label{MathItem_2.7}  Suppose that $\overline{M}$ is a smooth manifold with smooth boundary, $\rho$ is a boundary defining function, $\mathcal{E}$ is a smooth complex vector bundle over $\overline{M}$, and that $\nabla^\mathcal{E}$ is a smooth connection in $\mathcal{E}$.
If $N\geq0$ is an integer, we say that the ``connection symbols of $\nabla^\mathcal{E}$ are in $\rho^NC^\infty \left(\overline{M} \right)$ in any boundary coordinates (of $\overline{M}$) and frame (for $\mathcal{E}$)" if the following holds.
For any boundary coordinates $ \left(x^i \right)$ of $\overline{M}$ and any frame $ \left(b_k \right)$ of $\mathcal{E}$ over these coordinates' domain, the connection symbols ${}^\mathcal{E}\Gamma_{ij}^k$ in the expression
$$\nabla_v^\mathcal{E}u=v \left(u^k \right)b_k+{}^\mathcal{E}\Gamma_{ij}^kv^iu^jb_k$$
satisfy ${}^\mathcal{E}\Gamma_{ij}^k\in\rho^NC^\infty \left(\overline{M} \right)$, where $v^i$ are the components of $v$ with respect to $ \left(\sfrac{\partial }{\partial  x^i} \right)$ and $u^k$ are the components of $u$ with respect to $b_k$.
\end{definition}

As we mentioned earlier, the answer to our main problem is that we cannot recover the connection and Higgs field from the data (\ref{MathItem_2.5}) because such data can come from two distinct pairs $ \left(\nabla^\mathcal{E},\Phi \right)$ and $ \left({\widetilde{\nabla}}^\mathcal{E},\widetilde{\Phi} \right)$.
However, if that is the case, the two pairs $ \left(\nabla^\mathcal{E},\Phi \right)$ and $ \left({\widetilde{\nabla}}^\mathcal{E},\widetilde{\Phi} \right)$ are related by a well understood gauge relation.
The following is our main result on the matter.
To state it, we use the unit tangent bundle $SM$, the notion of pullback bundles (e.g. $\pi ^\ast\End{\mathcal{E}}$), and connections on endomorphism fields $\nabla^{\End{\mathcal{E}}}$, which we define in Sections \ref{Section_4.2}, \ref{Section_4.5}, and \ref{Section_4.7} respectively.
We also use the regularity spaces $\mathcal{R}^k \left(SM;\ldots \right)$ and the notion of \textit{nontrivial twisted conformal Killing tensor fields} (CKTs for short) which are defined in Section \ref{Section_5.1} and at the end of Section \ref{Section_5.3} respectively.

Note the ``skew-Hermitian" and ``unitary" assumptions in the theorem's statement.

\begin{theorem} \label{MathItem_2.8}  Suppose that $ \left(M\subseteq\overline{M},g \right)$ is a nontrapping asymptotically hyperbolic space, that the sectional curvatures of $g$ are negative, $\rho$ is a boundary defining function, and that $ \left(\mathcal{E},\langle \cdot,\cdot\rangle _\mathcal{E} \right)$ is a smooth complex Hermitian vector bundle over $\overline{M}$.
There exists an integer $N\geq0$ big enough dependent only on $ \left(M,g \right)$ such that the following holds.

Suppose we have a $\Phi\in\rho^{N+1}C^\infty \left(\overline{M};{\End}_{\mathrm{Sk}}{\mathcal{E}} \right)$, a unitary connection $\nabla^\mathcal{E}$ in $\mathcal{E}$ whose connection symbols are in $\rho^NC^\infty \left(\overline{M} \right)$ in any boundary coordinates and frame (in the sense of Definition \ref{MathItem_2.7}), and that we have another pair $\widetilde{\Phi}$ and ${\widetilde{\nabla}}^\mathcal{E}$ satisfying the same conditions.
Consider the unit tangent bundle $SM\subseteq TM$ with projection map $\pi  : SM\rightarrow M$.
Let $A={\widetilde{\nabla}}^\mathcal{E}-\nabla^\mathcal{E}$, consider the connection $\nabla^\prime U:=\nabla^{\End{\mathcal{E}}}U-UA$, and suppose that $\nabla^\prime$ has no nontrivial twisted CKTs in $\mathcal{R}^3 \left(SM;\pi ^\ast\End{\mathcal{E}} \right)$.

Lastly, suppose that the data (\ref{MathItem_2.5}) for all possible $\gamma$ and $e$ as above are the same for (\ref{MathItem_2.4}), and (\ref{MathItem_2.4}) with $\Phi$ and $\nabla^\mathcal{E}$ replaced by $\widetilde{\Phi}$ and ${\widetilde{\nabla}}^\mathcal{E}$ respectively.
Then there exists a unitary\footnote{ An endomorphism section $Q : \overline{M}\rightarrow\End{\mathcal{E}}$ is called unitary if $\left|Qe\right|_\mathcal{E}=\left|e\right|_\mathcal{E}$ for any $e\in\mathcal{E}$.} $Q\in C^0 \left(\overline{M};\End{\mathcal{E}} \right)\cap C^\infty \left(M;\End{\left.\mathcal{E}\right|_M} \right)$ such that $\left.Q\right|_{\partial \overline{M}}=\mathrm{id}$ and over $\overline{M}$ satisfies
\begin{equation} \label{MathItem_2.9} {\widetilde{\nabla}}^\mathcal{E}=Q^{-1}\nabla^\mathcal{E}Q,\ \ \ \ \ \ \ \ \ \ \widetilde{\Phi}=Q^{-1}\Phi Q.\end{equation}
Furthermore, $ \left(Q-\mathrm{id} \right)\in\mathcal{R}^3 \left(SM;\pi ^\ast\End{\mathcal{E}} \right)$ when $Q-\mathrm{id}$ is lifted from $M$ to $SM$ by setting $ \left(Q-\mathrm{id} \right) \left(v \right)= \left(Q-\mathrm{id} \right) \left(x \right)$ for any $v\in S_xM$ using the canonical identification $ \left(\pi ^\ast\End{\mathcal{E}} \right)_v\cong \left(\End{\mathcal{E}} \right)_x$.
\end{theorem}

\begin{remark} \label{MathItem_2.10}  \normalfont The notation ${\widetilde{\nabla}}^\mathcal{E}=Q^{-1}\nabla^\mathcal{E}Q$ means
$${\widetilde{\nabla}}_v^\mathcal{E}u=Q^{-1}\nabla_v^\mathcal{E} \left(Qu \right)$$
for any tangent vector $v\in TM$ and any section $u\in C^\infty \left(\overline{M};\mathcal{E} \right)$.
The equations (\ref{MathItem_2.9}) are called the \textbf{gauge relation} between the pairs $ \left(\nabla^\mathcal{E},\Phi \right)$ and $ \left({\widetilde{\nabla}}^\mathcal{E},\widetilde{\Phi} \right)$, and we will provide intuition below for where it comes from.
\end{remark}

A natural question arises of whether one can explicitly determine or estimate the value of $N$ in the above theorem's statement.
We do not attempt to answer this question or to prove an upper bound.
In this paper we prove that it only depends on the geometry of the space $ \left(M,g \right)$.
Its size comes up when proving the regularity theorem for the transport equation, and in particular we will need it to be big enough so that the decay rate of $\Phi$ and the connection symbols will overpower the growth of the derivatives of the geodesic flow.
Since the latter are difficult to compute, except perhaps in some special cases such as the hyperbolic space, our approach does not indicate exactly how big $N$ needs to be for the proof to work.

We point out that in the special case of when the connection is known and has zero curvature (see Section \ref{Section_4.6} below for the latter), then it is possible to recover the Higgs field.
Furthermore, in this case we can also drop the assumption about the nonexistence of twisted CKTs with suitable regularity because we will get it for free from the zero-curvature assumption.
Here is the precise statement:

\begin{corollary} \label{MathItem_2.11}  Suppose that $ \left(M\subseteq\overline{M},g \right)$ is a nontrapping asymptotically hyperbolic space, that the sectional curvatures of $g$ are negative,\textit{ }$\rho$ is a boundary defining function, and that $ \left(\mathcal{E},\langle \cdot,\cdot\rangle _\mathcal{E} \right)$ is a smooth complex Hermitian vector bundle over $\overline{M}$.
Suppose also that we have a unitary connection $\nabla^\mathcal{E}$ in $\mathcal{E}$ whose curvature is zero.
There exists an integer $N\geq0$ big enough dependent only on $ \left(M,g \right)$ such that the following holds.
Suppose that the connection symbols of $\nabla^\mathcal{E}$ are in $\rho^NC^\infty \left(\overline{M} \right)$ in any boundary coordinates and frame (in the sense of Definition \ref{MathItem_2.7}).
Suppose also that $\Phi,\widetilde{\Phi}\in\rho^{N+1}C^\infty \left(\overline{M};{\End}_{\mathrm{Sk}}{\mathcal{E}} \right)$ are such that the data (\ref{MathItem_2.5}) for all possible $\gamma$ and $e$ as above are the same for (\ref{MathItem_2.4}), and (\ref{MathItem_2.4}) with $\Phi$ replaced by $\widetilde{\Phi}$.
Then $\Phi=\widetilde{\Phi}$.
\end{corollary}

By the explicit formula for the curvature (\ref{MathItem_4.23}), a simple example of when the corollary applies is if $\overline{M}$ is a subset of $\mathbb{R}^{n+1}$, $\mathcal{E}=\overline{M}\times\mathbb{C}^d$ is the trivial bundle whose sections we write as column vectors, and the connection $\nabla^\mathcal{E}$ is simply given by $\nabla_v^\mathcal{E}u=\left[v \left(u^1 \right),\ldots,v \left(u^d \right)\right]$.

We end this section with discussion on intuition and outline of the proofs.
The relation (\ref{MathItem_2.9}) may look mysterious at first, so let us give intuition for it.
The proof of Theorem \ref{MathItem_2.8} essentially begins with taking (\ref{MathItem_2.4}) and rewriting it in terms of an endomorphism field $U : \mathbb{R}\rightarrow\End{\mathcal{E}}$ in a way that it encodes the same data:
\begin{equation} \label{MathItem_2.12} \nabla_{\dot{\gamma} \left(t \right)}^{\End{\mathcal{E}}}U \left(t \right)+\Phi \left(\gamma \left(t \right) \right)U \left(t \right)=0,\ \ \ \ \lim_{t\rightarrow-\infty}{U \left(t \right)}=\mathrm{id},\end{equation}
where ``$\nabla^{\End{\mathcal{E}}}$" is the natural connection on the space of endomorphism fields $\End{\mathcal{E}}$ induced by $\nabla^\mathcal{E}$ (see Section \ref{Section_4.7}).
Precisely, we read of the data (\ref{MathItem_2.5}) from (\ref{MathItem_2.12}) as follows: for any $e\in\mathcal{E}_{x_0}$ where $x_0\in\partial  M$ we solve $\nabla_{\dot{\gamma} \left(t \right)}^\mathcal{E}w \left(t \right)=0$ with $\lim_{t\rightarrow-\infty}{w\circ\gamma \left(t \right)}=e$ and show that the data (\ref{MathItem_2.5}) is given by $\lim_{t\rightarrow\infty}{U \left(t \right)w \left(t \right)}$. 

Then we ask the following question.
Having our unitary connection $\nabla^\mathcal{E}$, skew-Hermitian $\Phi$, and the data (\ref{MathItem_2.5}) that they generate, how can we ``come up" with another unitary ${\widetilde{\nabla}}^\mathcal{E}$ and skew-Hermitian $\widetilde{\Phi}$ that generates the same data set? To do this, we take an arbitrary unitary endomorphism field $Q$ and manipulate the above equation as follows:
\begin{align*}\nabla_{\dot{\gamma}}^{\End{\mathcal{E}}} \left(QQ^{-1}U \right)+\Phi \left(QQ^{-1}U \right)=0,& &\mathrm{since}\ QQ^{-1}=\mathrm{id} ,\\\nabla_{\dot{\gamma}}^{\End{\mathcal{E}}} \left(Q \right)Q^{-1}U+Q\nabla_{\dot{\gamma}}^{\End{\mathcal{E}}} \left(Q^{-1}U \right)+\Phi \left(QQ^{-1}U \right)=0,& &\mathrm{product\ rule},\\Q^{-1}\nabla_{\dot{\gamma}}^{\End{\mathcal{E}}} \left(Q \right)\widetilde{U}+\nabla_{\dot{\gamma}}^{\End{\mathcal{E}}}\widetilde{U}+ \left(Q^{-1}\Phi Q \right)\widetilde{U}=0,& &\mathrm{multiply\ by}\ Q^{-1}\ \mathrm{and\ set} \ \widetilde{U}=Q^{-1}U,\\\nabla_{\dot{\gamma}}^{\End{\mathcal{E}}}\widetilde{U}+A\widetilde{U}+\widetilde{\Phi}\widetilde{U}=0,& &\ \ \ \end{align*}
where in the last step we set $A=Q^{-1}\nabla_{\dot{\gamma}}^{\End{\mathcal{E}}} \left(Q \right)$ and $\widetilde{\Phi}:=Q^{-1}\Phi Q$.
As we will show, from the way we read off the data (\ref{MathItem_2.5}) from (\ref{MathItem_2.12}) it follows that ${\widetilde{\nabla}}^\mathcal{E}=\nabla^\mathcal{E}+A$ and $\widetilde{\Phi}$ generate the same data set (\ref{MathItem_2.5}) if $\left.Q\right|_{\partial \overline{M}}=\mathrm{id}$ and $Q$ satisfies (\ref{MathItem_2.9}).
The main point of Theorem \ref{MathItem_2.8} is that this is the only way that we can produce another unitary ${\widetilde{\nabla}}^\mathcal{E}$ and skew-Hermitian $\widetilde{\Phi}$ that generates the same data set.
This example is illustrative in the sense that, in our objective to prove the existence of $Q$ satisfying the desired gauge relation (\ref{MathItem_2.9}), the above tells us that we should try
\begin{equation} \label{MathItem_2.13} Q=U{\widetilde{U}}^{-1}\end{equation}
where $U$ and $\widetilde{U}$ are defined as above.

To elaborate more on the outline of the proof, we will show that $Q-\mathrm{id}$ satisfies a transport equation of the specific form
\begin{equation} \label{MathItem_2.14} \nabla_X^{\pi ^\ast\End{\mathcal{E}}}W+\Psi W=f\end{equation}
on the \textit{sphere bundle} $SM$ where $X$ is the geodesic vector field (c.f. (\ref{MathItem_5.86}) below).
The right-hand side $f$ will be sufficiently regular at infinity and have Fourier modes of order no bigger than one with respect to the vertical Laplacian.
In Section \ref{Section_5.4} we will prove a regularity theorem for transport equations that will imply that the solution $W$ is also sufficiently regular at infinity.
In Section \ref{Section_5.3} we will conduct a Fourier study of transport equations that, combined with the just mentioned result, will imply that $W$, and hence $Q$, have Fourier degree zero (i.e. are of the form $C^\infty \left(\overline{M};\End{\mathcal{E}} \right)$).
This, in particular, is the step where we use the nonexistence of nontrivial twisted conformal Killing tensor fields (CKTs) in $\mathcal{R}^3 \left(SM;\pi ^\ast\End{\mathcal{E}} \right)$.
From there it will quickly follow that $Q$ satisfies the conclusions in Theorem \ref{MathItem_2.8}.

One of the key steps in our investigation will be to generalize an identity called the ``Pestov identity" to vector bundles over asymptotically hyperbolic (AH) spaces (i.e. Theorem \ref{MathItem_5.9} below).
This is similar to the Pestov identity for scalar functions that was generalized to nontrapping AH spaces in \cite{Bibitem_15} (see bottom of page 2892 there) which the authors accomplished by also generalizing Santal\'o's formula to nontrapping AH spaces.
We will take a different approach which will allow us to not assume that the manifold nontrapping to prove this intermediate step.
Of course, this will be more than we need since in our final result Theorem \ref{MathItem_2.8} we do assume that $ \left(M,g \right)$ is nontrapping.

Another key feature of this paper is the formulation of regularity spaces to which solutions of transport equations over $SM$ of the form (\ref{MathItem_2.14}) belong to and to which we can apply vertical differential operators while maintaining sufficient regularity at infinity.
Precisely, we have chosen the regularity spaces $\mathcal{R}^k$ that we introduce in Section \ref{Section_5.1} because they are subsets of $L^2$ and the differential operators that we use in this paper have the mapping properties $\mathcal{R}^k\rightarrow\mathcal{R}^{k-1}$.
Both properties play a crucial role in our generalization of the Pestov identity to asymptotically hyperbolic spaces (Theorem \ref{MathItem_5.9}).
The approach that we take to study the existence and regularity of solutions to transport equations (see Proposition \ref{MathItem_5.56}) is to embed $SM$ as a subset of the $b$-cosphere bundle ${}^bS^\ast\overline{M}$ and the 0-cosphere bundle ${}^0S^\ast\overline{M}$ (defined in Section \ref{Section_4.2}).
The domain ${}^bS^\ast\overline{M}$ provides a natural setting to prove the existence of solutions to transport equations with boundary conditions at infinity due to the nice behavior of the geodesic vector field $X$ at the boundary $\partial {}^bS^\ast\overline{M}$.
The bundle ${}^0S^\ast\overline{M}$ on the other hand provides a convenient domain to prove the regularity of solutions at infinity due to its compactness.
\subsection{Existence of Examples} \label{Section_2.4} 

Here we address the question of the existence of manifolds $ \left(M,g \right)$ and connections $\nabla^\mathcal{E}$ that satisfy the assumptions of our main results Theorem \ref{MathItem_2.8} and Corollary \ref{MathItem_2.11}.
The existence of a nontrapping asymptotically hyperbolic (AH) space $ \left(M,g \right)$ with negative sectional curvature is provided by the Poincar\'e ball model.
More examples can be produced by (smoothly) deforming the metric slightly in small enough regions; we provide a brief argument for this.

Suppose that we have an AH space $ \left(M\subseteq\overline{M},g \right)$ that is nontrapping and has negative sectional curvature.
Let $\rho$ be a geodesic boundary defining function as defined in Section \ref{Section_2.2}.
It is clear that if one takes any open $U\subseteq M$ whose closure is compact and contained in the domain of local coordinates, then small enough deformations of the metric $g$ over $U$ with respect to the coordinates' $C^2$-norm will preserve its negative curvature.
So, let us demonstrate that if this region and deformation are small enough, then we can also preserve the nontrapping property of $ \left(M,g \right)$.
It follows from Lemma 2.3 in \cite{Bibitem_15} that there exists an $\varepsilon_0>0$ small enough such that if a $g$-geodesic $\gamma$ makes its way into the region $\left\{\rho<\varepsilon_0\right\}$ at some time $t_0$, then it will stay in $\left\{\rho<\varepsilon_0\right\}$ for all $t>t_0$ and will escape to infinity (i.e. $\liminf_{t\rightarrow\infty}{\rho \left(\gamma \left(t \right) \right)}=0$).
Now, let us suppose that our $U$ above is outside of this region: $U\subseteq\left\{\rho\geq\varepsilon_0\right\}$.
Let $\theta_g$ be the flow of the geodesic vector field $X_g$ over $TM$ (we wrote ``$g$" here for emphasis).
The continuous dependence of solutions to ordinary differential equations on parameters (e.g. see Theorem 7.4 of Chapter 1 in \cite{Bibitem_6}) implies that the map $\theta_g$ depends continuously on the values of $g$ and its first and second order partials in $U$.
By assumption $g$ is nontrapping, and hence all of its geodesics coming out of $U$ will eventually make their way to $\left\{\rho<\varepsilon_0\right\}$ and escape to infinity.
If we make $U$ precompact in coordinates, it is not hard to see then that small enough deformations of the metric $g$ over $U$ with respect to the $C^2$-norm will also satisfy that all of their geodesics coming out of $U$ will eventually make their way to $\left\{\rho<\varepsilon_0\right\}$, and hence such deformations will be nontrapping.
This proves our claim.

Finally, the following result addresses the question of the existence of connections $\nabla^\mathcal{E}$ on such AH spaces that satisfy the assumption that the connection $\nabla^\prime$ in Theorem \ref{MathItem_2.8} has no nontrivial twisted CKTs in $\mathcal{R}^3 \left(SM;\pi ^\ast\End{\mathcal{E}} \right)$.
To make the notation simpler, we prove instead a result for when $\nabla^\mathcal{E}$ has no nontrivial twisted CKTs in $\mathcal{R}^3 \left(SM;\pi ^\ast\mathcal{E} \right)$, from which the former will follow by substituting $\nabla^\prime$ into $\nabla^\mathcal{E}$.
Similar results can be found in [41, Corollary 3.6] and [8, Theorem 1.6].
We refer the reader to Section \ref{Section_4.6} for the definition of the curvature operator $F^\mathcal{E}$ of $\nabla^\mathcal{E}$.

\begin{theorem} \label{MathItem_2.15}  Suppose that $ \left(M\subseteq\overline{M},g \right)$ is a nontrapping asymptotically hyperbolic space, that the sectional curvatures of $g$ are negative, and that $ \left(\mathcal{E},\langle \cdot,\cdot\rangle _\mathcal{E} \right)$ is a smooth complex Hermitian vector bundle over $\overline{M}$.
Then the sectional curvatures are bounded above by $-\kappa$ for some $\kappa>0$.
For any unitary connection $\nabla^\mathcal{E}$ in $\mathcal{E}$ whose curvature norm satisfies $\left\lVert F^\mathcal{E}\right\rVert _{L^\infty}\le\kappa\sqrt n$, there are no nontrivial twisted CKTs in $\mathcal{R}^3 \left(SM;\pi ^\ast\mathcal{E} \right)$.
\end{theorem}

Similar to the remark made after Corollary \ref{MathItem_2.11}, an example of when this lemma applies is if $\overline{M}$ is a subset of $\mathbb{R}^{n+1}$, $\mathcal{E}=\overline{M}\times\mathbb{C}^d$ is the trivial bundle whose sections we write as column vectors, and the connection $\nabla^\mathcal{E}$ is given by $\nabla_v^\mathcal{E}u=\left[v \left(u^1 \right)+{}^\mathcal{E}\Gamma_{ij}^1v^iu^j,\ldots,v \left(u^d \right)+{}^\mathcal{E}\Gamma_{ij}^dv^iu^j\right]$ where the connection symbols ${}^\mathcal{E}\Gamma_{ij}^k$ decay fast enough at the boundary (i.e. ${}^\mathcal{E}\Gamma_{ij}^k\in\rho^NC^\infty \left(\overline{M} \right)$ for big enough $N\geq0$), and such that they and their first partials are small enough so that the curvature estimate in the above lemma is satisfied.
\subsection{Non-Abelian X-Ray Transform} \label{Section_2.5} 

We mention a way to formalize the operator that takes $ \left(\nabla,\Phi \right)$ to the map taking pairs $ \left(\gamma,e \right)$ as above to (\ref{MathItem_2.5}) using the $b$ and 0 cotangent bundles that we introduce in Sections \ref{Section_4.1} and \ref{Section_4.2} below.
We will not make use of this formulation, and only the material up to (\ref{MathItem_2.16}) and the two sentences afterwards here will be used later in the paper.
Suppose that $ \left(M\subseteq\overline{M},g \right)$ is a nontrapping asymptotically hyperbolic (AH) space, $\mathcal{E}$ is a smooth complex vector bundle over $\overline{M}$, $\Phi\in\rho C^\infty \left(\overline{M};\End{\mathcal{E}} \right)$, and that $\nabla^\mathcal{E}$ is a smooth connection in $\mathcal{E}$.

Consider the cotangent and $b$-cotangent bundles $T^\ast M$ and ${}^bT^\ast\overline{M}$ respectively, and their unit cosphere bundles $S^\ast M$ and ${}^bS^\ast\overline{M}$ respectively.
Suppose $ \left(\rho,y^1,\ldots,y^n \right)$ are boundary coordinates of $\overline{M}$ and consider the frame $ \left(\sfrac{d\rho}{\rho},dy^1,\ldots,dy^n \right)$ spanning covectors in $T^\ast M$.
On page 2865 of \cite{Bibitem_15} the authors remind the reader that this extends to the boundary to become a smooth frame of  ${}^bT^\ast\overline{M}$ and that furthermore if $\zeta=\eta_0\sfrac{d\rho}{\rho}+\eta_\mu dy^\mu\in\left.{}^bT^\ast\overline{M}\right|_{\partial \overline{M}}$ is over the boundary, then the map
$$\eta_0\frac{d\rho}{\rho}+\eta_\mu dy^\mu\longmapsto\eta_0$$
is well defined (i.e. independent of the coordinates $ \left(\rho,y^\mu \right)$ that we choose).
The boundary of the unit cosphere bundle ${}^bS^\ast\overline{M}\subseteq{}^bT^\ast\overline{M}$ turns out to have the following two components:
\begin{align*}\partial _-{}^bS^\ast\overline{M}=\left\{\zeta\in\left.{}^bT^\ast\overline{M}\right|_{\partial \overline{M}} : \eta_0=1\right\}& &\mathrm{called\ the\ ``incoming\ boundary,"}\\\partial _+{}^bS^\ast\overline{M}=\left\{\zeta\in\left.{}^bT^\ast\overline{M}\right|_{\partial \overline{M}} : \eta_0=-1\right\}& &\mathrm{called\ the\ ``outgoing\ boundary."}\end{align*}
Let $\pi  : SM\rightarrow M$ and $\pi _b : {}^bS^\ast\overline{M}\rightarrow\overline{M}$ denote the natural projection maps.
Recall that any unit-speed geodesic $\gamma :  \left(-\infty,\infty \right)\rightarrow M$ is the image under $\pi $ of an integral curve $\sigma :  \left(-\infty,\infty \right)\rightarrow SM$ of the geodesic vector field $X$.
Letting $X_b$ denote the pushforward of $X$ onto $\left.{}^bS^\ast\overline{M}\right|_M$ via the canonical identification between $TM$ and $\left.{}^bT^\ast\overline{M}\right|_M$, we have that $\gamma$ is the image under $\pi _b$ of an integral curve $\sigma_b :  \left(-\infty,\infty \right)\rightarrow\left.{}^bS^\ast\overline{M}\right|_M$ of $X_b$.
It follows from the proof of Corollary 2.5 in \cite{Bibitem_15} that the limit of any such curve exists in ${}^bS^\ast\overline{M}$:
\begin{equation} \label{MathItem_2.16} \begin{array}{r@{}l}\lim_{t\rightarrow-\infty}{\sigma_b}=\partial _-{}^bS^\ast\overline{M},\\\lim_{t\rightarrow\infty}{\sigma_b}=\partial _+{}^bS^\ast\overline{M}.\\\end{array}\end{equation}
Intuitively, the first limit here can be thought of as the ``initial velocity" of the geodesic as it ``enters" the AH space at infinity, while the second its ``exit velocity" as it ``leaves" at infinity.
Conversely, it follows from the same proof that every $\zeta\in\partial _-{}^bS^\ast\overline{M}$ (resp.
$\zeta\in\partial _+{}^bS^\ast\overline{M}$) is the limit in ${}^bS^\ast\overline{M}$ of a unique (up to reparameterization) such curve $\sigma_b$ as $t\rightarrow-\infty$ (resp.
$t\rightarrow\infty$).

Hence we may define the map
$$T^{ \left(\nabla^\mathcal{E},\Phi \right)} : \left.\pi _b^\ast\mathcal{E}\right|_{\partial _-{}^bS^\ast\overline{M}}\longrightarrow\left.\pi _b^\ast\mathcal{E}\right|_{\partial _+{}^bS^\ast\overline{M}}$$
as follows.
Take any $e\in\left.\pi _b^\ast\mathcal{E}\right|_{\partial _-{}^bS^\ast\overline{M}}$ whose base point we denote by $\zeta\in\partial _-{}^bS^\ast\overline{M}$.
Let $\sigma_b$ be an integral curve of $X_b$ with $\zeta=\lim_{t\rightarrow-\infty}{\sigma_b}$ and let $\zeta_{\mathrm{exit}}=\lim_{t\rightarrow\infty}{\sigma_b}$.
Take the geodesic $\gamma=\pi _b\circ\sigma_b$ and let $u$ be the solution to (\ref{MathItem_2.4}) where we let $e$ also denote the element in $\mathcal{E}_{\pi _b \left(\zeta \right)}$ that is canonically identified to $e\in \left(\pi _b^\ast\mathcal{E} \right)_\zeta$.
Then we set
$$T^{ \left(\nabla^\mathcal{E},\Phi \right)} \left(e \right)=\lim_{t\rightarrow\infty}{u \left(t \right)},$$
making the similar canonical identification $ \left(\pi _b^\ast\mathcal{E} \right)_{\zeta_{\mathrm{exit}}}\cong\mathcal{E}_{\pi _b \left(\zeta_{\mathrm{exit}} \right)}$.
We point out that this limit exists by Lemma \ref{MathItem_2.6} part 1), and that the ``$T$" here stands for ``transport equation."

\begin{definition} \label{MathItem_2.17}  Suppose that $ \left(M\subseteq\overline{M},g \right)$ is an asymptotically hyperbolic space, $\mathcal{E}$ is a smooth complex vector bundle over $\overline{M}$, $\Phi\in\rho C^\infty \left(\overline{M};\End{\mathcal{E}} \right)$, and that $\nabla^\mathcal{E}$ is a smooth connection in $\mathcal{E}$.
The operator
$$ \left(\nabla^\mathcal{E},\Phi \right)\longmapsto T^{ \left(\nabla^\mathcal{E},\Phi \right)}$$
is called the \textbf{non-abelian X-ray transform}.
This is well defined by the discussion above.
\end{definition}

For instance, another way to formulate Corollary \ref{MathItem_2.11} above is that for any $g$ and $\nabla^\mathcal{E}$ satisfying the conditions there, the non-abelian X-ray transform is injective over the set of all skew-Hermitian Higgs field satisfying the decay condition also described there.
\subsection{Prior Research Discussion and Applications} \label{Section_2.6} 

A standard approach for studying injectivity properties of X-ray transforms is via energy identities that was first introduced in \cite{Bibitem_35}.
The type of energy estimate that is used in this approach is called the \textbf{Pestov identity} (or \textbf{Muhometov-Pestov identity}) which over the years has taken many forms as authors apply them in various contexts - see for instance \cite{Bibitem_12}, \cite{Bibitem_40}, \cite{Bibitem_43}, and \cite{Bibitem_44}.
The mentioned paper \cite{Bibitem_12} furthermore explains the connection between X-ray transform over connections and inverse problems related to the wave equation.
Of recent works, in dimension two the authors of \cite{Bibitem_39} used a Pestov identity to prove solenoidal injectivity of the X-ray transform over tensors, and in their earlier work \cite{Bibitem_38} they proved a ``Pestov type identity" to study the attenuated ray transform with a connection and Higgs field.

The paper \cite{Bibitem_41} proceeded to generalize these methods to manifolds of dimensions greater than two, but it did not cover the case of connections.
In \cite{Bibitem_18} the authors generalized the setup in \cite{Bibitem_41} where they studied the X-ray transform for connections and Higgs fields together.
For instance, Theorem \ref{MathItem_2.8} above was proved in \cite{Bibitem_18} in the case when $ \left(M,g \right)$ is a compact Riemannian manifold, has strictly convex boundary, has negative sectional curvature, the boundary condition in (\ref{MathItem_2.4}) is changed to $u \left(\gamma \left(a \right) \right)$ where $\gamma : \left[a,b\right]\rightarrow M$ is a unit-speed geodesic traveling between boundary points, and (\ref{MathItem_2.5}) is changed to ``recording" $u \left(\gamma \left(b \right) \right)$.
In this paper we also generalize the Pestov identity proved in \cite{Bibitem_18} to asymptotically hyperbolic (AH) spaces, of which a similar formula also appears in \cite{Bibitem_45}.

We prove the nonexistence of nontrivial twisted conformal tensor fields (CKTs) in Theorem \ref{MathItem_2.15} under the condition that the curvature of the connection is small enough.
Our approach is based on [41, Corollary 3.6] that proved a similar result in the compact setting, which was also proved in [8, Theorem 1.6].
Additionally, the work \cite{Bibitem_18} proved that a twisted CKT that vanishes on $\pi ^{-1}\left[\Gamma\right]$ where $\Gamma$ is a hypersurface in $M$ and $\pi  : SM\rightarrow M$ is the projection from the unit-sphere bundle must be trivial.
By setting $\Gamma=\partial  M$, the authors of \cite{Bibitem_18} were able to state their analog of our Theorem \ref{MathItem_2.8} without any assumption of the nonexistence of nontrivial twisted CKTs.
The work \cite{Bibitem_18} also proved the nonexistence of nontrivial twisted CKTs on closed surfaces and for sufficiently high Fourier modes on closed manifolds under certain regularity conditions.
The nonexistence of nontrivial twisted CKTs on closed manifolds was further studied in \cite{Bibitem_5}.

We mention the early work \cite{Bibitem_1} that studied the injectivity of the X-ray transform for one-forms.
The work of \cite{Bibitem_47} studied injectivity on tensor fields of rank $m\le2$ for analytic simple metrics and a generic class of two-dimensional simple metrics, and proved a stability estimate for the normal operator.
Later, \cite{Bibitem_46} proved injectivity on two-tensors for all two-dimensional simple metrics which was then extended to tensors of all rank in \cite{Bibitem_39}.
The papers \cite{Bibitem_52} and \cite{Bibitem_50} proved injectivity for functions and two-tensors respectively on Riemannian manifolds that admit convex foliations.
The paper \cite{Bibitem_17} proved injectivity on tensors of all ranks over Riemannian manifolds with negative curvature and strictly convex boundary.
We mention that the work \cite{Bibitem_4} characterized the range of the non-abelian X-ray transform on simple surfaces in terms of boundary quantities and that \cite{Bibitem_3} and \cite{Bibitem_33} proved stability estimates for it over Higgs fields.
Microlocal techniques have also been applied to the study of the X-ray transform in the presence of conjugate points - we refer the reader to the works \cite{Bibitem_23}, \cite{Bibitem_34}, \cite{Bibitem_48}, and \cite{Bibitem_49}.

In the noncompact realm, injectivity for the scalar X-ray transform over hyperbolic spaces was proved in \cite{Bibitem_21}, and inversion formulas are given in \cite{Bibitem_2} and \cite{Bibitem_20}.
In \cite{Bibitem_27} the author proved analogous injectivity over Cartan-Hadamard manifolds and in \cite{Bibitem_28} the results were extended to higher dimensions and tensor fields.
The paper \cite{Bibitem_37} proved a gauge equivalence for the X-ray transform for connections on Euclidean space assuming a bound on the size of the connection in dimension two.

AH manifolds have gained interest in the past two decades partly due to their role in physics such as the AdS/CFT conjecture made in \cite{Bibitem_29}.
The work \cite{Bibitem_7} for instance described the role of integral geometry in the AdS/CFT correspondence.
In this setting, the paper \cite{Bibitem_15} proved injectivity of the X-ray transform for tensor of all orders on asymptotically hyperbolic spaces.
On simple AH manifolds, the work \cite{Bibitem_10} generalized their result for the scalar X-ray transform by proving a stability estimate for the normal operator.
Analogous to the local problem studied in \cite{Bibitem_52}, \cite{Bibitem_11} proved a local injectivity result for the scalar X-ray transform on AH spaces.

Regarding applications of the non-abelian X-ray transform, we also mention its appearance in the theory of solitons when studying the Bogomolny equations in dimensions $2+1$ - see \cite{Bibitem_30} and \cite{Bibitem_53} for details.
The paper \cite{Bibitem_24} describes its applications to coherent quantum tomography.
For a survey of the non-abelian X-ray transform and to read more about its applications, we refer the reader to \cite{Bibitem_36}.
\subsection{Acknowledgments} \label{Section_2.7} 

First of all, I would like to thank my advisor Gunther Uhlmann who from the beginning of my PhD studies has always supported and encouraged my curiosities, mathematical development, and growth as a human being and mathematician.
I am grateful for both his patience and the invaluable advice that he gave me during times of both progress and confusion, including identifying key references for inspiration on how to proceed forward whenever faced with difficulties.
I would also like to thank Robin Graham for helping me navigate my graduate career and having many fascinating mathematical discussions with me over the years.
I would like to thank Kelvin Lam for our regular meetings in which we discussed our work.

I am grateful to everyone in the math department at the University of Washington for creating a welcoming and productive environment in which I could study.
The author was partly supported by the NSF grant \# 2105956.
\section{Well-Definedness of the Non-abelian X-Ray Transform} \label{Section_3} 

In this section we prove Lemma \ref{MathItem_2.6}.
We start with part 1).
Take any $e\in\mathcal{E}_{x_0}$.
The plan is to do the following three tasks:

	\begin{enumerate} \item prove the existence and uniqueness of the solution to the initial value problem (\ref{MathItem_2.4}), on an interval of the form $ \left(-\infty,t_0\right]$ for some $t_0\in\mathbb{R}$,

	\item argue the existence and uniqueness on the rest of the interval $\left[t_0,\infty \right)$ (and hence everywhere),

	\item and finally prove that the limit (\ref{MathItem_2.5}) exists.\end{enumerate}

We begin with task 1), which we prove by mapping the infinite interval to a bounded one and then applying standard existence and uniqueness results of ordinary differential equations (ODEs).
Let $d=\rank{\mathcal{E}}$.
Let $ \left(\rho,y^1,\ldots,y^n \right)= \left(x^i \right)$ be asymptotic boundary normal coordinates of $\overline{M}$ containing $x_0$ in their domain and let $ \left(b_j \right)$ be a frame for $\mathcal{E}$ over the same domain.
Let ${}^\mathcal{E}\Gamma_{ij}^k$ denote the connection symbols of $\nabla^\mathcal{E}$ with respect to $ \left(\sfrac{\partial }{\partial  x^i} \right)$ and $ \left(b_k \right)$.
Let $t_0\in\mathbb{R}$ be a time such that the image of $\gamma$ is contained in these coordinates for all times $t\in \left(-\infty,t_0\right]$.
Then, writing $u=u^kb_k$, in these coordinates for $t\in \left(-\infty,t_0\right]$ we have that (\ref{MathItem_2.4}) becomes the following system of ODEs
\begin{equation} \label{MathItem_3.1} \frac{du^k}{dt}+{}^\mathcal{E}\Gamma_{ij}^k{\dot{\gamma}}^iu^j+\Phi_i^ku^i=0,\ \ \ \ \ \ \ \ \ \ \lim_{t\rightarrow-\infty}{u^k}=e^k,\ \ \ \ \ \ \ \ \ \ k\in\left\{1,\ldots,d\right\},\end{equation}
Let us look at the growth rate of the ${\dot{\gamma}}^i$'s.
By definition, $g=\sfrac{\overline{g}}{\rho^2}$ for some smooth metric $\overline{g}$ on $\overline{M}$.
Since $\gamma$ has a constant speed one, we have that
$${\overline{g}}_{ij}{\dot{\gamma}}^i{\dot{\gamma}}^j=\rho^2.$$
Clearly the closure of the image $\gamma \left(-\infty,t_0\right]$ is a compact subset of our coordinates' domain, and so the matrix in the bilinear form $v\mapsto{\overline{g}}_{ij}v^iv^j$ has a minimum positive eigenvalue along this set.
Hence from the above we get that there exists a $C>0$ such that each $\left|{\dot{\gamma}}^i\right|\le C\rho$.

Now, take the diffeomorphism $h :  \left(-\sfrac{\pi }{2},s_0\right]\rightarrow \left(-\infty,t_0\right]$ given by $h \left(s \right)=\tan{s}$.
Making the change of variables $t=h \left(s \right)$ in (\ref{MathItem_3.1}) gives that for each $k\in\left\{1,\ldots,d\right\}$
\begin{equation} \label{MathItem_3.2} \frac{du^k}{ds}+{}^\mathcal{E}\Gamma_{ij}^k{\dot{\gamma}}^i\frac{dh}{ds}u^j+\Phi_i^k\frac{dh}{ds}u^i=0\ \ \ \mathrm{on}  \ s\in \left(-\sfrac{\pi }{2},s_0 \right),\ \ \ \ \ u^k \left(-\sfrac{\pi }{2} \right)=e^k.\end{equation}
In other words, the existence and uniqueness of a continuous solution $u$ to this system of initial value problems will prove task 1).
This in turn will follow from standard results on linear ODEs (see for instance \cite{Bibitem_6}) if we show that the above coefficients ${}^\mathcal{E}\Gamma_{ij}^k{\dot{\gamma}}^ih^\prime$ and $\Phi_i^kh^\prime$ extend continuously to $s=-\sfrac{\pi }{2}$.

It follows from Lemma 2.3 in \cite{Bibitem_15} (specifically (2.11) there) that there exists a constant $C^\prime>0$ such that for $t\in \left(-\infty,t_0\right]$,
\begin{equation} \label{MathItem_3.3} \rho\circ\gamma \left(t \right)<C^\prime e^t.\end{equation}
Since by assumption ${}^\mathcal{E}\Gamma_{ij}^k\in C^\infty \left(\overline{M} \right)$, $\left|{\dot{\gamma}}^i\right|\le C\rho$, and $\Phi_i^k\in\rho C^\infty \left(\overline{M} \right)$, we have that there exists a constant $C^{\prime\prime}>0$ such that for $s\in \left(-\sfrac{\pi }{2},s_0\right]$ both $\left|{}^\mathcal{E}\Gamma_{ij}^k{\dot{\gamma}}^ih^\prime\right|$ and $\left|\Phi_i^kh^\prime\right|$ are bounded above by
$$C^{\prime\prime} \left(\rho\circ\gamma \left(h \left(s \right) \right) \right)h^\prime \left(s \right)\le C^{\prime\prime}C^\prime e^{\tan{ \left(s \right)}}\sec^2{ \left(s \right)}\rightarrow0\ \ \ \ \ \mathrm{as}  \ s\rightarrow-{\sfrac{\pi }{2}}^+$$
Hence indeed ${}^\mathcal{E}\Gamma_{ij}^k{\dot{\gamma}}^ih^\prime$ and $\Phi_i^kh^\prime$ extend continuously to $s=-\sfrac{\pi }{2}$.

Task 2) follows by applying standard existence and uniqueness theory of ODEs in coordinates and frames for $\mathcal{E}$ as one travels along the geodesic.
Task 3) is proved similarly to 1) except one uses Lemma 2.3 in \cite{Bibitem_15} in forward time (for instance, the $e^t$ in (\ref{MathItem_3.3}) will change to $e^{-t}$).

To prove part 2) of Lemma \ref{MathItem_2.6}, we recall the fact from ODE theory that for linear homogeneous systems of the form
$$\frac{dw}{dt} \left(t \right)+A \left(t \right)w \left(t \right)=0,\ \ \ \ \ w \left(t_0 \right)=e_0,$$
were $w$ is a column vector and $A$ is a continuous matrix, the map $e_0\mapsto w \left(t_1 \right)$ is an isomorphism for any fixed $t_1$ (e.g. see Chapter 3 Section \ref{Section_2} of \cite{Bibitem_6} - in particular Theorem 2.2).
Part 2) of Lemma \ref{MathItem_2.6} then follows by applying this result to (\ref{MathItem_3.2}).

\begin{flushright}$\blacksquare$\end{flushright}
\section{Geometric Preliminaries} \label{Section_4} 

Throughout this section we assume that $ \left(M\subseteq\overline{M},g \right)$ is an asymptotically hyperbolic space, $\rho$ is a boundary defining function, $ \left(\mathcal{E},\langle \cdot,\cdot\rangle _\mathcal{E} \right)$ is a smooth complex Hermitian vector bundle over $\overline{M}$, and that $\nabla^\mathcal{E}$ is a smooth connection in $\mathcal{E}$.

\subsection{The b and 0 Cotangent Bundles} \label{Section_4.1} 

In this section we introduce the $b$ and 0 cotangent bundles.
We will only state their properties, referring the reader to Section 2.2 in \cite{Bibitem_32} for more details.
We begin by recalling that lowering and raising an index with respect to $g$ provides a bundle isomorphism between the tangent and cotangent bundles over the interior:
$$\flat  : TM\longrightarrow T^\ast M,\ \ \ \ \ \ \ \ \ \ \sharp  : T^\ast M\longrightarrow TM.$$
The first bundle that we introduce is the {\boldmath $b$}\textbf{-tangent bundle} ``${}^bT\overline{M}$," which comes with a canonical smooth map $F : {}^bT\overline{M}\rightarrow T\overline{M}$ that has the following two properties:

	\begin{enumerate} \item $F$ induces a bijection between smooth sections of ${}^bT\overline{M}$ and smooth sections of $T\overline{M}$ that are tangent to the boundary $\partial \overline{M}$.

	\item For any fixed point $x\in\overline{M}$, $F$ restricts to a linear homomorphism $F_x : {}^bT_x\overline{M}\rightarrow T_x\overline{M}$ that is also an isomorphism when $x$ is in the interior $M$.\end{enumerate}

The second is the \textbf{0-tangent bundle} ``${}^0T\overline{M}$," which is defined similarly as coming with a smooth map $H : {}^0T\overline{M}\rightarrow T\overline{M}$ that has the following two properties:

	\begin{enumerate} \item $H$ induces a bijection between smooth sections of ${}^0T\overline{M}$ and smooth sections of $T\overline{M}$ that vanish at the boundary $\partial \overline{M}$.

	\item For any fixed point $x\in\overline{M}$, $H$ restricts to a linear homomorphism $H_x : {}^0T_x\overline{M}\rightarrow T_x\overline{M}$ that is also an isomorphism when $x$ is in the interior $M$.\end{enumerate}

Of more importance to us will be the dual bundles ${}^bT^\ast\overline{M}$ and ${}^0T^\ast\overline{M}$, which are called the {\boldmath $b$}\textbf{ and 0 cotangent bundles} respectively.
They naturally generate pullback maps $F^\ast : T^\ast\overline{M}\rightarrow{}^bT^\ast\overline{M}$ and $H^\ast : T^\ast\overline{M}\rightarrow{}^0T^\ast\overline{M}$ which are also bundle homomorphisms that are isomorphisms on fibers over the interior $M$ (c.f. points 2) above)

\begin{remark} \label{MathItem_4.1}  \normalfont Considering that $\flat ,\sharp ,F,H,F^\ast,H^\ast$ are all isomorphisms (on fibers) over the interior $M$, we will often identify two points in $TM$, $T^\ast M$, $\left.{}^bT\overline{M}\right|_M$, $\left.{}^bT^\ast\overline{M}\right|_M$, $\left.{}^0T\overline{M}\right|_M$, and $\left.{}^0T^\ast\overline{M}\right|_M$ as being the same if it is possible to go from one to the other by a composition of the ``canonical identification" maps mentioned above.
\end{remark}

We mention important frames for the $b$ and 0 cotangent bundles near the boundary $\partial \overline{M}$.
Suppose $ \left(\rho,y^1,\ldots,y^n \right)= \left(x^i \right)$ are boundary coordinates of $\overline{M}$.
Then it turns out that
$$F^\ast \left(\frac{d\rho}{\rho} \right),F^\ast \left(dy^1 \right),\ldots,F^\ast \left(dy^n \right)$$
$$H^\ast \left(\frac{dx^0}{\rho} \right),\ldots,H^\ast \left(\frac{dx^n}{\rho} \right)$$
extend smoothly to the boundary $\partial  M$ to be frames of ${}^bT^\ast\overline{M}$ and ${}^0T^\ast\overline{M}$ respectively.
It is standard to abuse notation by simply writing that $\sfrac{d\rho}{\rho},dy^1,\ldots,dy^n$ and $\sfrac{dx^0}{\rho},\ldots,\sfrac{dx^n}{\rho}$ are frames for ${}^bT^\ast\overline{M}$ and ${}^0T^\ast\overline{M}$ respectively.
Hence we often write coordinates of ${}^bT^\ast\overline{M}$ and ${}^0T^\ast\overline{M}$ as $\eta_0\sfrac{d\rho}{\rho}+\eta_\mu dy^\mu\mapsto \left(\rho,y^\mu,\eta_0,\eta_\mu \right)$ and ${\overline{\eta}}_i\sfrac{dx^i}{\rho}\mapsto \left(x^i,{\overline{\eta}}_i \right)$ respectively.
For example, if we consider the coordinates $v^i\sfrac{\partial }{\partial  x^i}\mapsto \left(x^i,v^i \right)$ of $TM$, then the canonical identification $H^\ast\circ\flat  : \left.TM\right|_M\rightarrow{}^0T^\ast\overline{M}$ is given by $v^i\sfrac{\partial }{\partial  x^i}\mapsto \left(\rho g_{ij}v^i \right)\sfrac{dx^j}{\rho}$.
\subsection{The b and 0 Cosphere Bundles} \label{Section_4.2} 

Suppose $ \left(\rho,y^1,\ldots,y^n \right)$ are asymptotic boundary normal coordinates of $\overline{M}$ as described in Section \ref{Section_2.2} above.
We know that $T^\ast M$ has a fiber metric $g$.
Thus the maps $F^\ast : T^\ast\overline{M}\rightarrow{}^bT^\ast\overline{M}$ and $H^\ast : T^\ast\overline{M}\rightarrow{}^0T^\ast\overline{M}$ push $g$ to become fiber metrics on $\left.{}^bT^\ast\overline{M}\right|_M$ and $\left.{}^0T^\ast\overline{M}\right|_M$, which we denote by $g_b$ and $g_0$ respectively.
If we consider the boundary frames for ${}^bT^\ast\overline{M}$ and ${}^0T^\ast\overline{M}$ introduced at the end of Section \ref{Section_4.1} above, we have that these metrics are given by (c.f. (\ref{MathItem_2.3}))
\begin{equation} \label{MathItem_4.2} \begin{array}{r@{}l}\left|\eta_0\frac{d\rho}{\rho}+\eta_\mu d y^\mu\right|_{g_b}^2=\eta_0^2+\rho^2h^{\mu\nu}\eta_\mu\eta_\nu,\\\left|{\overline{\eta}}_0\frac{d\rho}{\rho}+{\overline{\eta}}_\mu\frac{dy^\mu}{\rho}\right|_{g_0}^2={\overline{\eta}}_0^2+h^{\mu\nu}{\overline{\eta}}_\mu{\overline{\eta}}_\nu,\\\end{array}\end{equation}
where $ \left(h^{\mu\nu} \right)$ denotes the inverse matrix of $ \left(h_{\mu\nu} \right)$.
From here we see that both $g_b$ and $g_0$ extend smoothly to all of ${}^bT^\ast\overline{M}$ and ${}^0T^\ast\overline{M}$ respectively.
In analog of the \textbf{unit sphere bundle}:
$$SM=\left\{v\in T M : \left|v\right|_g=1\right\},$$
the equations (\ref{MathItem_4.2}) allow us to define the \textbf{unit cosphere bundles} in ${}^bT^\ast\overline{M}$ and ${}^0T^\ast\overline{M}$:
$${}^bS^\ast\overline{M}=\left\{\zeta\in{}^bT^\ast\overline{M} : \left|\zeta\right|_{{}^bg}=1\right\},$$
$${}^0S^\ast\overline{M}=\left\{\overline{\zeta}\in{}^0T^\ast\overline{M} : \left|\overline{\zeta}\right|_{{}^0g}=1\right\}.$$
We note that in \cite{Bibitem_15}, they use the notation ``$\overline{S^\ast M}$" for what we denote by ``${}^bS^\ast\overline{M}$." We let $\pi _b : {}^bS^\ast\overline{M}\rightarrow\overline{M}$ and $\pi _0 : {}^0S^\ast\overline{M}\rightarrow\overline{M}$ denote the natural projection maps.

\begin{remark} \label{MathItem_4.3}  \normalfont Similarly to the remark made in Remark \ref{MathItem_4.1}, we will often identify two points in $SM$, $S^\ast M$, ${}^bS^\ast\overline{M}$, and ${}^0S^\ast\overline{M}$ to be the same if it is possible to go from one to the other by a composition of the maps mentioned there.
\end{remark}

We point out that it is easy to see that both ${}^bS^\ast\overline{M}$ and ${}^0S^\ast\overline{M}$ are smooth embedded submanifolds with boundary of ${}^bT^\ast\overline{M}$ and ${}^0T^\ast\overline{M}$.
We also note that by (\ref{MathItem_4.2}), $g_b$ degenerates over $\partial \overline{M}$ (i.e. stops being positive definite) while $g_0$ does not.
In particular this implies that ${}^bS^\ast\overline{M}$ is not compact while ${}^0S^\ast\overline{M}$ is compact.
\subsection{Splitting the Tangent Bundle} \label{Section_4.3} 

Next we define a natural Riemannian metric on the tangent space $TM$, called the \textbf{Sasaki metric}, generated by $g$.
We recommend that when checking many of the claims below, to check them above the center of normal coordinates since in many cases the expressions simplify considerably due to the vanishing of the Christoffel symbols and the first order partials of $g$.
Consider the tangent bundle's projection map $\pi  : TM\rightarrow M$ and its differential $d\pi  : TTM\rightarrow TM$.
There is another natural map between these tangent spaces called the \textbf{connection map}: $\mathcal{K} : TTM\rightarrow TM$, which is defined as follows.
Take any $\omega\in T_vTM$ and let $\alpha :  \left(a,b \right)\rightarrow M$ be a smooth curve and $V :  \left(a,b \right)\rightarrow TM$ a smooth vector field along $\alpha$ such that $ \left(\alpha,V \right)^\prime \left(0 \right)=\omega$.
Then we set $\mathcal{K} \left(\omega \right)$ to be the covariant derivative
$$\mathcal{K} \left(\omega \right):=\frac{DV}{dt} \left(0 \right).$$
To check that this is independent of the $\alpha$ and $V$ that we choose, a quick computation shows that taking coordinates $ \left(x^i \right)$ of $M$ and the coordinates $v^i\sfrac{\partial }{\partial  x^i}\mapsto \left(x^i,v^i \right)$ of $TM$, $\mathcal{K}$ is given by
$$\mathcal{K} \left(\alpha^i\left.\frac{\partial }{\partial  x^i}\right|_{v^i\sfrac{\partial }{\partial  x^i}}+\beta^i\left.\frac{\partial }{\partial  v^i}\right|_{v^i\sfrac{\partial }{\partial  x^i}} \right)= \left(\beta^k+\Gamma_{ij}^k\alpha^iv^j \right)\frac{\partial }{\partial  x^k},$$
where $\Gamma_{ij}^k$ are the Christoffel symbols of $g$ with respect to $ \left(\sfrac{\partial }{\partial  x^i} \right)$.
Next, an easy exercise shows that the kernels of $d\pi $ and $\mathcal{K}$ partition the tangent bundle's tangent space at any $v\in T_xM$:
\begin{equation} \label{MathItem_4.4} T_vTM={\widetilde{\mathcal{H}}}_v\oplus{\widetilde{\mathcal{V}}}_v\end{equation}
where
$${\widetilde{\mathcal{H}}}_v=\ker{\left.\mathcal{K}\right|_{T_vTM}}\ \ \ \ \ \ \ \ \mathrm{and} \ \ \ \ \ \ \ {\widetilde{\mathcal{V}}}_v=\ker{\left.d\pi \right|_{T_vTM}}.$$
The ``$\widetilde{\mathcal{V}}$" stands for ``vertical" because it can be imagined as being a tangent subspace at $v$ standing vertically above $x$, while the ``$\widetilde{\mathcal{H}}$" stands for ``horizontal." As one can check, both spaces are canonically identified (i.e. isomorphically mapped to) with $T_xM$ by the restricted maps
$$d\pi  : {\widetilde{\mathcal{H}}}_v\longrightarrow T_xM,$$
$$\mathcal{K} : {\widetilde{\mathcal{V}}}_v\longrightarrow T_xM.$$
With this splitting in hand, the Sasaki metric $G$ on $TM$ is defined as follows: for any $\omega,\varsigma\in T_vTM$,
$$\langle \omega,\varsigma\rangle _G=\langle d\pi  \left(\omega \right),d\pi  \left(\varsigma \right)\rangle _g+\langle \mathcal{K} \left(\omega \right),\mathcal{K} \left(\varsigma \right)\rangle _g.$$
It follows immediately that (\ref{MathItem_4.4}) is an orthogonal decomposition with respect to $G$.

We will only work with unit speed geodesics and hence most of our work will be done on the unit sphere bundle
$$SM=\left\{v\in T M : \left|v\right|_g=1\right\}.$$
The Sasaki metric on $TM$ induces a metric on $SM$ which we will also call the Sasaki metric and denote by $G$, relying on context to differentiate the two.
It is not hard to see that at any $v\in SM$ the tangent space of $SM$ splits into the form
$$T_vSM={\widetilde{\mathcal{H}}}_v\oplus\mathcal{V}_v$$
where $\mathcal{V}_v$ is the subspace of ${\widetilde{\mathcal{V}}}_v$ that is $G$-perpendicular to the unit normals to the ``sphere" $S_xM$ above $x$.
Now, take the geodesic vector field $X$ over $SM$.
It is easy to check that $X$ always lies in ${\widetilde{\mathcal{H}}}_v$ for all $v\in SM$ and hence we obtain the splitting
\begin{equation} \label{MathItem_4.5} T_vSM= \left(\mathbb{R}X_v \right)\oplus\mathcal{H}_v\oplus\mathcal{V}_v\end{equation}
where $\mathcal{H}_v$ denotes the orthogonal complement of $\mathbb{R}X_v$ in ${\widetilde{\mathcal{H}}}_v$.
We emphasize that this is an orthogonal decomposition $T_vSM$.
For future use, we point out that restrictions of $d\pi $ and $\mathcal{K}$ map bijectively
\begin{equation} \label{MathItem_4.6} \begin{array}{r@{}l}d\pi  : \mathcal{H}_v\longrightarrow\left\{v^\bot\right\}\subseteq T_xM,\\\mathcal{K} : \mathcal{V}_v\longrightarrow\left\{v^\bot\right\}\subseteq T_xM.\\\end{array}\end{equation}
\subsection{Integration on the Sphere Bundle} \label{Section_4.4} 

Since $SM$ has a Riemannian metric $G$, it has a Riemannian density and hence the Lebesgue measure generated by it (the latter two are independent of orientation).
Hence we may perform Lebesgue integration on $SM$ with respect to $G$.
If $ \left(x^i \right)$ are local coordinates of $M$, $ \left(r_i \right)$ is a frame of $TM$, and we take the coordinates $v^ir_i\mapsto \left(x^i,v^i \right)$ of $TM$, then it turns out that the integral of any function $f\in L^1 \left(SM \right)$ supported over our coordinates is given by the iterated integral
\begin{equation} \label{MathItem_4.7} \int f=\int{\int_{S_xM}{f \left(x^0,\ldots,x^n,v^0,\ldots,v^n \right)dS_x \left(v^0,\ldots,v^n \right)}\sqrt{\det{g}}dx^0\ldots d x^n},\end{equation}
where $ \left(v^0,\ldots,v^n \right)$ are on the sphere $\left|v\right|_g^2=1$ and $dS_x$ is the Lebesgue measure on $S_xM$ induced by $T_xM$ with inner product $g_x$.
We refer the reader to Section 3.6.2 in \cite{Bibitem_42} for a proof.
We point out that the (total) measure of $S_xM$ is the Euclidean surface area of the Euclidean $n$-sphere for all $x\in M$, which for instance follows by looking at the center of normal coordinates.

One example of the usefulness of this observation is the following.
Since $g=\sfrac{\overline{g}}{\rho^2}$ for some smooth metric $\overline{g}$, $\sqrt{\det{g}}$ is $\rho^{- \left(n+1 \right)}$ times ``something smooth" on $\overline{M}$.
So by (\ref{MathItem_4.7}), it follows that any function of the form $\rho^{n+1}L^\infty \left(SM \right)$ is integrable.
\subsection{Splitting the Connection Over the Unit Tangent Bundle} \label{Section_4.5} 

Let us take the natural projection map $\pi  : SM\rightarrow M$.
The pullback (vector) bundle $\pi ^\ast\mathcal{E}$ over $SM$ is defined as the set obtained by taking any point $x\in M$ and attaching a copy of $\mathcal{E}_x$ to every point of the sphere $S_xM$ above it.
Formally,
$$\pi ^\ast\mathcal{E}:=\left\{ \left(v,e \right) : v\in S M,e\in\mathcal{E}_{\pi  \left(v \right)}\right\}.$$
We often canonically identify $ \left(v,e \right)\cong e$ for fixed $v$.
To every fiber $ \left(\pi ^\ast\mathcal{E} \right)_v$ we impose the inner product space structure of $ \left(\mathcal{E}_x,\langle \cdot,\cdot\rangle _{\mathcal{E}_x} \right)$.
If $ \left(b_i \right)$ is a smooth frame for $\mathcal{E}$, then we turn $\pi ^\ast\mathcal{E}$ into a smooth vector bundle over $SM$ (with smooth inner product) by declaring\footnote{ Applying $\pi ^\ast$ to $b_i$ means ``$b_i\circ\pi $."} $ \left(\pi ^\ast b_i \right)$ to be smooth frames for $\pi ^\ast\mathcal{E}$.
The pullback connection $\nabla^{\pi ^\ast\mathcal{E}}=\pi ^\ast\nabla^\mathcal{E}$ in $\pi ^\ast\mathcal{E}$ is defined to be the unique connection so that if $\omega\in TSM$ and $u_0 : M\rightarrow\mathcal{E}$ is smooth, then
\begin{equation} \label{MathItem_4.8} \nabla_\omega^{\pi ^\ast\mathcal{E}} \left(\pi ^\ast u_0 \right)=\pi ^\ast \left(\nabla_{d\pi  \left(\omega \right)}^\mathcal{E}u_0 \right).\end{equation}
For a smooth section $u=u^j\pi ^\ast b_j : SM\rightarrow\pi ^\ast\mathcal{E}$, the pullback connection is explicitly given by
\begin{equation} \label{MathItem_4.9} \nabla_\omega^{\pi ^\ast\mathcal{E}}u=\omega \left(u^j \right)\pi ^\ast b_j+u^j\pi ^\ast \left(\nabla_{d\pi  \left(\omega \right)}^\mathcal{E}b_j \right).\end{equation}
\begin{remark} \label{MathItem_4.10}  \normalfont In the same way we can define the pullback bundles $\pi _b^\ast\mathcal{E}$ and $\pi _0^\ast\mathcal{E}$ on ${}^bS^\ast\overline{M}$ and ${}^0S^\ast\overline{M}$ respectively.
In fact, over the interior we will often identify a section $u\in C^\infty \left(SM;\pi ^\ast\mathcal{E} \right)$ as an element of $C^\infty \left(\left.{}^bS^\ast\overline{M}\right|_M,\left.\pi _b^\ast\mathcal{E}\right|_M \right)$ and $C^\infty \left(\left.{}^0S^\ast\overline{M}\right|_M,\left.\pi _0^\ast\mathcal{E}\right|_M \right)$ via the natural identification $ \left(\pi ^\ast e \right)_v\cong \left(\pi _b^\ast e \right)_\zeta\cong \left(\pi _0^\ast e \right)_{\overline{\zeta}}\cong e$ where $v\cong\zeta\cong\overline{\zeta}$ are identified points on $S_xM$, ${}^0S_x^\ast\overline{M}$, and ${}^bS_x^\ast\overline{M}$ respectively.
\end{remark}

Having defined the splitting of the unit tangent bundle in (\ref{MathItem_4.5}), we now define a natural splitting of the connection of any section $u : SM\rightarrow\pi ^\ast\mathcal{E}$ in the following form:
\begin{equation} \label{MathItem_4.11} \nabla^{\pi ^\ast\mathcal{E}}u\cong{\buildrel X\over\nabla}{}^{\pi ^\ast\mathcal{E}}u+{\buildrel\mathrm{h}\over\nabla}{}^{\pi ^\ast\mathcal{E}}u+{\buildrel\mathrm{v}\over\nabla}{}^{\pi ^\ast\mathcal{E}}u.\end{equation}
The reason we put ``$\cong$" here is that the left-hand side is a tensor of the form $C^\infty \left(SM;T^\ast S M\otimes\pi ^\ast\mathcal{E} \right)$ while we will define the terms on the right-hand side to be tensors of the form $C^\infty \left(SM;TSM\otimes\pi ^\ast\mathcal{E} \right)$.

Let us start by defining ${\buildrel\mathrm{h}\over\nabla}{}^{\pi ^\ast\mathcal{E}}u$.
As we just mentioned, the full connection $\nabla^{\pi ^\ast\mathcal{E}}u$ is a tensor of the form $C^\infty \left(SM;T^\ast S M\otimes\pi ^\ast\mathcal{E} \right)$.
Now, consider the same tensor but with the first index raised with respect to $G$: $\left[\nabla^{\pi ^\ast\mathcal{E}}u\right]^{\sharp}\in C^\infty \left(SM;TSM\otimes\pi ^\ast\mathcal{E} \right)$.
Next, it is an easy exercise to check that there exists a unique linear map
\begin{equation} \label{MathItem_4.12} P_\mathcal{H} : C^\infty \left(SM;TSM\otimes\pi ^\ast\mathcal{E} \right)\longrightarrow C^\infty \left(SM;TSM\otimes\pi ^\ast\mathcal{E} \right)\end{equation}
that satisfies
$$P_\mathcal{H} \left(\omega\otimes e \right)= \left({\mathrm{proj}}_\mathcal{H}\omega \right)\otimes e$$
where ${\mathrm{proj}}_\mathcal{H} : T_vTM\rightarrow\mathcal{H}\subseteq T_vTM$ is the orthogonal projection map onto $\mathcal{H}$.
We then define
$${\buildrel\mathrm{h}\over\nabla}{}^{\pi ^\ast\mathcal{E}}u:=P_\mathcal{H} \left(\left[\nabla^{\pi ^\ast\mathcal{E}}u\right]^{\sharp} \right).$$
As noted above, this is a tensor of the form $C^\infty \left(SM;TSM\otimes\pi ^\ast\mathcal{E} \right)$.

We define ${\buildrel X\over\nabla}{}^{\pi ^\ast\mathcal{E}}u$ and ${\buildrel\mathrm{v}\over\nabla}{}^{\pi ^\ast\mathcal{E}}u$ in the same way but instead use analogous map $P_{\mathbb{R}X}$, ${\mathrm{proj}}_{\mathbb{R}X}$ and $P_\mathcal{V}$, ${\mathrm{proj}}_\mathcal{V}$ respectively.
However instead of using ${\buildrel X\over\nabla}{}^{\pi ^\ast\mathcal{E}}u$, it is more common to use the related quantity
\begin{equation} \label{MathItem_4.13} \mathbb{X}u:=\nabla_X^{\pi ^\ast\mathcal{E}}u.\end{equation}
The quantities $\mathbb{X}u$ and ${\buildrel X\over\nabla}{}^{\pi ^\ast\mathcal{E}}u$ are equivalent in the sense that knowing one allows you to compute the other.
Hence we often record the decomposition (\ref{MathItem_4.11}) instead as
\begin{equation} \label{MathItem_4.14} \nabla^{\pi ^\ast\mathcal{E}}u= \left(\mathbb{X}u,{\buildrel\mathrm{h}\over\nabla}{}^{\pi ^\ast\mathcal{E}}u,{\buildrel\mathrm{v}\over\nabla}{}^{\pi ^\ast\mathcal{E}}u \right).\end{equation}
The second two components are called the \textbf{horizontal and vertical derivatives} of $u$ respectively.
However, it is convenient to change the interpretation of the latter two derivatives as follows.

We define the bundle $N$ over $SM$ by attaching to every $v\in S_xM$ a copy of $\left\{v^\bot\right\}\subseteq T_xM$.
Formally,
$$N=\left\{ \left(v,w \right) : v\in S_xM\ \ \mathrm{where} \ \ x\in M\ \ \mathrm{and} \ \ w\in\left\{v^\bot\right\}\right\}.$$
To every fiber $N_v$ we impose the inner product space structure of $ \left(\left\{v^\bot\right\},g_x \right)$ which we denote by ``$\langle \cdot,\cdot\rangle _{N_v}$." It is an easy exercise to show that this is a smooth subbundle of $\pi ^\ast TM$.
By (\ref{MathItem_4.6}) we can think of $d\pi $ and $\mathcal{K}$ as mapping bijectively
\begin{equation} \label{MathItem_4.15} \begin{array}{r@{}l}d\pi  : \mathcal{H}_v\longrightarrow N_v,\\\mathcal{K} : \mathcal{V}_v\longrightarrow N_v.\\\end{array}\end{equation}
At every $v\in TM$, the maps $P_\mathcal{H}$ and $P_\mathcal{V}$ above map into $\mathcal{H}_v\otimes \left(\pi ^\ast\mathcal{E} \right)_v$ and $\mathcal{V}_v\otimes \left(\pi ^\ast\mathcal{E} \right)_v$ respectively, thus using the identification (\ref{MathItem_4.15}) we can think of the horizontal and vertical derivatives as both being $N\otimes\pi ^\ast\mathcal{E}$-valued:
$${\buildrel\mathrm{h}\over\nabla}{}^{\pi ^\ast\mathcal{E}}u\in C^\infty \left(SM;N\otimes\pi ^\ast\mathcal{E} \right)\ \ \ \ \ \ \ \ \mathrm{and} \ \ \ \ \ \ \ {\buildrel\mathrm{v}\over\nabla}{}^{\pi ^\ast\mathcal{E}}u\in C^\infty \left(SM;N\otimes\pi ^\ast\mathcal{E} \right).$$
We mention that we assign the natural inner product on $N\otimes\pi ^\ast\mathcal{E}$ (i.e. the unique one satisfying $\langle z\otimes e,z^\prime\otimes e^\prime\rangle _{N\otimes\pi ^\ast\mathcal{E}}=\langle z,z^\prime\rangle _N\langle e,e^\prime\rangle _{\pi ^\ast\mathcal{E}}$).
The reason that this interpretation is useful is that it becomes natural to apply well-known adjoint formulas for the horizontal and vertical derivatives over this space.
In particular, it turns out that there are differential operators
\begin{equation} \label{MathItem_4.16} \begin{array}{r@{}l}{\buildrel\mathrm{h}\over{\mathrm{div}}}{}^{\pi ^\ast\mathcal{E}} : C^\infty \left(SM;N\otimes\pi ^\ast\mathcal{E} \right)\longrightarrow C^\infty \left(SM;\pi ^\ast\mathcal{E} \right),\\{\buildrel\mathrm{v}\over{\mathrm{div}}}{}^{\pi ^\ast\mathcal{E}} : C^\infty \left(SM;N\otimes\pi ^\ast\mathcal{E} \right)\longrightarrow C^\infty \left(SM;\pi ^\ast\mathcal{E} \right),\\\end{array}\end{equation}
with the property that if $u\in C^\infty \left(SM;\pi ^\ast\mathcal{E} \right)$ and $w\in C^\infty \left(SM;N\otimes\pi ^\ast\mathcal{E} \right)$ are such that at least one of them is of compact support, then
\begin{equation} \label{MathItem_4.17} \begin{array}{r@{}l}\langle {\buildrel\mathrm{h}\over\nabla}{}^{\pi ^\ast\mathcal{E}}u,w\rangle _{L^2 \left(SM;N\otimes\pi ^\ast\mathcal{E} \right)}=-\langle u,{\buildrel\mathrm{h}\over{\mathrm{div}}}{}^{\pi ^\ast\mathcal{E}}w\rangle _{L^2 \left(SM;\pi ^\ast\mathcal{E} \right)},\\\langle {\buildrel\mathrm{v}\over\nabla}{}^{\pi ^\ast\mathcal{E}}u,w\rangle _{L^2 \left(SM;N\otimes\pi ^\ast\mathcal{E} \right)}=-\langle u,{\buildrel\mathrm{v}\over{\mathrm{div}}}{}^{\pi ^\ast\mathcal{E}}w\rangle _{L^2 \left(SM;\pi ^\ast\mathcal{E} \right)}.\\\end{array}\end{equation}
The operators in (\ref{MathItem_4.16}) are naturally called the \textbf{horizontal and vertical divergences} respectively.
Since we will only make use of the vertical divergence, we only justify its existence.
To do this, fix any $x\in M$ and consider the sphere $S_xM$ as a Riemannian submanifold of $ \left(SM,G \right)$.
Let ${\mathrm{div}}_{S_xM}$ denote the ``divergence" on $S_xM$ that takes any vector field $V\in C^\infty \left(S_xM;TS_xM \right)$ and outputs a scalar function ${\rm div}_{S_xM}{V}\in C^\infty \left(S_xM \right)$ (e.g. see (2.19) and Proposition 2.46 in \cite{Bibitem_26}).

Now, choose any coordinates $ \left(x^i \right)$ of $M$ and any \textit{orthonormal} frame $ \left(b_i \right)$ for $\mathcal{E}$ with the same domain.
Take any $u\in C^\infty \left(SM;\pi ^\ast\mathcal{E} \right)$ and $w\in C^\infty \left(SM;N\otimes\pi ^\ast\mathcal{E} \right)$ which in coordinates we can write as
$$u=u^i\pi ^\ast b_i\ \ \ \ \ \mathrm{and} \ \ \ \ w=w^{ij}\pi ^\ast\frac{\partial }{\partial  x^i}\otimes\pi ^\ast b_j.$$
Let us first look at $w$.
For fixed $x\in M$ and fixed $v\in S_xM$, we have by the definition of $N$ that for fixed $j\in\left\{1,\ldots,\dim{\mathcal{E}}\right\}$, $w^{ij}\sfrac{\partial }{\partial  x^i}\in\left\{v^\bot\right\}$.
An explicit expression for this is
\begin{equation} \label{MathItem_4.18}  \left(g_x \right)_{ki}v^kw^{ij}=0.\end{equation}
Thus if we consider the coordinates $ \left(x^i,v^i \right)$ of $TM$ given by $v^i\sfrac{\partial }{\partial  x^i}\mapsto \left(x^i,v^i \right)$, (\ref{MathItem_4.18}) tells us that $w^{ij}\sfrac{\partial }{\partial  v^i}$ is a smooth \textit{tangent} vector field to $S_xM$.
Hence we may apply ${\mathrm{div}}_{S_xM}$ to $w^{ij}\sfrac{\partial }{\partial  v^i}$.
We define
\begin{equation} \label{MathItem_4.19} {\buildrel\mathrm{v}\over{\mathrm{div}}}{}^{\pi ^\ast\mathcal{E}}w= \left({\rm div}_{S_xM}{ \left(w^{ij}\frac{\partial }{\partial  v^i} \right)} \right)\pi ^\ast b_j.\end{equation}
We will justify shortly why this is coordinate and frame invariant.
To check that this satisfies the desired second equation in (\ref{MathItem_4.17}), let us explicitly compute ${\buildrel\mathrm{v}\over\nabla}{}^{\pi ^\ast\mathcal{E}}u$.
Let ${}^\mathcal{E}\Gamma_{ij}^k$ denote the connection symbols of $\nabla^\mathcal{E}$ with respect to $ \left(\sfrac{\partial }{\partial  x^i} \right)$ and $ \left(b_k \right)$.
Using (\ref{MathItem_4.9}) we have that (here $ \left(dx^m \right)$ is the dual frame of $ \left(\sfrac{\partial }{\partial  x^i} \right)$ with respect to $g_x$, and $\sharp $ is associated to $G$)
$$\nabla^{\pi ^\ast\mathcal{E}}u=du^i\otimes\pi ^\ast b_i+{}^\mathcal{E}\Gamma_{mi}^ku^idx^m\otimes\pi ^\ast b_k$$
$$\Longrightarrow\ \ \ \ \ {\buildrel\mathrm{v}\over\nabla}{}^{\pi ^\ast\mathcal{E}}u=P_\mathcal{V} \left(du^i\otimes\pi ^\ast b_i+{}^\mathcal{E}\Gamma_{mi}^ku^idx^m\otimes\pi ^\ast b_k \right)^{\sharp}$$
$$=\mathcal{K} \left({\mathrm{proj}}_\mathcal{V} \left(du^i \right)^{\sharp} \right)\otimes\pi ^\ast b_i+{}^\mathcal{E}\Gamma_{mi}^ku^i\mathcal{K} \left({\mathrm{proj}}_\mathcal{V} \left(dx^m \right)^{\sharp} \right)\otimes\pi ^\ast b_k$$
where $\mathcal{K}$ is as in (\ref{MathItem_4.15}).
It is not hard to check that by construction $\mathcal{V}_v=T_vS_xM$ for any $v\in S_xM$.
Hence $ \left(dx^m \right)^{\sharp}\bot\mathcal{V}_v$ since for any $\omega$ tangent to $S_xM$, $dx^m \left(\omega \right)=0$.
It also quickly follows that ${\mathrm{proj}}_\mathcal{V} \left(du^i \right)^{\sharp}={\grad}_{S_xM}{u^i}$.
Thus
\begin{equation} \label{MathItem_4.20} \nabla^{\pi ^\ast\mathcal{E}}u=\mathcal{K} \left({\grad}_{S_xM}{u^i} \right)\otimes\pi ^\ast b_i.\end{equation}
Let us fix $x\in M$ and integrate over $S_xM$.
It is not hard to see that $\mathcal{K} \left(v^i\sfrac{\partial }{\partial  v^i} \right)=\sfrac{v^i\pi ^\ast\partial }{\partial  x^i}$ and so
$$\langle {\buildrel\mathrm{v}\over\nabla}{}^{\pi ^\ast\mathcal{E}}u,w\rangle _{L^2 \left(S_xM;N\otimes\pi ^\ast\mathcal{E} \right)}=\int_{S_xM}{\langle \mathcal{K} \left({\grad}_{S_xM}{u^r} \right),w^{ij}\pi ^\ast\frac{\partial }{\partial  x^i}\rangle _N\langle \pi ^\ast b_r,\pi ^\ast b_j\rangle _{\pi ^\ast\mathcal{E}}dS_xM}$$
$$=\int_{S_xM}{\langle {\grad}_{S_xM}{u^r},w^{ij}\frac{\partial }{\partial  v^i}\rangle _N\langle \pi ^\ast b_r,\pi ^\ast b_j\rangle _{\pi ^\ast\mathcal{E}}dS_xM}.$$
Using that ``${\mathrm{div}}_{S_xM}$" is the formal adjoint of ``${\mathrm{grad}}_{S_xM}$" (c.f. Exercise 16.2 (b) in \cite{Bibitem_25}), this is equal to
$$=\int_{S_xM}{\langle u^r,{\rm div}_{S_xM}{w^{ij}\frac{\partial }{\partial  v^i}}\rangle _N\langle \pi ^\ast b_r,\pi ^\ast b_j\rangle _{\pi ^\ast\mathcal{E}}dS_xM}=\langle u,{\buildrel\mathrm{v}\over{\mathrm{div}}}{}^{\pi ^\ast\mathcal{E}}w\rangle _{L^2 \left(S_xM;\pi ^\ast\mathcal{E} \right)}.$$
In other words, a stronger statement holds: $\nabla^{\pi ^\ast\mathcal{E}}$ and ${\buildrel\mathrm{v}\over{\mathrm{div}}}{}^{\pi ^\ast\mathcal{E}}$ are $L^2$ formal adjoints on each sphere $S_xM$ rather than simply on all of $SM$.
For later use, we write this down explicitly:
$$\int_{S_xM}{\langle {\buildrel\mathrm{v}\over\nabla}{}^{\pi ^\ast\mathcal{E}}u,w\rangle _{N\otimes\pi ^\ast\mathcal{E}}dS_xM}=\int_{S_xM}{\langle u,{\buildrel\mathrm{v}\over{\mathrm{div}}}{}^{\pi ^\ast\mathcal{E}}w\rangle _{\pi ^\ast\mathcal{E}}dS_xM}.$$
By varying $u$ in this equality, this shows that the vertical divergence is unique and hence (\ref{MathItem_4.19}) is indeed independent of coordinates and frame.
Because one can iterate integrals over $SM$ as in (\ref{MathItem_4.7}), integrating both sides of the above equation in $x\in M$ shows that the second equation in (\ref{MathItem_4.17}) indeed holds if either $u$ or $w$ has compact support.
We mention that the operators ``${\mathrm{grad}}_{S_xM}$" and ``${\mathrm{div}}_{S_xM}$" are known as ``${\buildrel\mathrm{v}\over\nabla}$" and ``${\buildrel\mathrm{v}\over{\mathrm{div}}}$" respectively in the literature (e.g. see \cite{Bibitem_15}, \cite{Bibitem_18}, and \cite{Bibitem_41}).

For future reference, we end this section with one more definition.
We define the differential operator
$$\mathbb{X} : C^\infty \left(SM;N\otimes\pi ^\ast\mathcal{E} \right)\rightarrow C^\infty \left(SM;N\otimes\pi ^\ast\mathcal{E} \right),$$
differentiated from the $\mathbb{X}$ introduced in (\ref{MathItem_4.13}) by context, to be the unique operator satisfying that for any $Z\otimes b\in C^\infty \left(SM;N\otimes\pi ^\ast\mathcal{E} \right)$,
\begin{equation} \label{MathItem_4.21} \mathbb{X} \left(Z\otimes b \right)=X \left(Z \right)\otimes b+Z\otimes\mathbb{X} \left(b \right),\end{equation}
where $X \left(Z \right)$ at any point $v\in SM$ denotes the covariant derivative of $Z$ along the unit speed geodesic $\gamma$ with initial velocity $v$ at time $t=0$:
\begin{equation} \label{MathItem_4.22} \left.X \left(Z \right)\right|_v=\frac{D_\gamma Z}{dt} \left(0 \right).\end{equation}
Since $Z\bot\dot{\gamma}$ implies that $\sfrac{D_\gamma Z}{dt}\ \bot\dot{\gamma}$, we see that $\mathbb{X}$ indeed maps into smooth sections of $N\otimes\pi ^\ast\mathcal{E}$ (i.e. not simply into $\pi ^\ast TM\otimes\pi ^\ast\mathcal{E}$).
\subsection{Curvatures} \label{Section_4.6} 

We now cover the curvature operators of $\mathcal{E}$, $\pi ^\ast\mathcal{E}$, and simply the metric $g$.
We start with the first one.
The operator $\nabla^\mathcal{E}$ maps between the following spaces of sections:
$$\nabla^\mathcal{E} : C^\infty \left(\overline{M};\mathcal{E} \right)\longrightarrow C^\infty \left(\overline{M};T^\ast\overline{M}\otimes\mathcal{E} \right)$$
Let $\Lambda^k \left(T^\ast\overline{M} \right)$ denote the bundle of covariant alternating $k$-tensors and let 
$$C^\infty \left(\overline{M};\Lambda^k \left(T^\ast\overline{M} \right)\otimes\mathcal{E} \right)$$
denote the space of smooth sections of $T^\ast\overline{M}\otimes\ldots\otimes T^\ast\overline{M}\otimes\mathcal{E}$ that are alternating in their first $k$ arguments.
The operators
$$\nabla^\mathcal{E} : C^\infty \left(\overline{M};\Lambda^k \left(T^\ast\overline{M} \right)\otimes\mathcal{E} \right)\longrightarrow C^\infty \left(\overline{M};\Lambda^{k+1} \left(T^\ast\overline{M} \right)\otimes\mathcal{E} \right)$$
are defined to be the unique operators that satisfy that for any $\theta\in C^\infty \left(\overline{M};\Lambda^k \left(T^\ast\overline{M} \right) \right)$ and any $u\in C^\infty \left(M;\mathcal{E} \right)$
$$\nabla^\mathcal{E} \left(\theta\otimes u \right)=d\theta\otimes u+ \left(-1 \right)^k\theta\land\nabla^\mathcal{E}u,$$
where $\theta\land\nabla^\mathcal{E}u$ denotes the wedge-like product:
$$ \left(\theta\land\nabla^\mathcal{E}u \right) \left(v_1,\ldots,v_{k+1},l \right)=\frac{1}{k!1!}\sum_{\sigma\in S_{k+1}}{\mathrm{sgn} \left(\sigma \right)\theta \left(v_{\sigma \left(1 \right)},\ldots,v_{\sigma \left(k \right)} \right)\nabla^\mathcal{E}u \left(v_{\sigma \left(k+1 \right)},l \right)},$$
where each $v_i\in T_x\overline{M}$, $l\in\mathcal{E}_x^\ast$, and $S_{k+1}$ denotes the set of permutations of $k+1$ elements.

The \textbf{curvature} of $\nabla^\mathcal{E}$ is defined to be
$$f^\mathcal{E}:=\nabla^\mathcal{E}\circ\nabla^\mathcal{E} : C^\infty \left(\overline{M};\mathcal{E} \right)\longrightarrow C^\infty \left(\overline{M};\Lambda^2 \left(T^\ast\overline{M} \right)\otimes\mathcal{E} \right).$$
A straightforward computation shows that in any coordinates $ \left(x^i \right)$ of $\overline{M}$ and any frame $ \left(b_i \right)$ of $\mathcal{E}$, the curvature $f^\mathcal{E}$ applied to any smooth section $u\in C^\infty \left(\overline{M};\mathcal{E} \right)$ is given by
\begin{equation} \label{MathItem_4.23} f^\mathcal{E}u=u^l \left(\frac{\partial {}^\mathcal{E}\Gamma_{jl}^k}{\partial  x^i}-{}^\mathcal{E}\Gamma_{il}^m{}^\mathcal{E}\Gamma_{jm}^k-\frac{\partial {}^\mathcal{E}\Gamma_{il}^k}{\partial  x^j}+{}^\mathcal{E}\Gamma_{jl}^m{}^\mathcal{E}\Gamma_{im}^k \right)dx^i\otimes dx^j\otimes b_k.\end{equation}
where ${}^\mathcal{E}\Gamma_{ij}^k$ are the connection symbols of $\nabla^\mathcal{E}$ with respect to $ \left(\sfrac{\partial }{\partial  x^i} \right)$ and $ \left(b_k \right)$.
The resemblance of this tensor to the Riemann curvature tensor is the motivation for the name of $f^\mathcal{E}$. 

Next we define a curvature operator associated to $f^\mathcal{E}$ which acts over $SM$.
Notice that $f^\mathcal{E}$ can be viewed as a $C^\infty \left(\overline{M};\Lambda^2 \left(T^\ast\overline{M} \right)\otimes\mathcal{E}\otimes\mathcal{E}^\ast \right)$ tensor field by thinking of the $u$ in (\ref{MathItem_4.23}) as the fourth argument of $f^\mathcal{E}$ (i.e. the argument of $\mathcal{E}^\ast$).
Hence it can also be canonically identified with a map, denoted by the same letter, of the form
$$f^\mathcal{E} : C^\infty \left(\overline{M};T\overline{M} \right)\times C^\infty \left(\overline{M};\mathcal{E} \right)\longrightarrow C^\infty \left(\overline{M};T^\ast\overline{M}\otimes\mathcal{E} \right).$$
In our coordinates and frames it is given by the following: if $f_{ij}{}^k{}_l$ denotes the tensor component written out in the parenthesis ``$ \left(\ldots \right)$" in (\ref{MathItem_4.23}), then for any $x\in\overline{M}$, $v\in T_x\overline{M}$, and $e\in\mathcal{E}_x$,
\begin{equation} \label{MathItem_4.24} f_x^\mathcal{E} \left(v,e \right)=e^lf_{ij}{}^k{}_lv^idx^j\otimes b_k.\end{equation}
We define the \textbf{curvature operator} associated to $f^\mathcal{E}$ to be the map
$$F^\mathcal{E} : C^\infty \left(SM;\pi ^\ast\mathcal{E} \right)\longrightarrow C^\infty \left(SM;N\otimes\pi ^\ast\mathcal{E} \right)$$
given by the following.
For any $x\in M$, $v\in S_xM$, and $e\in \left(\pi ^\ast\mathcal{E} \right)_v$,
\begin{equation} \label{MathItem_4.25} F_v^\mathcal{E} \left(e \right):=\left[f_x^\mathcal{E} \left(v,e \right)\right]^{\sharp}.\end{equation}
where $\sharp $ raises the first index of $f_x^\mathcal{E} \left(v,e \right)$ (i.e. $j$ in (\ref{MathItem_4.24})).
From (\ref{MathItem_4.23}) we see that the component $f_{ij}{}^k{}_l$ in (\ref{MathItem_4.24}) is antisymmetric in $i$ and $j$, from which a quick computation shows that (\ref{MathItem_4.25}) is perpendicular to $v\otimes e^\prime\in \left(\pi ^\ast T M \right)_v\otimes \left(\pi ^\ast\mathcal{E} \right)_v$ for any $e^\prime\in \left(\pi ^\ast\mathcal{E} \right)_v$.
In particular, this implies that $F^\mathcal{E}$ is indeed $N\otimes\pi ^\ast\mathcal{E}$-valued.

The last curvature quantity that we want to establish notation for is the ordinary curvature of $g$.
Let
$$R : C^\infty \left(M;TM \right)\times C^\infty \left(M;TM \right)\times C^\infty \left(M;TM \right)\longrightarrow C^\infty \left(M;TM \right)$$
denote the Riemann curvature endomorphism given by
$$R \left(X,Y \right)Z=\nabla_X\nabla_YZ-\nabla_Y\nabla_XZ-\nabla_{\left[X,Y\right]}Z,$$
where $\nabla$ is the Levi-Civita connection and $\left[\cdot,\cdot\right]$ is the Lie bracket.
Recall that $R$ is multi-linear over $C^\infty \left(M \right)$ and thus has well-defined restrictions to $ \left(T_xM \right)^3\rightarrow T_xM$ for any fixed $x\in M$.
We define operators, denoted by the same letter,
$$R : C^\infty \left(SM;N \right)\longrightarrow C^\infty \left(SM;N \right)$$
$$R : C^\infty \left(SM;N\otimes\pi ^\ast\mathcal{E} \right)\longrightarrow C^\infty \left(SM;N\otimes\pi ^\ast\mathcal{E} \right)$$
to be the unique ones satisfying that for any $x\in M$, $v\in S_xM$, $w\in N_v$, and $e\in \left(\pi ^\ast\mathcal{E} \right)_v$,
$$R_v \left(w \right):=R_x \left(w,v \right)v\ \ \ \ \ \ \ \ \mathrm{and} \ \ \ \ \ \ \ R_v \left(w\otimes e \right):=\left[R_x \left(w,v \right)v\right]\otimes e.$$

\subsection{Connection on Endomorphism Fields} \label{Section_4.7} 

The connection $\nabla^\mathcal{E}$ in $\mathcal{E}$ over $\overline{M}$ generates a natural connection $\nabla^{\End{\mathcal{E}}}$ in the endomorphism bundle $\End{\mathcal{E}}$ over $\overline{M}$ as follows.
For any $U\in C^\infty \left(\overline{M};\End{\mathcal{E}} \right)$ and any $v\in T_x\overline{M}$, we define $\nabla_v^{\End{\mathcal{E}}}U$ to be the unique element of $\End{\mathcal{E}_x}$ satisfying
$$ \left(\nabla_v^{\End{\mathcal{E}}}U \right)h=\left[\nabla_v^\mathcal{E},U\right]h,$$
for any $h\in C^\infty \left(\overline{M};\mathcal{E} \right)$.
This is considered a natural connection on $\End{\mathcal{E}}$ because it immediately follows that it satisfies the product rule $\nabla_v^\mathcal{E} \left(Uh \right)= \left(\nabla_v^{\End{\mathcal{E}}}U \right)h+U\nabla_v^\mathcal{E}h$.
Let us see what this connection looks like in coordinates.
Take coordinates $ \left(x^i \right)$ of $\overline{M}$, a frame $ \left(b_i \right)$ for $\mathcal{E}$ over their domain, and let ${}^\mathcal{E}\Gamma_{ij}^k$ denote the connection symbols of $\nabla^\mathcal{E}$ with respect to $ \left(\sfrac{\partial }{\partial  x^i} \right)$ and $ \left(b_k \right)$.
Then, a quick computation shows that
\begin{equation} \label{MathItem_4.26} \nabla_v^{\End{\mathcal{E}}}U=vU+ \left({}^\mathcal{E}\Gamma \right)U-U \left({}^\mathcal{E}\Gamma \right),\end{equation}
where on the right-hand side $U$ is thought of as a matrix in the basis $ \left(b_i \right)$, $vU$ denotes applying $v$ to every entry of $U$, and ${}^\mathcal{E}\Gamma$ represents the matrix with entry ${}^\mathcal{E}\Gamma_{ij}^kv^i$ in the $k^{\mathrm{th}}$ row and $j^{\mathrm{th}}$ column.
We will use the Frobenius inner product $\langle \cdot,\cdot\rangle _{\End{\mathcal{E}}}$ on $\End{E}$ given by
\begin{equation} \label{MathItem_4.27} \langle A,B\rangle _{\End{\mathcal{E}}}=\mathrm{Trace} \left(AB^\ast \right).\end{equation}
Of more importance to us will be the pullback bundle $\pi ^\ast\End{\mathcal{E}}$ on $SM$ with inner product $\langle \cdot,\cdot\rangle _{\End{\mathcal{E}}}$ imposed on its fibers, whose connection we denote by $\nabla^{\pi ^\ast\End{\mathcal{E}}}:=\pi ^\ast\nabla^{\End{\mathcal{E}}}$.
We analogously define the operator $\mathbb{X}=\nabla_X^{\pi ^\ast\End{\mathcal{E}}}$, using context to differentiate it from our other operators also denoted by ``$\mathbb{X}$."
\section{Gauge Equivalence of Connections and Higgs fields} \label{Section_5} 

In this section we build the necessary tools to prove Theorem \ref{MathItem_2.8}, and then prove it at the end.
Throughout this section we suppose that $ \left(M\subseteq\overline{M},g \right)$ is an asymptotically hyperbolic (AH) space, $\rho$ is a boundary defining function, $ \left(\mathcal{E},\langle \cdot,\cdot\rangle _\mathcal{E} \right)$ is a smooth complex Hermitian vector bundle over $\overline{M}$, and that $\nabla^\mathcal{E}$ is a smooth connection in $\mathcal{E}$.

\subsection{The Regularity Spaces} \label{Section_5.1} 

For the majority of this paper we will be working with the solution to a transport equation over $SM$.
In particular, we will be making use of its $L^2$ norm after we apply several differential operators to it that extend smoothly to the boundary of the 0-cosphere bundle.
For this reason, we will make use of the following regularity spaces:

\begin{definition} \label{MathItem_5.1}  Suppose that $ \left(M\subseteq\overline{M},g \right)$ is an asymptotically hyperbolic space,\textit{ }$\rho$ is a boundary defining function, $ \left(\mathcal{E},\langle \cdot,\cdot\rangle _\mathcal{E} \right)$ is a smooth complex Hermitian vector bundle over $\overline{M}$, and that $\nabla^\mathcal{E}$ is a smooth connection in $\mathcal{E}$.
For any fixed integer $k\geq0$, we define the spaces $\mathcal{R}^k \left(SM;\pi ^\ast\mathcal{E} \right)$ and $\mathcal{R}^k \left(SM;N\otimes\pi ^\ast\mathcal{E} \right)$ of order $k$ to be the spaces of smooth sections $u\in C^\infty \left(SM;\pi ^\ast\mathcal{E} \right)$ and $w\in C^\infty \left(SM;N\otimes\pi ^\ast\mathcal{E} \right)$ respectively that satisfy the following.

We say that $u\in\mathcal{R}^k \left(SM;\pi ^\ast\mathcal{E} \right)$ if for any smooth vector fields $V_1,...,V_k\in C^\infty \left({}^0S^\ast\overline{M};T{}^0S^\ast\overline{M} \right)$ over the 0-cosphere bundle, any frame $ \left(b_i \right)$ of $\mathcal{E}$, and any compact subset $K\subseteq\dom{ \left(b_i \right)}\subseteq\overline{M}$, $u=u^j\pi ^\ast b_j$ must satisfy that
\begin{equation} \label{MathItem_5.2} u^j,\ \ \ V_1u^j,\ \ \ V_2V_1u^j,\ \ldots\ \ \ ,\ V_k\ldots V_1u^j\ \ \mathrm{are\ all\ in}\ \ \rho^{\sfrac{ \left(n+1 \right)}{2}}L^\infty\left[\pi ^{-1}\left[K\right]\right].\end{equation}
We say that $w\in\mathcal{R}^k \left(SM;N\otimes\pi ^\ast\mathcal{E} \right)$ if for any smooth vector fields $V_1,...,V_k\in C^\infty \left({}^0S^\ast\overline{M};T{}^0S^\ast\overline{M} \right)$, any boundary coordinates $ \left(\rho,y^\mu \right)= \left(x^i \right)$ or interior coordinates $ \left(x^i \right)$ of $\overline{M}$, any frame $ \left(b_i \right)$ of $\mathcal{E}$ over these coordinates' domain, $w$ is of the form
\begin{equation} \label{MathItem_5.3} w=\rho w^{ij}\pi ^\ast\frac{\partial }{\partial  x^i}\otimes\pi ^\ast b_j\end{equation}
(note the $\rho$ in front), and satisfies that for any compact subset $K\subseteq\dom{ \left(x^i \right)}\subseteq\overline{M}$,
\begin{equation} \label{MathItem_5.4} w^{ij},\ \ \ V_1w^{ij},\ \ \ V_2V_1w^{ij},\ldots\ \ \ ,\ V_k\ldots V_1w^{ij}\ \ \mathrm{are\ all\ in}\ \ \rho^{\sfrac{ \left(n+1 \right)}{2}}L^\infty\left[\pi ^{-1}\left[K\right]\right].\end{equation}
\end{definition}

We sometimes omit writing the domain of our space and simply write $\mathcal{R}^k$ when either the domain is of no importance or clear from context.
The power ``$\sfrac{ \left(n+1 \right)}{2}$" above is chosen to ensure that both $u$ and $w$ as above are in $L^2$ - see Remark \ref{MathItem_5.8} below.

\begin{remark} \label{MathItem_5.5}  \normalfont We point out that (\ref{MathItem_5.2}) and (\ref{MathItem_5.4}) only need to be checked in an atlas of $\overline{M}$ and do not depend on the chosen boundary defining function $\rho$.
This is clear for interior charts.
To see why this holds near the boundary, take any pair of boundary coordinates and frames $ \left(\rho,y^\mu \right)= \left(x^i \right), \left(b_j \right)$ and $ \left(\widetilde{\rho},{\widetilde{y}}^\mu \right)= \left({\widetilde{x}}^i \right), \left({\widetilde{b}}_j \right)$ whose domains intersect.
The coefficients $\alpha_i^j$ and $\beta_i^j$ in the transformation laws $\sfrac{\partial }{\partial  x^i}=\alpha_i^j\sfrac{\partial }{\partial {\widetilde{x}}^i}$ and $b_i=\beta_i^j{\widetilde{b}}_j$ will be smooth on $\overline{M}$ and hence smooth on ${}^0S^\ast\overline{M}$ (when lifted).
Thus $\alpha_i^j$ and $\beta_i^j$ will be bounded over $\pi _0^{-1}\left[K\right]$ for any compact subset $K\subseteq\dom{ \left(x^i \right)}\cap\dom{ \left({\widetilde{x}}^i \right)}$ of $\overline{M}$, and the same thing will hold if we apply smooth vectors $V_1,\ldots,V_k$ to them as in the above definition.
We also have that any two boundary defining functions are comparable (i.e. both $\sfrac{\rho}{\widetilde{\rho}},\sfrac{\widetilde{\rho}}{\rho}\in C^\infty \left(\overline{M} \right)$).
From this it quickly follows that (\ref{MathItem_5.2}) and (\ref{MathItem_5.4}) only need to be checked in an atlas of $\overline{M}$ and do not depend on the chosen boundary defining function.
\end{remark}

\begin{remark} \label{MathItem_5.6}  \normalfont Regarding continuity at infinity, from (\ref{MathItem_5.2}) it follows that any element of $\mathcal{R}^k \left(SM;\pi ^\ast\mathcal{E} \right)$ extends continuously to ${}^0S^\ast\overline{M}$ if we identify it as an element of $C^k \left(\left.{}^0S^\ast\overline{M}\right|_M,\left.\pi _0^\ast\mathcal{E}\right|_M \right)$ (see Remark \ref{MathItem_4.10}) by setting it to be zero on the boundary $\partial {}^0S^\ast\overline{M}$.

Although we will not make use of this, we also mention that one can similarly extend elements of $\mathcal{R}^k \left(SM;N\otimes\pi ^\ast\mathcal{E} \right)$ to ${}^0S^\ast\overline{M}$ as follows.
First one defines the following vector bundle over ${}^0S^\ast\overline{M}$ analogous to $N$:
$${}^0N=\left\{ \left(\overline{\zeta},\vartheta \right) : \overline{\zeta}\in{}^0S_x^\ast\overline{M}\ \ \ \mathrm{where} \ \ \ x\in\overline{M}\ \ \ \mathrm{and} \ \ \ \vartheta\in\ker{\overline{\zeta}}\subseteq{}^0T_x\overline{M}\right\},$$
where in the last part we view $\overline{\zeta}$ as a linear functional $\overline{\zeta} : {}^0T_x\overline{M}\rightarrow\mathbb{R}$.
To see the analogy with $N$, it not hard to see that over $M$, $\ker{\overline{\zeta}}=\{({\overline{\zeta}}^{\sharp}){}^\bot\}$ where $\sharp $ and $\bot$ are computed using $g$ and the canonical identifications $\left.{}^0T^\ast\overline{M}\right|_M\cong TM$ and $\left.{}^0T\overline{M}\right|_M\cong TM$.
It's not hard to check that ${}^0N$ is a smooth subbundle of $\pi _0^\ast \left({}^0T\overline{M} \right)$.
With respect to the coordinates and frames in Definition \ref{MathItem_5.1}, the latter is spanned by $ \left(\rho\sfrac{\partial }{\partial  x^i} \right)$ and hence any section $w\in C^k \left({}^0S^\ast\overline{M};{}^0N\otimes\pi _0^\ast\mathcal{E} \right)$ is of the form
$$w=w^{ij}\pi _0^\ast \left(\rho\sfrac{\partial }{\partial  x^i} \right)\otimes\pi _0^\ast b_i.$$
Up to canonical identification, this is the same as (\ref{MathItem_5.3}) over $SM$ and thus indeed any element of $\mathcal{R}^k \left(SM;N\otimes\pi ^\ast\mathcal{E} \right)$ extends continuously to ${}^0S^\ast\overline{M}$ as an element of $C^k \left({}^0S^\ast\overline{M};{}^0N\otimes\pi _0^\ast\mathcal{E} \right)$.

\end{remark}

\begin{example} \label{MathItem_5.7}  \normalfont Building off of the previous remark, one can construct examples of $u\in\mathcal{R}^k \left(SM,\ldots \right)$ that are not in $\mathcal{R}^{k+1} \left(SM,\ldots \right)$ where ``..." is $\pi ^\ast\mathcal{E}$ or $N\otimes\pi ^\ast\mathcal{E}$ as follows.
One can consider sections of the form $u=\rho^{k+\sfrac{ \left(n+1 \right)}{2}}\widetilde{u}$ where $\widetilde{u}\in C^\infty \left({}^0S^\ast\overline{M};\ldots \right)$ and $\widetilde{u}\neq0$ everywhere on $\partial {}^0S^\ast\overline{M}$.
In this case, for any smooth $V_1,\ldots,V_k\in C^\infty \left({}^0S^\ast\overline{M};T{}^0S^\ast\overline{M} \right)$, in the coordinates and frames of Definition \ref{MathItem_5.1} it will hold that $V_r\ldots V_1u^j$ is in $\rho^{k-r+\sfrac{ \left(n+1 \right)}{2}}C^\infty \left({}^0S^\ast\overline{M};\ldots \right)$ and hence in $\rho^{\sfrac{ \left(n+1 \right)}{2}}L^\infty \left({}^0S^\ast\overline{M};\ldots \right)$ for $r=1,\ldots,k$.
However, if one takes $V_1,\ldots,V_{k+1}=\sfrac{\partial }{\partial \rho}$, then $V_1\ldots V_{k+1}u^j=\rho^{-1+\sfrac{ \left(n+1 \right)}{2}}{\hat{u}}^j$ for some $\hat{u}\in C^\infty \left({}^0S^\ast\overline{M};\ldots \right)$ with $\hat{u}\neq0$ everywhere on $\partial {}^0S^\ast\overline{M}$ and hence it will not hold that $V_1\ldots V_{k+1}u^j\in\rho^{\sfrac{ \left(n+1 \right)}{2}}L^\infty\left[\pi ^{-1}\left[K\right]\right]$.
Thus $u\notin\mathcal{R}^{k+1}$.
\end{example}

\begin{remark} \label{MathItem_5.8}  \normalfont We make a remark regarding integrability.
In Definition \ref{MathItem_5.1} we have that each $\left|\pi ^\ast\sfrac{\partial }{\partial  x^i}\right|_N=\left|\sfrac{\partial }{\partial  x^i}\right|_g$ is $\rho^{-1}$ times something smooth on $\overline{M}$ (since $g$ is $\rho^{-2}$ times something smooth on $\overline{M}$).
Furthermore, each $\left|\pi ^\ast b_j\right|_{\pi ^\ast\mathcal{E}}=\left|b_j\right|_\mathcal{E}$ is bounded over $K$ since $\mathcal{E}$ is a smooth bundle over $\overline{M}$ by our standing assumption.
Hence, it follows that both $\left|u\right|_{\pi ^\ast\mathcal{E}},\left|w\right|_{N\otimes\pi ^\ast\mathcal{E}}\in\rho^{\sfrac{ \left(n+1 \right)}{2}}L^\infty\left[K\right]$.
Since $\overline{M}$ is compact and hence can be covered by a finite collection of such sets $K$, it follows that $\mathcal{R}^0\subseteq L^2$ by the comment made at the end of Section \ref{Section_4.4} above.
Furthermore, it follows straight from the definition that $\mathcal{R}^k\subseteq\mathcal{R}^{k^\prime}$ if $k\geq k^\prime$.
Hence, $\mathcal{R}^k\subseteq L^2$ for all $k\geq0$.
\end{remark}

The main reason for considering these spaces $\mathcal{R}^k$ is that, in addition to being subsets of $L^2$, all differential operator that we use in this paper have the mapping properties $\mathcal{R}^k\rightarrow\mathcal{R}^{k-1}$.
\subsection{Pestov Identity} \label{Section_5.2} 

In this section we prove the following version of the \textbf{Pestov Identity} with a connection on asymptotically hyperbolic (AH) spaces, which appears as Proposition 3.3 in \cite{Bibitem_18} in the setting of compact manifolds.
It will be used in the Fourier analysis study of transport equations.
Note the ``unitary" assumption on $\nabla^\mathcal{E}$ in the following theorem's statement.

\begin{theorem} \label{MathItem_5.9}  Suppose that $ \left(M\subseteq\overline{M},g \right)$ is an asymptotically hyperbolic space, $ \left(\mathcal{E},\langle \cdot,\cdot\rangle _\mathcal{E} \right)$ is a smooth complex Hermitian vector bundle over $\overline{M}$, and that $\nabla^\mathcal{E}$ is a smooth unitary connection in $\mathcal{E}$.
Suppose also that $u\in\mathcal{R}^2 \left(SM;\pi ^\ast\mathcal{E} \right)$.
Then\footnote{ Here we are implicitly restricting to the interior so that we may apply the differential operators involved and integrate.}
$$\left\lVert {\buildrel\mathrm{v}\over\nabla}{}^{\pi ^\ast\mathcal{E}}\mathbb{X}u\right\rVert _{L^2}^2=\left\lVert \mathbb{X}{\buildrel\mathrm{v}\over\nabla}{}^{\pi ^\ast\mathcal{E}}u\right\rVert _{L^2}^2-\langle R{\buildrel\mathrm{v}\over\nabla}{}^{\pi ^\ast\mathcal{E}}u,{\buildrel\mathrm{v}\over\nabla}{}^{\pi ^\ast\mathcal{E}}u\rangle _{L^2}-\langle F^\mathcal{E}u,{\buildrel\mathrm{v}\over\nabla}{}^{\pi ^\ast\mathcal{E}}u\rangle _{L^2}+n\left\lVert \mathbb{X}u\right\rVert _{L^2}^2$$
where $L^2$ stands for $L^2 \left(SM;N\otimes\pi ^\ast\mathcal{E} \right)$ in the first four quantities and $L^2 \left(SM;\pi ^\ast\mathcal{E} \right)$ in the last one.
\end{theorem}

Intuitively speaking, the above Pestov identity studies how the ``energy" (i.e. $L^2$-norm squared) changes when one switches the order of ${\buildrel\mathrm{v}\over\nabla}{}^{\pi ^\ast\mathcal{E}}$ and $\mathbb{X}$.
We require that $u$ is in $\mathcal{R}^2 \left(SM;\pi ^\ast\mathcal{E} \right)$ to ensure that all of the $L^2$ norms and inner products in the above equation make sense.
The order of the space $\mathcal{R}^2$ is chosen to be ``2" because both in the statement of the theorem and its proof we will not be applying more than two smooth vector fields over the 0-cosphere bundle to the components of $u$ at any one time.
The theorem is proved by simply starting with $\lVert {\buildrel\mathrm{v}\over\nabla}{}^{\pi ^\ast\mathcal{E}}\mathbb{X}u\rVert _{L^2}^2$ and then applying $L^2$-adjoint relations and commutator formulas until one arrives at $\lVert \mathbb{X}{\buildrel\mathrm{v}\over\nabla}{}^{\pi ^\ast\mathcal{E}}u\rVert _{L^2}^2$.
The following lemma provides us with the required set of adjoint relations.

\begin{lemma} \label{MathItem_5.10}  Suppose that $\nabla^\mathcal{E}$ is unitary.
The following are true, where $m,m^\prime\geq1$ are integers and all $L^2$ stand for appropriate $L^2 \left(SM;\ldots \right)$ spaces.

	\begin{enumerate} \item If $u\in\mathcal{R}^m \left(SM;\pi ^\ast\mathcal{E} \right)$, then $\mathbb{X}u\in\mathcal{R}^{m-1} \left(SM;\pi ^\ast\mathcal{E} \right)$.
Furthermore, if $w\in\mathcal{R}^{m^\prime} \left(SM;\pi ^\ast\mathcal{E} \right)$, then

$$\langle \mathbb{X}u,w\rangle _{L^2}=-\langle u,\mathbb{X}w\rangle _{L^2}.$$
	\item If $u\in\mathcal{R}^m \left(SM;N\otimes\pi ^\ast\mathcal{E} \right)$, then $\mathbb{X}u\in\mathcal{R}^{m-1} \left(SM;N\otimes\pi ^\ast\mathcal{E} \right)$.
Furthermore, if $w\in\mathcal{R}^{m^\prime} \left(SM;N\otimes\pi ^\ast\mathcal{E} \right)$, then

$$\langle \mathbb{X}u,w\rangle _{L^2}=-\langle u,\mathbb{X}w\rangle _{L^2}.$$
	\item If $u\in\mathcal{R}^m \left(SM;\pi ^\ast\mathcal{E} \right)$ and $w\in\mathcal{R}^{m^\prime} \left(SM;N\otimes\pi ^\ast\mathcal{E} \right)$, then ${\buildrel\mathrm{v}\over\nabla}{}^{\pi ^\ast\mathcal{E}}u\in\mathcal{R}^{m-1} \left(SM;N\otimes\pi ^\ast\mathcal{E} \right)$ and ${\buildrel\mathrm{v}\over{\mathrm{div}}}{}^{\pi ^\ast\mathcal{E}}w\in\mathcal{R}^{m^\prime-1} \left(SM;\pi ^\ast\mathcal{E} \right)$.
Furthermore,

$$\langle {\buildrel\mathrm{v}\over\nabla}{}^{\pi ^\ast\mathcal{E}}u,w\rangle _{L^2}=-\langle u,{\buildrel\mathrm{v}\over{\mathrm{div}}}{}^{\pi ^\ast\mathcal{E}}w\rangle _{L^2}.$$\end{enumerate}
\end{lemma}

In other words, $\mathbb{X}$, ${\buildrel\mathrm{v}\over\nabla}{}^{\pi ^\ast\mathcal{E}}$, and ${\buildrel\mathrm{v}\over{\mathrm{div}}}{}^{\pi ^\ast\mathcal{E}}$ map $\mathcal{R}^m\rightarrow\mathcal{R}^{m-1}$ and their well-known adjoint relations are also satisfied on AH spaces as well.
To prove the above lemma, we will use the following compactly supported version of it:

\begin{lemma} \label{MathItem_5.11}  Suppose that $\nabla^\mathcal{E}$ is unitary.
The following are true.

	\begin{enumerate} \item If $u,w\in C^\infty \left(SM;\pi ^\ast\mathcal{E} \right)$ are such that at least one of them is compactly supported in the interior $M$, then

$$\langle \mathbb{X}u,w\rangle _{L^2}=-\langle u,\mathbb{X}w\rangle _{L^2}.$$
	\item If $u,w\in C^\infty \left(SM;N\otimes\pi ^\ast\mathcal{E} \right)$ are such that at least one of them is compactly supported in the interior $M$, then

$$\langle \mathbb{X}u,w\rangle _{L^2}=-\langle u,\mathbb{X}w\rangle _{L^2}.$$\end{enumerate}
\end{lemma}

\noindent\textbf{Proof:} To prove point 1), let $ \left(x^i \right)$ be coordinates of $M$, let $ \left(b_i \right)$ be a frame for $\mathcal{E}$ over these coordinates' domain, and consider the coordinates $v^i\sfrac{\partial }{\partial  x^i}\mapsto \left(x^i,v^i \right)$ of $TM$.
Let ${}^\mathcal{E}\Gamma_{ij}^k$ denote the connection symbols of $\nabla^\mathcal{E}$ with respect to $ \left(\sfrac{\partial }{\partial  x^i} \right)$ and $ \left(b_k \right)$.
Now, suppose first that $u$ or $w$ is compactly supported over our coordinates' domain and write $u=u^k\pi ^\ast b_k$ and $w=w^k\pi ^\ast b_k$.
Then
$$\mathbb{X}u=\left[X \left(u^k \right)+{}^\mathcal{E}\Gamma_{ij}^kv^iu^j\right]\pi ^\ast b_k.$$
For convenience, assume that $ \left(b_i \right)$ is orthonormal so that we may write
$$\langle \mathbb{X}u,w\rangle _{L^2}=\sum_{k=1}^{d}\int_{SM}{\left[X \left(u^k \right)w^k+{}^\mathcal{E}\Gamma_{ij}^kv^iu^jw^k\right]dv_{SM}}.$$
In Appendix A of \cite{Bibitem_41} the authors prove that the $L^2$ adjoint of $X : C^\infty \left(SM \right)\rightarrow C^\infty \left(SM \right)$ is $-X$.
Furthermore, since $\nabla^\mathcal{E}$ is unitary it follows that the connection symbols are anti-symmetric: ${}^\mathcal{E}\Gamma_{ij}^k=-{}^\mathcal{E}\Gamma_{ik}^j$.
Applying these identities to the right-hand side above gives $-\langle u,\mathbb{X}w\rangle _{L^2}$.
Point 1) then follows by a partition of unity argument.
Point 2) is proved similarly where instead one uses the fact that the $L^2$ adjoint of the operator $X : C^\infty \left(SM;N \right)\rightarrow C^\infty \left(SM;N \right)$ defined in (\ref{MathItem_4.22}) is also $-X$, which is proved in Appendix A of \cite{Bibitem_41}.

\begin{flushright}$\blacksquare$\end{flushright}

\noindent\textit{Proof of Lemma \ref{MathItem_5.10} part 1):}

Let $u$ be as described in part 1).
Suppose $\rho : \overline{M}\rightarrow\left[0,\infty \right)$ is a geodesic boundary defining function as defined in Section \ref{Section_2.2} and let $\varepsilon>0$ be as described there.
Let $ \left(y^1,\ldots,y^n \right)$ be coordinates of $\partial \overline{M}$ and let $ \left(\rho,y^\mu \right)= \left(x^i \right)$ be the asymptotic boundary normal coordinates of $\overline{M}$ that they generate on $\left[0,\varepsilon \right)\times\dom{ \left(y^\mu \right)}$ as in (\ref{MathItem_2.3}).
Let $ \left(b_i \right)$ be a frame for $\mathcal{E}$ over these coordinates' domain.
Consider the coordinates $v^i\sfrac{\partial }{\partial  x^i}\mapsto \left(x^i,v^i \right)$ of $TM$.
Let ${}^\mathcal{E}\Gamma_{ij}^k$ denote the connection symbols of $\nabla^\mathcal{E}$ with respect to $ \left(\sfrac{\partial }{\partial  x^i} \right)$ and $ \left(b_k \right)$ and let $\Gamma_{ij}^k$ denote the Christoffel symbols with respect to $ \left(\sfrac{\partial }{\partial  x^i} \right)$.
Consider also the coordinates ${\overline{\eta}}_i\sfrac{dx^i}{\rho}\rightarrow \left(x^i,{\overline{\eta}}_i \right)$ of ${}^0T^\ast\overline{M}$ and observe that the canonical identification $H^\ast\circ\flat  : TM\rightarrow\left.{}^0T^\ast\overline{M}\right|_M$ is given by $v^i=\rho^{-1}g^{ii^\prime}{\overline{\eta}}_{i^\prime}$.
Recall that $g_{ij}$ and $g^{ij}$ are respectively $\rho^{-2}$ and $\rho^2$ times something smooth on $\overline{M}$.

Writing $u=u^k\pi ^\ast b_k$, we have that
\begin{equation} \label{MathItem_5.12} \mathbb{X}u=\left[X \left(u^k \right)+{}^\mathcal{E}\Gamma_{ij}^kv^iu^j\right]\pi ^\ast b_k.\end{equation}
Pulling $v^i$ to ${}^0S^\ast\overline{M}$ via the canonical identification gives $\rho^{-1}g^{ii^\prime}{\overline{\eta}}_{i^\prime}$, which is smooth over ${}^0S^\ast\overline{M}$.
The terms ${}^\mathcal{E}\Gamma_{ij}^k$ are smooth over $\overline{M}$ by our standing assumption and hence on ${}^0S^\ast\overline{M}$ when lifted.
Hence it quickly follows that the term ${}^\mathcal{E}\Gamma_{ij}^kv^iu^j\pi ^\ast b_k$ on the right-hand side of the above equation is in $\mathcal{R}^m\subseteq\mathcal{R}^{m-1}$.
So, let us take a look at the other term: $X \left(u^k \right)\pi ^\ast b_k$.

Consider the coordinates $\xi_idx^i\mapsto \left(x^i,\xi_i \right)$ of $T^\ast\overline{M}$.
In (2.3) of \cite{Bibitem_15} the authors write out an explicit equation for $X$ over $T^\ast M$ in asymptotic boundary normal coordinates (recall the convention about Greek and Latin indices):
$$X=\rho^2\xi_0\frac{\partial }{\partial \rho}+\rho^2h^{\mu\nu}\xi_\mu\frac{\partial }{\partial  y^\nu}-\left[\rho \left(\xi_0^2+\left|\xi\right|_h^2 \right)+\frac{1}{2}\rho^2\partial _\rho\left|\xi\right|_h^2\right]\frac{\partial }{\partial \xi_0}-\frac{1}{2}\rho^2\partial _{y^k}\left|\xi\right|_h^2\frac{\partial }{\partial \xi_k}$$
where $\left|\xi\right|_h^2=h^{\mu\nu}\xi_\mu\xi_\nu$.
Since canonical identification is given by ${\overline{\eta}}_i=\rho\xi_i$, it is a quick calculation to show that pushing $X$ to ${}^0S^\ast\overline{M}$ gives
\begin{equation} \label{MathItem_5.13} X=\rho{\overline{\eta}}_0\frac{\partial }{\partial \rho}+\rho h^{\mu\nu}{\overline{\eta}}_\mu\frac{\partial }{\partial  y^\nu}-\left[\left|\overline{\eta}\right|_h^2+\frac{1}{2}\rho\partial _\rho\left|\overline{\eta}\right|_h^2\right]\frac{\partial }{\partial {\overline{\eta}}_0}+\left[{\overline{\eta}}_0{\overline{\eta}}_\mu-\frac{1}{2}\rho\partial _{y^\mu}\left|\overline{\eta}\right|_h^2\right]\frac{\partial }{\partial {\overline{\eta}}_\mu}.\end{equation}
In particular, $X$ extends to be a smooth vector field on all of ${}^0S^\ast\overline{M}$ and hence $X \left(u^k \right)\pi ^\ast b_k$ is in $\mathcal{R}^{m-1}$.
Thus by (\ref{MathItem_5.12}) we have that indeed $\mathbb{X}u\in\mathcal{R}^{m-1}$.

Now suppose that $w$ is as in 1).
We will prove the equality in 1) by multiplying $w$ by a compactly supported (smooth) bump function, use Lemma \ref{MathItem_5.11} above, and then let the support of the bump function go out to infinity.
To construct the suitable family of bump functions, let $f_1 : \left[0,\infty \right)\rightarrow\left[0,\infty \right)$ be a smooth function that is identically zero on $\left[0,\sfrac{1}{2}\right]$, increasing on $\left[\sfrac{1}{2},1\right]$, and then identically one on $\left[1,\infty \right)$ (see Lemma 2.21 in \cite{Bibitem_25} for an explicit construction).
For any $\delta>0$, let $f_\delta : \left[0,\infty \right)\rightarrow\left[0,\infty \right)$ denote the function $f_\delta \left(x \right)=f \left(\sfrac{x}{\delta} \right)$.
Finally, for $\delta<\varepsilon$ let $\phi_\delta : \overline{M}\rightarrow\left[0,\infty \right)$ denote the one parameter family of bump functions given by
$$\phi_\delta \left(x \right)=\begin{cases}f_\delta\circ\rho \left(x \right)&\rho \left(x \right)<\delta\\1&\mathrm{otherwise}\\\end{cases}.$$
By Lemma \ref{MathItem_5.11} we have that
$$\langle \mathbb{X}u,\phi_\delta w\rangle _{L^2}=-\langle u,\mathbb{X} \left(\phi_\delta w \right)\rangle _{L^2}$$
since $\phi_\delta w$ is compactly supported.
Applying the product rule on the right-hand side gives
\begin{equation} \label{MathItem_5.14} \langle \mathbb{X}u,\phi_\delta w\rangle _{L^2}=-\langle u,\phi_\delta\mathbb{X}w\rangle _{L^2}-\langle u,X \left(\phi_\delta \right)w\rangle _{L^2}.\end{equation}
We now let $\delta\rightarrow0^+$ and show that this equation tends to the equality in 1).
We have that $u\in\mathcal{R}^m$, $w\in\mathcal{R}^{m^\prime}$, $\mathbb{X}u\in\mathcal{R}^{m-1}$, $\mathbb{X}w\in\mathcal{R}^{m^\prime-1}$ and hence all are in $L^2 \left(SM,\ldots \right)$.
Thus both $\langle \mathbb{X}u,w\rangle _{\pi ^\ast\mathcal{E}}$ and $\langle u,\mathbb{X}w\rangle _{\pi ^\ast\mathcal{E}}$ are in $L^1 \left(SM \right)$.
Next, differentiating in $\delta$ demonstrates that $\phi_\delta$ monotonely increases to the identically one function as $\delta\rightarrow0^+$.
Hence by the dominated convergence theorem, we get that the first two terms in (\ref{MathItem_5.14}) tend to $\langle \mathbb{X}u,w\rangle _{L^2}$ and $-\langle u,\mathbb{X}w\rangle _{L^2}$ respectively as $\delta\rightarrow0^+$.

Hence we will have proved 1) if we can show that the third term in (\ref{MathItem_5.14}) tends to zero as $\delta\rightarrow0^+$.
This will follow if we show that for any compact set $K\subseteq\overline{M}$ contained in the domain of some interior coordinates $ \left(x^i \right)$ of $\overline{M}$ or boundary coordinates $ \left(\rho,y^\mu \right)= \left(x^i \right)$ as above,
$$\int_{\pi ^{-1}\left[K\right]}\langle u,X \left(\phi_\delta \right)w\rangle _{\pi ^\ast\mathcal{E}}\longrightarrow0\ \ \ \ \ \mathrm{as}  \ \ \delta\rightarrow0^+.$$
If $K$ is contained in the domain of interior coordinates, then this follows immediately since $\phi_\delta\equiv1$ on $K$ for sufficiently small $\delta>0$.
So suppose that $K$ is contained in the domain of our boundary coordinates $ \left(\rho,y^\mu \right)= \left(x^i \right)$.
Writing the above integral in these coordinates as in Section \ref{Section_4.4} gives (here $d\hat{x}=dx^0\ldots dx^n$)
$$\left|\int_{K}{\int_{S_xM}{\langle u,X \left(\phi_\delta \right)w\rangle _{\pi ^\ast\mathcal{E}}dS_x \left(v \right)}\sqrt{\det{g}}d\hat{x}}\right|\le\sup_ K{\left|\langle u,w\rangle _{\pi ^\ast\mathcal{E}}\sqrt{\det{g}}\right|}\int_{K}{\int_{S_xM}{\left|X \left(\phi_\delta \right)\right|dS_x \left(v \right)}d\hat{x}}.$$
The $\sup{\ldots}$ is finite because $\left|u\right|_{\pi ^\ast\mathcal{E}},\left|w\right|_{\pi ^\ast\mathcal{E}}\in\rho^{\sfrac{ \left(n+1 \right)}{2}}L^\infty\left[\pi ^{-1}\left[K\right]\right]$ and $\sqrt{\det{g}}\in\rho^{- \left(n+1 \right)}C^\infty \left(\overline{M} \right)$.
Now, the explicit equation for $X$ in coordinates of $TM$ (e.g. see page 104 in \cite{Bibitem_26}) gives that $X \left(\phi_\delta \right)=v^0f_\delta^\prime \left(\rho \right)$, which we note is supported in $\left\{\rho\le\delta\right\}$.
Since $g=\sfrac{ \left(d\rho^2+h_{\mu\nu}dy^\mu d y^\nu \right)}{\rho^2}$ and $\left|v\right|_g=1$, we have that $\left|v^0\right|\le\rho$.
Letting $K_y$ denote the (compact) ordinary projection of $K$ onto the set $\left\{ \left(0,y^\mu \right)\right\}$ in our coordinates, we can bound
$$\int_{K}{\int_{S_xM}{\left|X \left(\phi_\delta \right)\right|dS_x \left(v \right)}d\hat{x}}\le\omega_n\int_{K_y}{\int_{0}^{\delta}{\delta f_\delta^\prime \left(\rho \right)d\rho}dy}=\delta\int_{K_y} d y\longrightarrow0\ \ \ \ \ \mathrm{as}  \ \ \delta\rightarrow0^+,$$
where $\omega_n$ denotes the surface area of the Euclidean $n$-sphere.
Hence the third term in (\ref{MathItem_5.14}) indeed tends to zero as $\delta\rightarrow0^+$.

\begin{flushright}$\blacksquare$\end{flushright}

\noindent\textit{Proof of Lemma \ref{MathItem_5.10} parts 2), 3):}

Let us begin with proving 2).
Let $u$ be as described there.
We keep working in the same coordinates that we used in part 1) above.
Writing $u=\rho u^{ij}\pi ^\ast\sfrac{\partial }{\partial  x^i}\otimes\pi ^\ast b_j$, we have that
$$\mathbb{X}u$$
\begin{equation} \label{MathItem_5.15} =X \left(\rho u^{kj} \right)\pi ^\ast\frac{\partial }{\partial  x^k}\otimes\pi ^\ast b_j\end{equation}
\begin{equation} \label{MathItem_5.16} +\Gamma_{i^\prime i}^kv^{i^\prime}\rho u^{ij}\pi ^\ast\frac{\partial }{\partial  x^k}\otimes\pi ^\ast b_j\end{equation}
\begin{equation} \label{MathItem_5.17} +{}^\mathcal{E}\Gamma_{i^\prime j}^kv^{i^\prime}\rho u^{ij}\pi ^\ast\frac{\partial }{\partial  x^i}\otimes\pi ^\ast b_k.\end{equation}
In part 1) we observed that $v^{i^\prime}$ is $\rho$ times something smooth on ${}^0S^\ast\overline{M}$ and hence it follows that the term (\ref{MathItem_5.17}) is in $\mathcal{R}^m \left(SM;N\otimes\pi ^\ast\mathcal{E} \right)\subseteq\mathcal{R}^{m-1} \left(SM;N\otimes\pi ^\ast\mathcal{E} \right)$.
Next, we have that since $g$ is $\rho^{-2}$ times a smooth metric on $\overline{M}$, the conformal transformation law of Christoffel symbols (e.g. see Proposition 7.29 in \cite{Bibitem_26}) give that $\Gamma_{ij}^k$ are $\rho^{-1}$ times something smooth on $\overline{M}$.
Thus the term (\ref{MathItem_5.16}) is also in $\mathcal{R}^m\subseteq\mathcal{R}^{m-1}$.
Finally, from (\ref{MathItem_5.13}) we have that $X \left(\rho \right)$ is $\rho$ times something smooth on ${}^0S^\ast\overline{M}$ and thus it follows from the product rule that the term (\ref{MathItem_5.15}) is in $\mathcal{R}^{m-1}$.
Hence indeed $\mathbb{X}u\in\mathcal{R}^{m-1}$.
The equality in 2) follows the same way that we proved the equality in 1).

Finally, let us prove 3).
Let $u=u^j\pi ^\ast b_j$ and $w=\rho w^{ij}\pi ^\ast\sfrac{\partial }{\partial  x^i}\otimes\pi ^\ast b_j$ be as described there.
In the proof of Lemma 3.2 in \cite{Bibitem_18} the authors give an equation for the vertical derivative\footnote{ We remark that in their work they write ``${\buildrel\mathrm{v}\over\nabla}{}^\mathcal{E}$" for what we denote by ``${\buildrel\mathrm{v}\over\nabla}{}^{\pi ^\ast\mathcal{E}}$."} of $u$ in terms of an operator ``${\buildrel\mathrm{v}\over\nabla}$" for which an explicit equation is given on page 350 of \cite{Bibitem_41}.
As the authors do in \cite{Bibitem_18}, we assume that $ \left(b_i \right)$ is orthonormal so that we may use their formula to write that
\begin{equation} \label{MathItem_5.18} {\buildrel\mathrm{v}\over\nabla}{}^{\pi ^\ast\mathcal{E}}u\ \mathrm{over}\ SM=\partial ^iu^j\pi ^\ast\frac{\partial }{\partial  x^i}\otimes\pi ^\ast b_j,\end{equation}
where for any $f\in C^\infty \left(SM \right)$
$$\partial _if:=\left.\left[\frac{\partial }{\partial  v^i} \left(f\circ p \right)\right]\right|_{SM},$$
$$\partial ^if:=g^{ii^\prime}\partial _{i^\prime}f,$$
where $p : TM\setminus\left\{0\right\}\rightarrow SM$ is the radial projection map $v\mapsto\sfrac{v}{\left|v\right|_g}$ over the tangent bundle minus the zero section.
Alternatively, the equality (\ref{MathItem_5.18}) also follows from (\ref{MathItem_4.20}).
Since the canonical identification$\left.\ {}^0S^\ast\overline{M}\right|_M\cong SM$ is given by ${\overline{\eta}}_i=\rho g_{ii^\prime}v^{i^\prime}$, pushing $\sfrac{\partial }{\partial  v^i}$ to ${}^0S^\ast\overline{M}$ gives
$$\frac{\partial }{\partial  v^i}\longmapsto\rho g_{ii^\prime}\frac{\partial }{\partial {\overline{\eta}}_{i^\prime}}.$$
Thus $\sfrac{\partial }{\partial  v^i}$ extends to $\rho^{-1}$ times a smooth vector field over ${}^0S^\ast\overline{M}$ and hence so does $\partial _i$.
Hence $\partial ^i$ is $\rho$ times a smooth vector field over ${}^0S^\ast\overline{M}$.
Thus by (\ref{MathItem_5.18}) we have that the vertical derivative of $u$ is indeed in $\mathcal{R}^{m-1} \left(SM;N\otimes\pi ^\ast\mathcal{E} \right)$.

Next let us take a look at $w$.
Using (\ref{MathItem_5.18}) above, a straightforward generalization of the derivation of the equation for ``${\buildrel\mathrm{v}\over{\mathrm{div}}}\ Z$" given on page 352 of \cite{Bibitem_41} gives
\begin{equation} \label{MathItem_5.19} {\buildrel\mathrm{v}\over{\mathrm{div}}}{}^{\pi ^\ast\mathcal{E}}w=\partial _i \left(\rho w^{ij} \right)\pi ^\ast b_j.\end{equation}
Alternatively, this equality also follows from (\ref{MathItem_4.19}) (the $w^{ij}$ there is our $\rho w^{ij}$ here).
Since $\partial _i$ only involves derivatives in $v^i$, we can pull $\rho$ out of the derivative on the right-hand side.
Since $\partial _i$ is $\rho^{-1}$ times a smooth vector field over ${}^0S^\ast\overline{M}$, it follows that this is indeed in $\mathcal{R}^{m^\prime-1} \left(SM;\pi ^\ast\mathcal{E} \right)$.

The equality in 3) follow essentially the same way we proved the equality in 1).
An example of a minor change that is needed is that the analog of (\ref{MathItem_5.14}) will be
$$\langle {\buildrel\mathrm{v}\over\nabla}{}^{\pi ^\ast\mathcal{E}}u,\phi_\delta w\rangle _{L^2}=-\langle u,\phi_\delta{\buildrel\mathrm{v}\over{\mathrm{div}}}{}^{\pi ^\ast\mathcal{E}}w\rangle _{L^2},$$
which we note does not have an analogous ``third term" as in (\ref{MathItem_5.14}) because $\phi_\delta$ only depends on position and thus is not affected by the vertical divergence.
From here one proceeds as before.

\begin{flushright}$\blacksquare$\end{flushright}

For use in Section \ref{Section_5.3} below, we record the $\mathcal{R}^m$-mapping property of the horizontal derivative as well.

\begin{lemma} \label{MathItem_5.20}  Suppose that $\nabla^\mathcal{E}$ is unitary.
If $u\in\mathcal{R}^m \left(SM;\pi ^\ast\mathcal{E} \right)$ with $m\geq1$, then ${\buildrel\mathrm{h}\over\nabla}{}^{\pi ^\ast\mathcal{E}}u\in\mathcal{R}^{m-1} \left(SM;N\otimes\pi ^\ast\ \mathcal{E} \right)$.
\end{lemma}

\noindent\textbf{Proof:} Take any such $u$.
Let $ \left(\rho,y^\mu \right)= \left(x^i \right)$, $ \left(b_i \right)$, ${}^\mathcal{E}\Gamma_{ij}^k$, and $\Gamma_{ij}^k$ be as in the beginning of the proof of Lemma \ref{MathItem_5.10} part 1) and consider the coordinates $v^i\sfrac{\partial }{\partial  x^i}\mapsto \left(x^i,v^i \right)$ of $TM$ and ${\overline{\eta}}_i\sfrac{dx^i}{\rho}\rightarrow \left(x^i,{\overline{\eta}}_i \right)$ of ${}^0T^\ast\overline{M}$ described there as well.

We write $u=u^j\pi ^\ast b_j$.
By the equations for the horizontal and vertical derivatives given in the proof of Lemma 3.2 in \cite{Bibitem_18} and on page 350 of \cite{Bibitem_41},
$${\buildrel\mathrm{h}\over\nabla}{}^{\pi ^\ast\mathcal{E}}u$$
\begin{equation} \label{MathItem_5.21} = \left(\delta^iu^j- \left(v^k\delta_ku^j \right)v^i \right)\pi ^\ast\frac{\partial }{\partial  x^i}\otimes\pi ^\ast b_j\end{equation}
\begin{equation} \label{MathItem_5.22} +u^l{\buildrel\mathrm{v}\over\nabla}{}^{\pi ^\ast\mathcal{E}} \left({}^\mathcal{E}\Gamma_{kl}^jv^k\pi ^\ast b_j \right)\end{equation}
where for any $f\in C^\infty \left(SM \right)$
\begin{equation} \label{MathItem_5.23} \delta_if:=\left.\left[ \left(\frac{\partial }{\partial  x^i}-\Gamma_{ij}^kv^j\frac{\partial }{\partial  v^k} \right) \left(f\circ p \right)\right]\right|_{SM},\end{equation}
\begin{equation} \label{MathItem_5.24} \delta^if:=g^{ii^\prime}\delta_{i^\prime}f,\end{equation}
where $p : TM\setminus\left\{0\right\}\rightarrow SM$ is the radial projection map $v\mapsto\sfrac{v}{\left|v\right|_g}$.

Let us start by taking a look at the term (\ref{MathItem_5.22}).
Recall our standing assumption that each ${}^\mathcal{E}\Gamma_{kl}^j$ is smooth on $\overline{M}$.
Next, we observed in the proof of Lemma \ref{MathItem_5.10} part 1) that each $v^k$ is $\rho$ times something smooth over ${}^0S^\ast\overline{M}$.
By (\ref{MathItem_5.18}) the vertical derivative only involves derivatives in $v^i$ and so we can pull the just mentioned factor of $\rho$ out of the vertical derivative in (\ref{MathItem_5.22}).
Thus the term (\ref{MathItem_5.22}) is equal to $\rho u^l$ times $\partial ^i$ of something smooth on ${}^0S^\ast\overline{M}$ times $\pi ^\ast\sfrac{\partial }{\partial  x^i}\otimes\pi ^\ast b_j$.
We observed in the proof of Lemma \ref{MathItem_5.10} parts 2) and 3) that $\partial ^i$ is a smooth vector field over ${}^0S^\ast\overline{M}$ and thus $\partial ^i$ of something smooth on ${}^0S^\ast\overline{M}$ is again smooth on ${}^0S^\ast\overline{M}$.
From here it follows that the term (\ref{MathItem_5.22}) is in $\mathcal{R}^m\subseteq\mathcal{R}^{m-1}$.

Finally let us take a look at the term (\ref{MathItem_5.21}).
Recall that $g_{ij}$ and $g^{ij}$ are respectively $\rho^{-2}$ and $\rho^2$ times something smooth on $\overline{M}$.
We observed in the proof of Lemma \ref{MathItem_5.10} parts 2) and 3) that each $\Gamma_{ij}^kv^j$ is smooth on ${}^0S^\ast\overline{M}$.
Next, canonical identification is given by ${\overline{\eta}}_i=\rho g_{ij}v^j$, from which it is a quick computation to show that the differential of this canonical identification takes (recall that $x^0=\rho$)
$$\frac{\partial }{\partial  x^0}\longmapsto\frac{\partial }{\partial  x^0}+ \left(\rho^{-1}{\overline{\eta}}_i+\frac{\partial  g_{ij}}{\partial \rho}g^{jj^\prime}{\overline{\eta}}_{j^\prime} \right)\frac{\partial }{\partial {\overline{\eta}}_i}$$
$$\frac{\partial }{\partial  x^\lambda}\longmapsto\frac{\partial }{\partial  x^\lambda}+ \left(\frac{\partial  g_{ij}}{\partial  x^\lambda}g^{jj^\prime}{\overline{\eta}}_{j^\prime} \right)\frac{\partial }{\partial {\overline{\eta}}_i}\ \ \ \ \ \ \ \ \mathrm{for} \ \ \lambda=1,\ldots,n,$$
$$\frac{\partial }{\partial  v^i}\longmapsto\rho g_{ii^\prime}\frac{\partial }{\partial {\overline{\eta}}_{i^\prime}}.$$
The important observation is that these are all $\rho^{-1}$ times vector fields that are smooth over ${}^0S^\ast\overline{M}$.
Thus by (\ref{MathItem_5.23}) and (\ref{MathItem_5.24}), the $\delta_i$ and $\delta^i$ are respectively $\rho^{-1}$ and $\rho$ times vector fields that are smooth over ${}^0S^\ast\overline{M}$.
Plugging all of these observations into (\ref{MathItem_5.21}) finally gives that indeed the horizontal derivative of $u$ is in $\mathcal{R}^{m-1}$.

\begin{flushright}$\blacksquare$\end{flushright}

Next we need the following lemma that tells us that the curvature operators have the mapping property $\mathcal{R}^m\rightarrow\mathcal{R}^m$.
Note that we do not need $\nabla^\mathcal{E}$ to be unitary for this lemma.

\begin{lemma} \label{MathItem_5.25}  Suppose that $u\in\mathcal{R}^m \left(SM;N\otimes\pi ^\ast\mathcal{E} \right)$ and $w\in\mathcal{R}^m \left(SM;\pi ^\ast\mathcal{E} \right)$ for integers $m\geq0$.
Then $Ru\in\mathcal{R}^m \left(SM;N\otimes\pi ^\ast\mathcal{E} \right)$ and $F^\mathcal{E}w\in\mathcal{R}^m \left(SM;N\otimes\pi ^\ast\mathcal{E} \right)$.
\end{lemma}

\noindent\textbf{Proof:} Let $ \left(\rho,y^\mu \right)= \left(x^i \right)$, $ \left(b_i \right)$, ${}^\mathcal{E}\Gamma_{ij}^k$, and $\Gamma_{ij}^k$ be as in the beginning of the proof of Lemma \ref{MathItem_5.10} part 1) and consider the coordinates $v^i\sfrac{\partial }{\partial  x^i}\mapsto \left(x^i,v^i \right)$ of $TM$ and ${\overline{\eta}}_i\sfrac{dx^i}{\rho}\rightarrow \left(x^i,{\overline{\eta}}_i \right)$ of ${}^0T^\ast\overline{M}$ described there as well.

We write $u=\rho u^{ij}\pi ^\ast\sfrac{\partial }{\partial  x^i}\otimes\pi ^\ast b_j$ and $w=w^j\pi ^\ast b_j$.
We have that
$$Ru= \left(\frac{\partial \Gamma_{j^\prime k}^l}{\partial  x^i}-\frac{\partial \Gamma_{ik}^l}{\partial  x^{j^\prime}}+\Gamma_{j^\prime k}^m\Gamma_{im}^l-\Gamma_{ik}^m\Gamma_{j^\prime m}^l \right)\rho u^{ij}v^{j^\prime}v^k\pi ^\ast\frac{\partial }{\partial  x^l}\otimes\pi ^\ast b_j.$$
In the proof of Lemma \ref{MathItem_5.10} part 1) we observed that each $v^i$ is $\rho$ times something smooth on ${}^0S^\ast\overline{M}$.
In the proof of Lemma \ref{MathItem_5.10} parts 2) and 3) we observed that each $\Gamma_{ij}^k$ is $\rho^{-1}$ times something smooth on $\overline{M}$ and hence its first partials $\sfrac{\partial \Gamma_{ij}^k}{\partial  x^l}$ are $\rho^{-2}$ times something smooth on $\overline{M}$.
From this it follows that indeed $Ru\in\mathcal{R}^m$.

Next, looking at (\ref{MathItem_4.24}) and (\ref{MathItem_4.25}) we have that
$$F^\mathcal{E}w=w^lg^{jj^\prime}f_{ij}{}^k{}_lv^i\pi ^\ast\frac{\partial }{\partial  x^{j^\prime}}\otimes\pi ^\ast b_k.$$
By (\ref{MathItem_4.23}), we have that each $f_{ij}{}^k{}_l$ is smooth on $\overline{M}$, and recall that each $g^{ij}$ is $\rho^2$ times something smooth on $\overline{M}$.
From here it follows that $F^\mathcal{E}w\in\mathcal{R}^m$ as well.

\begin{flushright}$\blacksquare$\end{flushright}

We need one final lemma that provides the needed commutator formulas to prove Theorem \ref{MathItem_5.9}.
The following lemma is Lemma 3.2 in \cite{Bibitem_18}, where one can also find a proof.

\begin{lemma} \label{MathItem_5.26}  The following are true, where $\left[\ldots,\ldots\right]$ denotes the commutator bracket.
\begin{equation} \label{MathItem_5.27} \left[\mathbb{X},{\buildrel\mathrm{v}\over\nabla}{}^{\pi ^\ast\mathcal{E}}\right]=-{\buildrel\mathrm{h}\over\nabla}{}^{\pi ^\ast\mathcal{E}},\end{equation}
\begin{equation} \label{MathItem_5.28} \left[\mathbb{X},{\buildrel\mathrm{h}\over\nabla}{}^{\pi ^\ast\mathcal{E}}\right]=R{\buildrel\mathrm{v}\over\nabla}{}^{\pi ^\ast\mathcal{E}}+F^\mathcal{E},\end{equation}
\begin{equation} \label{MathItem_5.29} {\buildrel\mathrm{h}\over{\mathrm{div}}}{}^{\pi ^\ast\mathcal{E}}{\buildrel\mathrm{v}\over\nabla}{}^{\pi ^\ast\mathcal{E}}-{\buildrel\mathrm{v}\over{\mathrm{div}}}{}^{\pi ^\ast\mathcal{E}}{\buildrel\mathrm{h}\over\nabla}{}^{\pi ^\ast\mathcal{E}}=n\mathbb{X},\end{equation}
\begin{equation} \label{MathItem_5.30} \left[\mathbb{X},{\buildrel\mathrm{v}\over{\mathrm{div}}}{}^{\pi ^\ast\mathcal{E}}\right]=-{\buildrel\mathrm{h}\over{\mathrm{div}}}{}^{\pi ^\ast\mathcal{E}}.\end{equation}
\end{lemma}

\noindent\textit{Proof of Theorem \ref{MathItem_5.9}:}

Let $u$ be as described in the theorem.
By Lemma \ref{MathItem_5.10} we have that
$$\langle {\buildrel\mathrm{v}\over\nabla}{}^{\pi ^\ast\mathcal{E}}\mathbb{X}u,{\buildrel\mathrm{v}\over\nabla}{}^{\pi ^\ast\mathcal{E}}\mathbb{X}u\rangle _{L^2}=\langle \mathbb{X}{\buildrel\mathrm{v}\over{\mathrm{div}}}{}^{\pi ^\ast\mathcal{E}}{\buildrel\mathrm{v}\over\nabla}{}^{\pi ^\ast\mathcal{E}}\mathbb{X}u,u\rangle _{L^2}.$$
We get that this is equal to (see right after for justifications)
$$\langle -{\buildrel\mathrm{h}\over{\mathrm{div}}}{}^{\pi ^\ast\mathcal{E}}{\buildrel\mathrm{v}\over\nabla}{}^{\pi ^\ast\mathcal{E}}\mathbb{X}u+{\buildrel\mathrm{v}\over{\mathrm{div}}}{}^{\pi ^\ast\mathcal{E}}\mathbb{X}{\buildrel\mathrm{v}\over\nabla}{}^{\pi ^\ast\mathcal{E}}\mathbb{X}u,u\rangle _{L^2},$$
$$=\langle -{\buildrel\mathrm{h}\over{\mathrm{div}}}{}^{\pi ^\ast\mathcal{E}}{\buildrel\mathrm{v}\over\nabla}{}^{\pi ^\ast\mathcal{E}}\mathbb{X}u+{\buildrel\mathrm{v}\over{\mathrm{div}}}{}^{\pi ^\ast\mathcal{E}}\mathbb{X}{\buildrel\mathrm{h}\over\nabla}{}^{\pi ^\ast\mathcal{E}}u+{\buildrel\mathrm{v}\over{\mathrm{div}}}{}^{\pi ^\ast\mathcal{E}}\mathbb{XX}{\buildrel\mathrm{v}\over\nabla}{}^{\pi ^\ast\mathcal{E}}u,u\rangle _{L^2},$$
$$=\langle -{\buildrel\mathrm{h}\over{\mathrm{div}}}{}^{\pi ^\ast\mathcal{E}}{\buildrel\mathrm{v}\over\nabla}{}^{\pi ^\ast\mathcal{E}}\mathbb{X}u+{\buildrel\mathrm{v}\over{\mathrm{div}}}{}^{\pi ^\ast\mathcal{E}} \left(R{\buildrel\mathrm{v}\over\nabla}{}^{\pi ^\ast\mathcal{E}}+F^\mathcal{E} \right)u+{\buildrel\mathrm{v}\over{\mathrm{div}}}{}^{\pi ^\ast\mathcal{E}}{\buildrel\mathrm{h}\over\nabla}{}^{\pi ^\ast\mathcal{E}}\mathbb{X}u+{\buildrel\mathrm{v}\over{\mathrm{div}}}{}^{\pi ^\ast\mathcal{E}}\mathbb{XX}{\buildrel\mathrm{v}\over\nabla}{}^{\pi ^\ast\mathcal{E}}u,u\rangle _{L^2},$$
$$=\langle -n\mathbb{XX}u+{\buildrel\mathrm{v}\over{\mathrm{div}}}{}^{\pi ^\ast\mathcal{E}} \left(R{\buildrel\mathrm{v}\over\nabla}{}^{\pi ^\ast\mathcal{E}}+F^\mathcal{E} \right)u+{\buildrel\mathrm{v}\over{\mathrm{div}}}{}^{\pi ^\ast\mathcal{E}}\mathbb{XX}{\buildrel\mathrm{v}\over\nabla}{}^{\pi ^\ast\mathcal{E}}u,u\rangle _{L^2},$$
where in the above four lines we used respectively (\ref{MathItem_5.30}), (\ref{MathItem_5.27}), (\ref{MathItem_5.28}), and (\ref{MathItem_5.29}).
Applying Lemma \ref{MathItem_5.10} again gives that this is equal to
$$n\langle \mathbb{X}u,\mathbb{X}u\rangle _{L^2}+\langle {\buildrel\mathrm{v}\over{\mathrm{div}}}{}^{\pi ^\ast\mathcal{E}}R{\buildrel\mathrm{v}\over\nabla}{}^{\pi ^\ast\mathcal{E}}u+{\buildrel\mathrm{v}\over{\mathrm{div}}}{}^{\pi ^\ast\mathcal{E}}F^\mathcal{E}u,u\rangle _{L^2}+\langle \mathbb{X}{\buildrel\mathrm{v}\over\nabla}{}^{\pi ^\ast\mathcal{E}}u,\mathbb{X}{\buildrel\mathrm{v}\over\nabla}{}^{\pi ^\ast\mathcal{E}}u\rangle _{L^2}.$$
Splitting the second inner product over the ``+" sign and then applying Lemma \ref{MathItem_5.10} to the resultant middle two terms proves the theorem.

\begin{flushright}$\blacksquare$\end{flushright}
\subsection{Finite Degree of Solutions to Transport Equations} \label{Section_5.3} 

Throughout this section only we make the additional assumption that $\nabla^\mathcal{E}$ is a smooth \textit{unitary} connection in $\mathcal{E}$.

In the proof of Theorem \ref{MathItem_2.8} we will end up showing that $Q-\mathrm{id}$ satisfies an equation of a form similar to
$$\mathbb{X}u+\Phi u=f$$
over $SM$, which is called a ``transport equation." It turns out that this equation has good behavior with respect to vertical Fourier analysis, which we now introduce.
Consider the \textbf{vertical Laplacian}:
$$\Delta^{\pi ^\ast\mathcal{E}}=-{\buildrel\mathrm{v}\over{\mathrm{div}}}{}^{\pi ^\ast\mathcal{E}}{\buildrel\mathrm{v}\over\nabla}{}^{\pi ^\ast\mathcal{E}} : C^\infty \left(SM;\pi ^\ast\mathcal{E} \right)\longrightarrow C^\infty \left(SM;\pi ^\ast\mathcal{E} \right).$$
By Lemma \ref{MathItem_5.10} 3) this operator has the $\mathcal{R}^l$-mapping property $\mathcal{R}^l \left(SM;\pi ^\ast\mathcal{E} \right)\rightarrow\mathcal{R}^{l-2} \left(SM;\pi ^\ast\mathcal{E} \right)$ for $l\geq2$.
Let us see what this operator looks like in coordinates.
Let $ \left(x^i \right)$ be coordinates of $M$, let $ \left(r_i \right)$ be an \textit{orthonormal} frame of $TM$ over their domain, let $ \left(b_j \right)$ denote an orthonormal frame of $\mathcal{E}$ over their domain, and consider the coordinates
\begin{equation} \label{MathItem_5.31} v^ir_i\longmapsto \left(x^i,v^i \right)\end{equation}
of $TM$.
We claim that for any smooth section $u=u^j\pi ^\ast b_j$,
\begin{equation} \label{MathItem_5.32} \Delta^{\pi ^\ast\mathcal{E}}u= \left(-\Delta^{\mathbb{S}^n}u^j \right)\pi ^\ast b_j\end{equation}
where ``$-\Delta^{\mathbb{S}^n}$" is the positive Laplacian on the $n$-sphere in the variables $v^i$.
This is most easily seen as follows.
Pick an arbitrary point $x_0\in M$ in the domain of our coordinates, choose normal coordinates $ \left({\hat{x}}^i \right)$ of $M$ centered at $x_0$, and consider the coordinates ${\hat{v}}^i\sfrac{\partial }{\partial {\hat{x}}^i}\mapsto \left({\hat{x}}^i,{\hat{v}}^i \right)$ of $TM$.
Then observe that (\ref{MathItem_5.18}) and (\ref{MathItem_5.19}) tell us that on the sphere $S_{x_0}M$, the operator $\Delta^{\pi ^\ast\mathcal{E}}$ applied to $u={\hat{u}}^j\pi ^\ast b_j$ is given by $ \left(-\Delta^{\mathbb{S}^n}{\hat{u}}^j \right)\pi ^\ast b_j$.
The claim then follows by pushing this expression through the change of variables $ \left({\hat{x}}^i,{\hat{v}}^i \right)\mapsto \left(x^i,v^i \right)$.

From this observation and the theory of spherical harmonics (c.f. Section 2.H in \cite{Bibitem_13} for the latter), we obtain several important implications regarding the vertical Laplacian.
First, we get that the eigenvalues of $\Delta^{\pi ^\ast\mathcal{E}}$ match those of $-\Delta^{\mathbb{S}^n}$, which are explicitly given by
$$\lambda_m=m \left(m+n-1 \right)\ \ \ \mathrm{for\ integers}\ \ \ m\geq0.$$
Furthermore, letting $\Omega_m$ denote the set of smooth eigenfunctions of $\Delta^{\pi ^\ast\mathcal{E}}$ with eigenvalue $\lambda_m$, any $u\in C^\infty \left(SM;\pi ^\ast\mathcal{E} \right)$ can be uniquely decomposed as the (pointwise converging) ``Fourier series"
$$u=\sum_{m=0}^{\infty}u_m,\ \ \ \ \ u_m\in\Omega_m.$$
Furthermore, for any fixed $x\in M$ this convergence also holds in $L^2 \left(S_xM; \left(\pi ^\ast\mathcal{E} \right)_x \right)$.
The $u_m$'s are called $u$'s \textbf{Fourier modes}.
The maximum index $m$ for which $u_m\neq0$ is called the \textbf{degree} of $u$ and is denoted by ``$\deg{u}$" (which could be infinity).
Naturally, we say that $u$ is of \textbf{finite degree} if its degree is finite.
We can write an explicit equation for the Fourier modes as follows.
For each $m\in\mathbb{Z}_+$ we let
$$\left\{Y_k^m : k=1,\ldots,l_m\right\}$$
denote a real-valued orthonormal basis of eigenfunctions of $-\Delta^{\mathbb{S}^n}$ with eigenvalue $\lambda_m$.
Then in the coordinates (\ref{MathItem_5.31}) and frame $ \left(b_i \right)$ there
\begin{equation} \label{MathItem_5.33} u_m^j \left(x^i,v^i \right)=\sum_{k=1}^{l_m}{\left[\int_{\mathbb{S}^n}{u^j \left(x^i,w^i \right)Y_k^m \left(w^i \right)dw_{\mathbb{S}^n}}\right]Y_k^m \left(v^i \right)}\end{equation}
(no implicit summation meant in $m$ here).
An important property of the vertical Laplacian eigenspaces is that they are orthogonal (see right after for justification):

\begin{proposition} \label{MathItem_5.34}  For any $x\in M$, the spaces $\Omega_m\cap L^2 \left(S_xM;\pi ^\ast\mathcal{E} \right)$ and $\Omega_{m^\prime}\cap L^2 \left(S_xM;\pi ^\ast\mathcal{E} \right)$ are orthogonal with respect to $L^2 \left(S_xM;\pi ^\ast\mathcal{E} \right)$ when $m\neq m^\prime$.
The same holds for $\Omega_m\cap L^2 \left(SM;\pi ^\ast\mathcal{E} \right)$ and $\Omega_{m^\prime}\cap L^2 \left(SM;\pi ^\ast\mathcal{E} \right)$.
\end{proposition}

From (\ref{MathItem_4.19}) and (\ref{MathItem_4.20}) (observing that $w^{ij}\sfrac{\partial }{\partial  v^i}=\mathcal{K}^{-1} \left(w^{ij}\sfrac{\pi ^\ast\partial }{\partial  x^i} \right)$ there) it follows that another way to express the vertical Laplacian is
$$\Delta^{\pi ^\ast\mathcal{E}}u= \left({\rm div}_{S_xM}{{\grad}_{S_xM}{u^i}} \right)\pi ^\ast b_i= \left(-\Delta_{S_xM}u^i \right)\pi ^\ast b_i.$$
The first part of the above proposition then follows from the fact that eigenspaces of Laplacian-Beltrami operators, such as $-\Delta_{S_xM}$, are orthogonal (e.g. see Problem 16-15 in \cite{Bibitem_25}).
The second part of the proposition follows from this and the fact that integrals over $SM$ can be partitioned as described in Section \ref{Section_4.4} above.

Another important property is that the vertical Laplacian commutes with taking the $m^{\mathrm{th}}$ Fourier mode:

\begin{proposition} \label{MathItem_5.35}  If $u\in C^\infty \left(SM;\pi ^\ast\mathcal{E} \right)$, then $\Delta^{\pi ^\ast\mathcal{E}} \left(u_m \right)= \left(\Delta^{\pi ^\ast\mathcal{E}}u \right)_m$.
\end{proposition}

\noindent\textbf{Proof:} Using the coordinate expression for the vertical Laplacian (\ref{MathItem_5.32}), this follows by taking (\ref{MathItem_5.33}) and integrating by parts:
$$\left[\Delta^{\pi ^\ast\mathcal{E}} \left(u_m \right)\right]^j=\lambda_mu_m^j=\sum_{k=1}^{l_m}{\left[\int_{\mathbb{S}^n}{u^j\lambda_mY_k^mdw_{\mathbb{S}^n}}\right]Y_k^m}=\sum_{k=1}^{l_m}{\left[\int_{\mathbb{S}^n}{u^j \left(-\Delta^{\mathbb{S}^n} \right)Y_k^mdw_{\mathbb{S}^n}}\right]Y_k^m}$$
$$=\sum_{k=1}^{l_m}{\left[\int_{\mathbb{S}^n}{ \left( \left(-\Delta^{\mathbb{S}^n} \right)u^j \right)Y_k^mdw_{\mathbb{S}^n}}\right]Y_k^m}= \left(\Delta^{\pi ^\ast\mathcal{E}}u \right)_m^j.$$
\begin{flushright}$\blacksquare$\end{flushright}

Next we will need the fact that Fourier modes have the same regularity as their original section:

\begin{proposition} \label{MathItem_5.36}  If $u\in\mathcal{R}^l \left(SM;\pi ^\ast\mathcal{E} \right)$ for $l\geq0$, then each Fourier mode $u_m\in\mathcal{R}^l \left(SM;\pi ^\ast\mathcal{E} \right)$ as well.
\end{proposition}

\noindent\textbf{Proof:} Take any $u\in\mathcal{R}^l \left(SM;\pi ^\ast\mathcal{E} \right)$ with $l\geq0$ and fix $m\geq0$.
Let $ \left(x^i \right)$ be coordinates of $\overline{M}$, let $ \left(r_i \right)$ be an \textit{orthonormal} frame of $TM$ over $\dom{ \left(x^i \right)}$, and consider the coordinates $v^ir_i\mapsto \left(x^i,v^i \right)$ of $TM$.
Let $ \left(b_j \right)$ denote an orthonormal frame of $\mathcal{E}$ over $\dom{ \left(x^i \right)}$.
We write $u=u^j\pi ^\ast b_j$.
Furthermore, let $K\subseteq\dom{ \left(x^i \right)}\subseteq\overline{M}$ be any compact subset.

First suppose that $l=0$.
Since $u^j\in\rho^{\sfrac{ \left(n+1 \right)}{2}}L^\infty\left[\pi ^{-1}\left[K\right]\right]$, (\ref{MathItem_5.33}) tells us that each $u_m^j\in\rho^{\sfrac{ \left(n+1 \right)}{2}}L^\infty\left[\pi ^{-1}\left[K\right]\right]$ as well and so indeed $u_m\in\mathcal{R}^l$.

Next suppose that $l=1$.
Suppose that the frame $ \left(r_i \right)$ was obtained by mapping a local orthonormal frame $({\overline{\zeta}}^i)$ of ${}^0T^\ast\overline{M}$ via the canonical identification $\sharp \circ H : \left.{}^0T^\ast\overline{M}\right|_M\rightarrow TM$, and consider the coordinates ${\overline{\eta}}_i{\overline{\zeta}}^i\mapsto \left(x^i,{\overline{\eta}}_i \right)$ of ${}^0T^\ast\overline{M}$.
Observe that this identification is given by ${\overline{\eta}}_i=v^i$.
Now, take any smooth vector field $V_1\in C^\infty \left({}^0S^\ast\overline{M};T{}^0S^\ast\overline{M} \right)$ which we write as $V= \left(V_1 \right)^r\sfrac{\partial }{\partial  x^r}+ \left(V_1 \right)_r\sfrac{\partial }{\partial {\overline{\eta}}_r}$.
Pushing $V$ to $SM$ gives\footnote{ We cannot use the Einstein summation convention on $ \left(V_1 \right)_r\sfrac{\partial }{\partial  v^r}$ because it has two lower indices.} $V= \left(V_1 \right)^r\sfrac{\partial }{\partial  x^r}+\sum_{r=0}^{n}{ \left(V_1 \right)_r\sfrac{\partial }{\partial  v^r}}$.
Thus by (\ref{MathItem_5.33}) 
$$V_1 \left(u_m^j \right)= \left(V_1 \right)^r\sum_{k=1}^{l_m}{\left[\int_{\mathbb{S}^n}{\frac{\partial  u^j}{\partial  x^r}Y_k^mdw_{\mathbb{S}^n}}\right]Y_k^m}+\sum_{k=1}^{l_m}{\left[\int_{\mathbb{S}^n}{u^jY_k^mdw_{\mathbb{S}^n}}\right]\sum_{r=0}^{n}{ \left(V_1 \right)_r\frac{\partial  Y_k^m}{\partial  v^r}}}.$$
Now, the components $ \left(V_1 \right)^i$, $ \left(V_1 \right)_i$ are bounded over ${}^0S^\ast\overline{M}\cap\pi _0^{-1}\left[K\right]$ because they are smooth over this compact set.
For future use, we point out that this also holds for their $\sfrac{\partial }{\partial  x^i}$, $\sfrac{\partial }{\partial  v^i}$ first and higher order partials as well.
Furthermore, since $\sfrac{\partial }{\partial  x^r}$ is a smooth vector field over the 0-cosphere bundle, we have that each $\sfrac{\partial  u^j}{\partial  x^r}$ in the above expression is in $\rho^{\sfrac{ \left(n+1 \right)}{2}}L^\infty\left[\pi ^{-1}\left[K\right]\right]$ by assumption.
Hence the above expression tells us that $V_1 \left(u_m^j \right)\in\rho^{\sfrac{ \left(n+1 \right)}{2}}L^\infty\left[\pi ^{-1}\left[K\right]\right]$ and so indeed $u_m\in\mathcal{R}^l$.
The cases $l\geq2$ are handled similarly.

\begin{flushright}$\blacksquare$\end{flushright}

One of the central properties of $\mathbb{X}=\nabla_X^{\pi ^\ast\mathcal{E}}$ is that it maps
\begin{equation} \label{MathItem_5.37} \mathbb{X} : \Omega_m\longrightarrow\Omega_{m-1}\oplus\Omega_{m+1}.\end{equation}
This is proven in Section 3.4 of \cite{Bibitem_18}.
Similarly, multiplication on the left by $\Phi$ maps $\Omega_m\rightarrow\Omega_m$ since $\Phi$ has no dependence on the vertical variable ``$v$." In particular, we see that the operator in the transport equation ``$\mathbb{X}+\Phi$" maps sections of finite degree to sections of finite degree.
The converse is also true, which is the main result of this section (note the assumptions of ``skew-Hermitian" and ``unitary"):

\begin{theorem} \label{MathItem_5.38}  Suppose that $ \left(M\subseteq\overline{M},g \right)$ is an asymptotically hyperbolic space, $ \left(\mathcal{E},\langle \cdot,\cdot\rangle _\mathcal{E} \right)$ is a smooth complex Hermitian vector bundle over $\overline{M}$, $\Phi\in C^\infty \left(\overline{M};{\End}_{\mathrm{sk}}{\mathcal{E}} \right)$, and that $\nabla^\mathcal{E}$ is a smooth unitary connection in $\mathcal{E}$.
Assume also that the sectional curvatures of $g$ are negative.
If $u\in\mathcal{R}^3 \left(SM;\pi ^\ast\mathcal{E} \right)$ solves
\begin{equation} \label{MathItem_5.39} \mathbb{X}u+\Phi u=f\end{equation}
for some $f\in C^\infty \left(SM;\pi ^\ast\mathcal{E} \right)$ of finite degree, then $u$ is also of finite degree.
\end{theorem}

We note that the above theorem appears as Theorem 4.6 in \cite{Bibitem_18} in the setting of compact manifolds.
To prove Theorem \ref{MathItem_5.38}, we need several preliminary results.

\begin{lemma} \label{MathItem_5.40}  It holds that
$$\left[\mathbb{X},\Delta^{\pi ^\ast\mathcal{E}}\right]=2{\buildrel\mathrm{v}\over{\mathrm{div}}}{}^{\pi ^\ast\mathcal{E}}{\buildrel\mathrm{h}\over\nabla}{}^{\pi ^\ast\mathcal{E}}+n\mathbb{X}.$$
\end{lemma}

The above lemma is stated as Lemma 3.4 of \cite{Bibitem_18}, whose proof is essentially identical to that of Lemma 3.5 in \cite{Bibitem_41}.

To state the next preliminary result we observe that because of (\ref{MathItem_5.37}), over each $\Omega_m$ we can decompose $\mathbb{X}=\mathbb{X}_-+\mathbb{X}_+$ where
\begin{equation} \label{MathItem_5.41} \mathbb{X}_\pm : \Omega_m\longrightarrow\Omega_{m\pm1}.\end{equation}
We point out that the maps $\mathbb{X}_\pm$ are distinct for different $\Omega_m$ even though we use the same notation to denote them. 

\begin{remark} \label{MathItem_5.42}  \normalfont Like $\mathbb{X}$, the operators $\mathbb{X}_\pm$ map $\mathcal{R}^l\rightarrow\mathcal{R}^{l-1}$ for $l\geq1$ by Lemma \ref{MathItem_5.10} and Proposition \ref{MathItem_5.36} above.
\end{remark}

We mention that the idea of splitting the action of the geodesic vector field as above was first introduced by Guillemin and Kazhdan - see \cite{Bibitem_19}.
The following preliminary result is a special case of the Pestov identity with a connection (Theorem \ref{MathItem_5.9}):

\begin{proposition} \label{MathItem_5.43}  Suppose that $u\in\Omega_m\cap\mathcal{R}^2 \left(SM;\pi ^\ast\mathcal{E} \right)$.
Then
$$ \left(2m+n \right)\left\lVert \mathbb{X}_+u\right\rVert _{L^2}^2$$
$$=\left\lVert {\buildrel\mathrm{h}\over\nabla}{}^{\pi ^\ast\mathcal{E}}u\right\rVert _{L^2}^2+ \left(2m+n-2 \right)\left\lVert \mathbb{X}_-u\right\rVert _{L^2}^2-\langle R{\buildrel\mathrm{v}\over\nabla}{}^{\pi ^\ast\mathcal{E}}u,{\buildrel\mathrm{v}\over\nabla}{}^{\pi ^\ast\mathcal{E}}u\rangle _{L^2}-\langle F^\mathcal{E}u,{\buildrel\mathrm{v}\over\nabla}{}^{\pi ^\ast\mathcal{E}}u\rangle _{L^2}.$$
\end{proposition}

\noindent\textbf{Proof:} We have that $u$ satisfies the equation in Theorem \ref{MathItem_5.9}.
Let us take a look at the term
\begin{align*}\left\lVert \mathbb{X}{\buildrel\mathrm{v}\over\nabla}{}^{\pi ^\ast\mathcal{E}}u\right\rVert _{L^2}^2=\langle -{\buildrel\mathrm{h}\over\nabla}{}^{\pi ^\ast\mathcal{E}}u+{\buildrel\mathrm{v}\over\nabla}{}^{\pi ^\ast\mathcal{E}}\mathbb{X}u,-{\buildrel\mathrm{h}\over\nabla}{}^{\pi ^\ast\mathcal{E}}u+{\buildrel\mathrm{v}\over\nabla}{}^{\pi ^\ast\mathcal{E}}\mathbb{X}u\rangle _{L^2}& &\mathrm{by\ (\ref{MathItem_5.27})},\\=\left\lVert {\buildrel\mathrm{h}\over\nabla}{}^{\pi ^\ast\mathcal{E}}u\right\rVert _{L^2}^2+2\langle \mathbb{X}u,{\buildrel\mathrm{v}\over{\mathrm{div}}}{}^{\pi ^\ast\mathcal{E}}{\buildrel\mathrm{h}\over\nabla}{}^{\pi ^\ast\mathcal{E}}u\rangle _{L^2}+\left\lVert {\buildrel\mathrm{v}\over\nabla}{}^{\pi ^\ast\mathcal{E}}\mathbb{X}u\right\rVert _{L^2}^2& &\mathrm{Lemma\ \ref{MathItem_5.10}\ 3)\ and\ Lemma\ \ref{MathItem_5.20}}.\end{align*}
Applying Lemma \ref{MathItem_5.40}, we see that the middle term in the last quantity is equal to
$$\langle \mathbb{X}u,\mathbb{X}\Delta^{\pi ^\ast\mathcal{E}}u-\Delta^{\pi ^\ast\mathcal{E}}\mathbb{X}u-n\mathbb{X}u\rangle _{L^2}.$$
Splitting $\mathbb{X}u=\mathbb{X}_-u+\mathbb{X}_+u\in\Omega_{m-1}\oplus\Omega_{m+1}$, using that
$$\Delta^{\pi ^\ast\mathcal{E}}u=\lambda_mu,\ \ \ \ \ \ \ \ \Delta^{\pi ^\ast\mathcal{E}}\mathbb{X}_-u=\lambda_{m-1}\mathbb{X}_-u,\ \ \ \ \ \ \ \ \Delta^{\pi ^\ast\mathcal{E}}\mathbb{X}_+u=\lambda_{m+1}\mathbb{X}_+u,$$
using the orthogonality of the vertical Laplacian eigenspaces, and then plugging the result into the equation in Theorem \ref{MathItem_5.9} proves the proposition after several cancellations.

\begin{flushright}$\blacksquare$\end{flushright}

The following lemma provides the contraction property that is needed in the proof of Theorem \ref{MathItem_5.38}.

\begin{lemma} \label{MathItem_5.44}  Suppose that the sectional curvatures of $g$ are negative.
Then there exist real constants $c_m\rightarrow\infty$ such that for sufficiently large $m$,
\begin{equation} \label{MathItem_5.45} \begin{cases}\left\lVert \mathbb{X}_-u\right\rVert _{L^2}^2+c_m\left\lVert u\right\rVert _{L^2}^2\le\left\lVert \mathbb{X}_+u\right\rVert _{L^2}^2&\mathrm{if}\ \ \ n\neq2,\\\left\lVert \mathbb{X}_-u\right\rVert _{L^2}^2+c_m\left\lVert u\right\rVert _{L^2}^2\le d_m\left\lVert \mathbb{X}_+u\right\rVert _{L^2}^2&\mathrm{if}\ \ \ n=2,\\\end{cases}\end{equation}
for all $u\in\Omega_m\cap\mathcal{R}^2 \left(SM;\pi ^\ast\mathcal{E} \right)$ where $d_m=1+\sfrac{1}{\left[ \left(2m-1 \right) \left(m+1 \right)^2\right]}$.
\end{lemma}

\noindent\textbf{Proof:} We begin by using the fact that the sectional curvatures of $g$ tend to $-1$ at $\partial \overline{M}$.
Precisely, by the remark after Proposition 1.10 in \cite{Bibitem_31} there exists an $\varepsilon>0$ so that the sectional curvatures of $g$ are less than $-\kappa^\prime$ for some $\kappa^\prime>0$ over the region $\left\{\rho<\varepsilon\right\}$.
Hence this, the compactness of $\left\{\rho\geq\varepsilon\right\}$, and the negative curvature assumption imply that there exists a $\kappa>0$ such that the sectional curvatures of $g$ are bounded above by $-\kappa$ on all of $M$.

Now, take any $u$ as in the statement of the lemma.
We have that $u$ satisfies the equation in Proposition \ref{MathItem_5.43} above.
We begin by estimating the $L^2$ norm of the term ${\buildrel\mathrm{h}\over\nabla}{}^{\pi ^\ast\mathcal{E}}u$ by utilizing the trick of looking at its vertical divergence.
By Lemma \ref{MathItem_5.40} we have that
$${\buildrel\mathrm{v}\over{\mathrm{div}}}{}^{\pi ^\ast\mathcal{E}}{\buildrel\mathrm{h}\over\nabla}{}^{\pi ^\ast\mathcal{E}}u=\frac{1}{2}\mathbb{X}\Delta^{\pi ^\ast\mathcal{E}}u-\frac{1}{2}\Delta^{\pi ^\ast\mathcal{E}}\mathbb{X}u-\frac{n}{2}\mathbb{X}u,$$
$$=\frac{1}{2} \left(\mathbb{X}_+\lambda_mu+\mathbb{X}_-\lambda_mu \right)-\frac{1}{2} \left(\lambda_{m+1}\mathbb{X}_+u+\lambda_{m-1}\mathbb{X}_-u \right)-\frac{n}{2} \left(\mathbb{X}_+u+\mathbb{X}_-u \right),$$
$$=- \left(m+n \right)\mathbb{X}_+u+ \left(m-1 \right)\mathbb{X}_-u,$$
$$=-{\buildrel\mathrm{v}\over{\mathrm{div}}}{}^{\pi ^\ast\mathcal{E}} \left(-\frac{m+n}{\lambda_{m+1}}{\buildrel\mathrm{v}\over\nabla}{}^{\pi ^\ast\mathcal{E}}\mathbb{X}_+u+\frac{m-1}{\lambda_{m-1}}{\buildrel\mathrm{v}\over\nabla}{}^{\pi ^\ast\mathcal{E}}\mathbb{X}_-u \right).$$
Plugging the expression for the $\lambda_k$'s into this, we conclude that
$${\buildrel\mathrm{h}\over\nabla}{}^{\pi ^\ast\mathcal{E}}u=\frac{1}{m+1}{\buildrel\mathrm{v}\over\nabla}{}^{\pi ^\ast\mathcal{E}}\mathbb{X}_+u-\frac{1}{m+n-2}{\buildrel\mathrm{v}\over\nabla}{}^{\pi ^\ast\mathcal{E}}\mathbb{X}_-u+Z$$
where $Z\in C^\infty \left(SM;N\otimes\pi ^\ast\mathcal{E} \right)$ is such that ${\buildrel\mathrm{v}\over{\mathrm{div}}}{}^{\pi ^\ast\mathcal{E}}Z=0$.
Since ${\buildrel\mathrm{v}\over{\mathrm{div}}}{}^{\pi ^\ast\mathcal{E}}$ is the formal $L^2$ adjoint of ${\buildrel\mathrm{v}\over\nabla}{}^{\pi ^\ast\mathcal{E}}$ over each sphere $S_xM$, it follows that $Z$ is perpendicular to the other two terms on the right-hand side with respect to $L^2 \left(S_xM;\pi ^\ast\mathcal{E} \right)$.
From the orthogonality of the vertical Laplacian eigenspaces with respect to $L^2 \left(S_xM;\pi ^\ast\mathcal{E} \right)$ we get the following estimate over each sphere:
$$\left\lVert {\buildrel\mathrm{h}\over\nabla}{}^{\pi ^\ast\mathcal{E}}u\right\rVert _{L^2 \left(S_xM;\pi ^\ast\mathcal{E} \right)}^2\geq\frac{m+n}{m+1}\left\lVert \mathbb{X}_+u\right\rVert _{L^2 \left(S_xM;\pi ^\ast\mathcal{E} \right)}^2+\frac{m-1}{m+n-2}\left\lVert \mathbb{X}_-u\right\rVert _{L^2 \left(S_xM;\pi ^\ast\mathcal{E} \right)}^2.$$
Integrating this in $x\in M$ and using (\ref{MathItem_4.7}) implies that this holds for $L^2 \left(SM;\pi ^\ast\mathcal{E} \right)$.
Plugging this observation into the equation in Proposition \ref{MathItem_5.43} gives (here $L^2=L^2 \left(SM;\pi ^\ast\mathcal{E} \right)$)
\begin{equation} \label{MathItem_5.46}  \left(2m+n-\frac{m+n}{m+1} \right)\left\lVert \mathbb{X}_+u\right\rVert _{L^2}^2\end{equation}
$$\geq \left(2m+n-2+\frac{m-1}{m+n-2} \right)\left\lVert \mathbb{X}_-u\right\rVert _{L^2}^2-\langle R{\buildrel\mathrm{v}\over\nabla}{}^{\pi ^\ast\mathcal{E}}u,{\buildrel\mathrm{v}\over\nabla}{}^{\pi ^\ast\mathcal{E}}u\rangle _{L^2}-\langle F^\mathcal{E}u,{\buildrel\mathrm{v}\over\nabla}{}^{\pi ^\ast\mathcal{E}}u\rangle _{L^2}.$$
We have that the term
\begin{equation} \label{MathItem_5.47} -\langle R{\buildrel\mathrm{v}\over\nabla}{}^{\pi ^\ast\mathcal{E}}u,{\buildrel\mathrm{v}\over\nabla}{}^{\pi ^\ast\mathcal{E}}u\rangle _{L^2}\geq\kappa\left\lVert {\buildrel\mathrm{v}\over\nabla}{}^{\pi ^\ast\mathcal{E}}u\right\rVert _{L^2}^2=\kappa\lambda_m\left\lVert u\right\rVert _{L^2}^2.\end{equation}
Furthermore,
\begin{equation} \label{MathItem_5.48} -\langle F^\mathcal{E}u,{\buildrel\mathrm{v}\over\nabla}{}^{\pi ^\ast\mathcal{E}}u\rangle _{L^2}\geq-\left\lVert F^\mathcal{E}\right\rVert _{L^\infty}\left\lVert u\right\rVert _{L^2}\left\lVert {\buildrel\mathrm{v}\over\nabla}{}^{\pi ^\ast\mathcal{E}}u\right\rVert _{L^2}=-\left\lVert F^\mathcal{E}\right\rVert _{L^\infty}\lambda_m^{\sfrac{1}{2}}\left\lVert u\right\rVert _{L^2}^2.\end{equation}
Hence, letting $a_{m,n}$ and $b_{m,n}$ denote the coefficients of $\left\lVert \mathbb{X}_+u\right\rVert _{L^2}^2$ and $\left\lVert \mathbb{X}_-u\right\rVert _{L^2}^2$ in (\ref{MathItem_5.46}) above respectively we get that
$$\frac{a_{m,n}}{b_{m,n}}\left\lVert \mathbb{X}_+u\right\rVert _{L^2}^2\geq\left\lVert \mathbb{X}_-u\right\rVert _{L^2}^2+\frac{\kappa\lambda_m-\left\lVert F^\mathcal{E}\right\rVert _{L^\infty}\lambda_m^{\sfrac{1}{2}}}{b_{m,n}}\left\lVert u\right\rVert _{L^2}^2.$$
Elementary algebra shows that $\sfrac{a_{m,n}}{b_{m,n}}$ is less than or equal to 1 if $n\neq2$ and $m>1$, and is equal to $d_m$ if $n=2$.
Since $\lambda_m=O \left(m^2 \right)$, the lemma follows.

\begin{flushright}$\blacksquare$\end{flushright}

We need one last technical lemma:

\begin{lemma} \label{MathItem_5.49}  If $u\in\mathcal{R}^3 \left(SM;\pi ^\ast\mathcal{E} \right)$, then $\left\lVert \mathbb{X}_+u_m\right\rVert _{L^2}\rightarrow0$ as $m\rightarrow\infty$.
\end{lemma}

\noindent\textbf{Proof:} Take any $u\in\mathcal{R}^3 \left(SM;\pi ^\ast\mathcal{E} \right)$.
Let $ \left(x^i \right)$ be coordinates of $\overline{M}$, let $ \left(b_j \right)$ denote a frame of $\mathcal{E}$ over $\dom{ \left(x^i \right)}$, and consider the coordinates $v^i\sfrac{\partial }{\partial  x^i}\mapsto \left(x^i,v^i \right)$ of $TM$.
Take any compact subset $K\subseteq\dom{ \left(x^i \right)}\subseteq\overline{M}$ and suppose for convenience that $ \left(b_j \right)$ is orthonormal.
We will show that $\left\lVert \mathbb{X}_+u_m\right\rVert _{L^2 \left(\pi ^{-1}\left[K\right] \right)}\rightarrow0$ as $m\rightarrow\infty$, from which the lemma will follow by covering the compact $\overline{M}$ by a finite number of such sets $K$.

We have that (we omit writing $\pi ^{-1}\left[K\right]$ for the rest of the proof)
\begin{equation} \label{MathItem_5.50} \left\lVert \mathbb{X}_+u_m\right\rVert _{L^2}\le\left\lVert \mathbb{X}u_m\right\rVert _{L^2}=\left\lVert \left[X \left(u_m^k \right)+{}^\mathcal{E}\Gamma_{ij}^kv^iu_m^j\right]\pi ^\ast b_k\right\rVert _{L^2}\end{equation}
$$\le\sum_{k=1}^{d} \left(\left\lVert X \left(u_m^k \right)\right\rVert _{L^2}+\sup_{v\in\pi ^{-1}\left[K\right]}{\left|{}^\mathcal{E}\Gamma_{ij}^kv^i\right|}\left\lVert u_m^j\right\rVert _{L^2} \right).$$
Since $g_{ij}v^iv^j\equiv1$ and $g_{ij}$ blows up like $\rho^{-2}$ at the boundary, we have that each $v^i$ is bounded over the compact $K$ and hence each $\sup{\ldots}$ above is finite.
From Proposition \ref{MathItem_5.35} it follows that
$$\left\lVert u_m^j\right\rVert _{L^2}\le\left\lVert u_m\right\rVert _{L^2}=\frac{1}{\lambda_m}\left\lVert  \left(\Delta^{\pi ^\ast\mathcal{E}}u \right)_m\right\rVert _{L^2}\le\frac{1}{\lambda_m}\left\lVert \Delta^{\pi ^\ast\mathcal{E}}u\right\rVert _{L^2}\rightarrow0\ \ \ \ \ \mathrm{as}  \ m\rightarrow\infty,$$
where we have used that $\Delta^{\pi ^\ast\mathcal{E}}u\in\mathcal{R}^{3-2} \left(SM;\pi ^\ast\mathcal{E} \right)\subseteq L^2 \left(SM;\pi ^\ast\mathcal{E} \right)$ and so the last $L^2$-norm is finite.
Thus by (\ref{MathItem_5.50}), the lemma will be proved if we show that $\left\lVert X \left(u_m^k \right)\right\rVert _{L^2}\rightarrow0$ as $m\rightarrow\infty$.

Let $ \left(r_i \right)$ be an \textit{orthonormal} frame of $TM$ over $\dom{ \left(x^i \right)}$ and now consider the coordinates $v^ir_i\mapsto \left(x^i,v^i \right)$ of $TM$.
Suppose furthermore that the frame $ \left(r_i \right)$ was obtained by mapping a local orthonormal frame $({\overline{\zeta}}^i)$ of ${}^0T^\ast\overline{M}$ via the canonical identification $\sharp \circ H : \left.{}^0T^\ast\overline{M}\right|_M\rightarrow TM$, and consider the coordinates ${\overline{\eta}}_i{\overline{\zeta}}^i\mapsto \left(x^i,{\overline{\eta}}_i \right)$ of ${}^0T^\ast\overline{M}$.
Recall that this identification is given by ${\overline{\eta}}_i=v^i$.
In (\ref{MathItem_5.13}) we pointed out that $X$ extends to a smooth vector field $X_0$ over ${}^0S^\ast\overline{M}$.
Let us write $X_0= \left(X_0 \right)^i\sfrac{\partial }{\partial  x^i}+ \left(X_0 \right)_i\sfrac{\partial }{\partial {\overline{\eta}}_i}$ over ${}^0S^\ast\overline{M}$, which if we push to $SM$ we get $X=\sfrac{ \left(X_0 \right)^i\partial }{\partial  x^i}+\sum_{i=0}^{n}{ \left(X_0 \right)_i\sfrac{\partial }{\partial  v^i}}$.
Thus
$$Xu_m^k= \left(X_0 \right)^i\frac{\partial  u_m^k}{\partial  x^i}+\sum_{i=0}^{n}{ \left(X_0 \right)_i\frac{\partial  u_m^k}{\partial  v^i}}.$$
where $\sfrac{\partial  u_m^k}{\partial  v^i}$ denotes $\sfrac{\partial }{\partial  v^i}$ applied to any smooth extension of $u_m^k$ from $SM$.
For definitiveness, let us say that we extend $u_m^k$ to $TM$ minus the zero section by making $u_m^k$ constant along the radial lines $t\mapsto tv$ for $t\in \left(0,\infty \right)$ and $v\in SM$.
Now, the components $ \left(X_0 \right)^i$, $ \left(X_0 \right)_i$ are bounded above $K$ because they are smooth over the compact ${}^0S^\ast\overline{M}\cap\pi _0^{-1}\left[K\right]$.
Hence by the above equation, the lemma will be proved if we show that both $\left\lVert \sfrac{\partial  u_m^k}{\partial  x^i}\right\rVert _{L^2},\left\lVert \sfrac{\partial  u_m^k}{\partial  v^i}\right\rVert _{L^2}\rightarrow0$ as $m\rightarrow\infty$.

First let us take a look at $\left\lVert \sfrac{\partial  u_m^k}{\partial  x^i}\right\rVert _{L^2}$.
By tucking $\sfrac{\partial }{\partial  x^i}$ under the integral sign in (\ref{MathItem_5.33}), it follows that
$$\frac{\partial  u_m^k}{\partial  x^i}= \left(\frac{\partial }{\partial  x^i}u \right)_m^k.$$
where $\sfrac{\partial  u}{\partial  x^i}$ denotes $ \left(\sfrac{\partial  u^k}{\partial  x^i} \right)\pi ^\ast b_k$.
Thus
$$\left\lVert \frac{\partial  u_m^k}{\partial  x^i}\right\rVert _{L^2}^2=\left\lVert  \left(\frac{\partial }{\partial  x^i}u \right)_m^k\right\rVert _{L^2}^2\le\left\lVert \frac{\partial }{\partial  x^i}u\right\rVert _{L^2}^2=\frac{1}{\lambda_m^2}\left\lVert \frac{\partial }{\partial  x^i}\Delta^{\pi ^\ast\mathcal{E}}u\right\rVert _{L^2}^2\rightarrow0\ \ \ \ \ \ \ \ \mathrm{as}  \ m\rightarrow\infty,$$
where we have used the fact that $\sfrac{\partial }{\partial  x^i}$ is a smooth vector field over the 0-cosphere bundle and hence $\sfrac{\partial }{\partial  x^i} \left(\Delta^{\pi ^\ast\mathcal{E}}u \right)\in\mathcal{R}^{3-3}\subseteq L^2$.

Finally, let us take a look at $\left\lVert \sfrac{\partial  u_m^k}{\partial  v^i}\right\rVert _{L^2}$.
Since we extended $u_m^k$ to be constant radially, for any fixed $x\in M$ we have that
$$\int_{S_xM}{\left|\frac{\partial  u_m^k}{\partial  v^i}\right|^2dS_xM}=\int_{\mathbb{S}^n}{\left|\frac{\partial  u_m^k}{\partial  v^i} \left(x^i,v^i \right)\right|^2dv_{\mathbb{S}^n}}\le\int_{\mathbb{S}^n}{\left|{\grad}_{v\in\mathbb{S}^n}{u_m^k}\right|^2dv_{\mathbb{S}^n}}$$
where ${\mathrm{grad}}_{v\in\mathbb{S}^n}$ means the gradient in $ \left(v^i \right)$ over the Euclidean sphere $\left\{ \left(v^1 \right)^2+\ldots \left(v^n \right)^2=1\right\}$ with respect to usual spherical metric.
Looking at the last integral, integrating by parts gives
$$\int_{\mathbb{S}^n}{\left|{\grad}_{v\in\mathbb{S}^n}{u_m^k}\right|^2dv_{\mathbb{S}^n}}=\int_{\mathbb{S}^n}{u_m^k \left(-\Delta^{\mathbb{S}^n}u_m^k \right)dv_{\mathbb{S}^n}}=\lambda_m\int_{S_xM}{\left|u_m^k\right|^2dS_xM}.$$
Integrating in $x$ then gives that
$$\left\lVert \frac{\partial  u_m^k}{\partial  v^i}\right\rVert _{L^2}^2\le\lambda_m\left\lVert u_m^k\right\rVert _{L^2}^2\le\frac{1}{\lambda_m}\left\lVert \Delta^{\pi ^\ast\mathcal{E}}u\right\rVert _{L^2}^2\rightarrow0\ \ \ \ \ \mathrm{as}  \ m\rightarrow\infty.$$
As discussed above, this proves the lemma.

\begin{flushright}$\blacksquare$\end{flushright}

\noindent\textit{Proof of Theorem \ref{MathItem_5.38}:}

We start by assuming that $n\neq2$ since the proof of the case $n=2$ requires a slight modification.
Comparing Fourier modes of order $m>\deg{f}$ in (\ref{MathItem_5.39}) gives that
\begin{equation} \label{MathItem_5.51} \mathbb{X}_+u_{m-1}+\mathbb{X}_-u_{m+1}+\Phi u_m=0.\end{equation}
Using this relation and plugging $u_{m+1}$ into $u$ in (\ref{MathItem_5.45}) gives that
\begin{equation} \label{MathItem_5.52} \left\lVert \mathbb{X}_+u_{m+1}\right\rVert _{L^2}^2\geq\left\lVert \mathbb{X}_+u_{m-1}\right\rVert _{L^2}^2+\left\lVert \Phi u_m\right\rVert _{L^2}^2+c_{m+1}\left\lVert u_{m+1}\right\rVert _{L^2}^2+2\mathrm{Re}{\langle \Phi u_m,\mathbb{X}_+u_{m-1}\rangle _{L^2}}.\end{equation}
The idea of the proof is the following.
One can plug the relation (\ref{MathItem_5.51}) with ``$m+1$" replaced by ``$m-1$" into $\left\lVert \mathbb{X}_+u_{m-1}\right\rVert _{L^2}^2$ in (\ref{MathItem_5.52}), use (\ref{MathItem_5.45}) again, and then proceed recursively.
One will get a long expression on the right which one needs to cleverly manipulate to bound $\left\lVert u_{m_0}\right\rVert _{L^2}^2$ for some fixed index $m_0$.
Then using that $\left\lVert \mathbb{X}_+u_{m+1}\right\rVert _{L^2}^2$ goes to zero as $m\rightarrow\infty$ by Lemma \ref{MathItem_5.49} will force $\left\lVert u_{m_0}\right\rVert _{L^2}^2=0$.
We will show that this holds for large enough $m_0$, from which the theorem will follow.
The first obstacle to accomplishing this is the inner product term $2\mathrm{Re}{\langle \Phi u_m,\mathbb{X}_+u_{m-1}\rangle _{L^2}}$ above.
The following claim helps resolve this.

Recall from Section \ref{Section_4.7} the definition of $\mathbb{X}\Phi=\nabla_X^{\pi ^\ast\End{\mathcal{E}}}\Phi$.
In particular, since $\Phi$ is smooth on $\overline{M}$ and hence on ${}^0S^\ast\overline{M}$ (when lifted), it follows from (\ref{MathItem_4.26}), (\ref{MathItem_5.12}) and the sentence after, and (\ref{MathItem_5.13}) that the quantity $\mathbb{X}\Phi$ extends smoothly to ${}^0S^\ast\overline{M}$.

\noindent\textit{Claim:} The following identity is true (because $\Phi$ is skew-Hermitian):
$$\langle \mathbb{X}_+u_{m-1},\Phi u_m\rangle _{L^2}+{\overline{\langle \mathbb{X}_+u_{m-2},\Phi u_{m-1}\rangle }}_{L^2}=-\langle u_{m-1}, \left(\mathbb{X}\Phi \right)u_m\rangle _{L^2}-\left\lVert \Phi u_{m-1}\right\rVert _{L^2}^2.$$
\noindent\textit{Proof of claim:} This is simply a computation:
\begin{align*}\langle \mathbb{X}_+u_{m-1},\Phi u_m\rangle _{L^2}=\langle \mathbb{X}u_{m-1},\Phi u_m\rangle _{L^2}& &\mathrm{use}\ \mathbb{X}=\mathbb{X}_-+\mathbb{X}_+\ \mathrm{and\ Proposition\ \ref{MathItem_5.34}} ,\\=-\langle u_{m-1},\mathbb{X} \left(\Phi u_m \right)\rangle _{L^2}& &\mathrm{Lemma\ \ref{MathItem_5.10}},\\=-\langle u_{m-1}, \left(\mathbb{X}\Phi \right)u_m+\Phi\mathbb{X}u_m\rangle _{L^2}& &\mathrm{definition\ of}\ \mathbb{X}\Phi,\\=-\langle u_{m-1}, \left(\mathbb{X}\Phi \right)u_m\rangle _{L^2}-\langle u_{m-1},\Phi\mathbb{X}_-u_m\rangle _{L^2}& &\mathbb{X}=\mathbb{X}_-+\mathbb{X}_+\ \mathrm{and\ Proposition\ \ref{MathItem_5.34}} ,\\=-\langle u_{m-1}, \left(\mathbb{X}\Phi \right)u_m\rangle _{L^2}+\langle u_{m-1},\Phi\left[\Phi u_{m-1}+\mathbb{X}_+u_{m-2}\right]\rangle _{L^2}& &\mathrm{used\ (\ref{MathItem_5.51})}.\end{align*}
From here the claim follows by rearranging and using that $\Phi$ is skew-Hermitian.

\begin{flushright}\textit{End of proof of claim.}\end{flushright}

We return to proving the theorem.
The above claim tells us that to get rid of the last term in (\ref{MathItem_5.52}), we can add (\ref{MathItem_5.52}) and (\ref{MathItem_5.52}) with ``$m$" replaced by ``$m-1$":
$$\left\lVert \mathbb{X}_+u_{m+1}\right\rVert _{L^2}^2+\left\lVert \mathbb{X}_+u_m\right\rVert _{L^2}^2\geq\left\lVert \mathbb{X}_+u_{m-1}\right\rVert _{L^2}^2+\left\lVert \mathbb{X}_+u_{m-2}\right\rVert _{L^2}^2+\left\lVert \Phi u_m\right\rVert _{L^2}^2+\left\lVert \Phi u_{m-1}\right\rVert _{L^2}^2$$
$$+c_{m+1}\left\lVert u_{m+1}\right\rVert _{L^2}^2+c_m\left\lVert u_m\right\rVert _{L^2}^2+2\mathrm{Re}{\langle \Phi u_m,\mathbb{X}_+u_{m-1}\rangle _{L^2}}+2\mathrm{Re}{\langle \Phi u_{m-1},\mathbb{X}_+u_{m-2}\rangle _{L^2}},$$
and then apply the equation in the claim to get the following inequality, where for brevity $a_m=\left\lVert \mathbb{X}_+u_m\right\rVert _{L^2}^2+\left\lVert \mathbb{X}_+u_{m-1}\right\rVert _{L^2}^2$:
$$a_{m+1}\geq a_{m-1}+\left\lVert \Phi u_m\right\rVert _{L^2}^2+\left\lVert \Phi u_{m-1}\right\rVert _{L^2}^2+c_{m+1}\left\lVert u_{m+1}\right\rVert _{L^2}^2+c_m\left\lVert u_m\right\rVert _{L^2}^2$$
$$-2\mathrm{Re}{\langle u_{m-1}, \left(\mathbb{X}\Phi \right)u_m\rangle _{L^2}}-2\left\lVert \Phi u_{m-1}\right\rVert _{L^2}^2$$
$$\geq a_{m-1}+c_{m+1}\left\lVert u_{m+1}\right\rVert _{L^2}^2+c_m\left\lVert u_m\right\rVert _{L^2}^2-\left\lVert u_{m-1}\right\rVert _{L^2}^2-\left\lVert  \left(\mathbb{X}\Phi \right)u_m\right\rVert _{L^2}^2-\left\lVert \Phi u_{m-1}\right\rVert _{L^2}^2.$$
Since $\Phi$ and $\mathbb{X}\Phi$ are smooth (and hence continuous) on the compact manifold ${}^0S^\ast\overline{M}$, there exist constants $B,C>0$ such that
$$\left\lVert \Phi h\right\rVert _{L^2}^2\le B\left\lVert h\right\rVert _{L^2}^2,$$
$$\left\lVert  \left(\mathbb{X}\Phi \right)h\right\rVert _{L^2}^2\le C\left\lVert h\right\rVert _{L^2}^2,$$
for any $h\in C^\infty \left(SM;\pi ^\ast\mathcal{E} \right)\cap L^2 \left(SM;\pi ^\ast\mathcal{E} \right)$.
Hence
\begin{equation} \label{MathItem_5.53} a_{m+1}\geq a_{m-1}+r_m\end{equation}
where
$$r_m:=c_{m+1}\left\lVert u_{m+1}\right\rVert _{L^2}^2+ \left(c_m-C \right)\left\lVert u_m\right\rVert _{L^2}^2- \left(1+B \right)\left\lVert u_{m-1}\right\rVert _{L^2}^2.$$
Applying (\ref{MathItem_5.53}) recursively gives
$$a_{m+1}=a_{m_0-1}+r_{m_0}+r_{m_0+2}+\ldots+r_m$$
for any pair of indices $m,m_0>\deg{f}$ such that $m=m_0+2k$ for some integer $k\geq0$.
We seek to bound the resultant tail of $r_i$'s.
To do so, choose $m_0$ big enough so that for $m\geq m_0$, $c_m$ is bigger than both $C$ and $B+1$.
Hence if $m=m_0+2k$,
$$r_{m_0}+r_{m_0+2}+\ldots+r_m=c_{m+1}\left\lVert u_{m+1}\right\rVert _{L^2}^2+\sum_{i=0}^{k-1}{ \left(c_{m_0+1+2i}- \left(1+B \right) \right)\left\lVert u_{m_0+1+2i}\right\rVert _{L^2}^2}$$
$$+\sum_{i=0}^{k}{ \left(c_{m_0+2i}-C \right)\left\lVert u_{m_0+2i}\right\rVert _{L^2}^2}- \left(1+B \right)\left\lVert u_{m_0-1}\right\rVert _{L^2}^2$$
$$\geq c_{m+1}\left\lVert u_{m+1}\right\rVert _{L^2}^2- \left(1+B \right)\left\lVert u_{m_0-1}\right\rVert _{L^2}^2.$$
Hence for any such $m=m_0+2k$ we get that
$$a_{m+1}\geq a_{m_0-1}- \left(1+B \right)\left\lVert u_{m_0-1}\right\rVert _{L^2}^2.$$
By the definition of $a_{m_0-1}$ and (\ref{MathItem_5.45}) we have that $a_{m_0-1}\geq c_{m_0-1}\left\lVert u_{m_0-1}\right\rVert _{L^2}^2$, and hence we finally arrive at
$$a_{m+1}\geq \left(c_{m_0-1}- \left(1+B \right) \right)\left\lVert u_{m_0-1}\right\rVert _{L^2}^2.$$
Assume we defined $m_0$ before so that $c_{m_0-1}\geq1+B$ as well.
Observe that $a_{m+1}\rightarrow0$ as $m\rightarrow\infty$ by Lemma \ref{MathItem_5.49}.
Hence we get that $\left\lVert u_{m_0-1}\right\rVert _{L^2}^2=0$ and thus $u_{m_0-1}=0$ for all such large enough $m_0$.
This proves the theorem in the case $n\neq2$.

Finally, let us discuss the modification needed in the case $n=2$.
In this case we instead use the second equality in (\ref{MathItem_5.45}) and hence proceeding as above arrive at that for sufficiently large $m$
$$d_ma_{m+1}\geq a_{m-1}+r_m,$$
where we have used that $d_m\geq d_{m+1}$ and that all $d_i\geq1$ for $i\geq1$.
Multiplying through by $d_{m-2}$, applying the same inequality with ``$m$" replaced by ``$m-2$" on the right-hand side, and then repeating recursively gives
$$ \left(d_m\cdot\ldots\cdot d_{m_0} \right)a_{m+1}=a_{m_0-1}+r_{m_0}+ \left(d_{m_0} \right)r_{m_0+2}+\ldots+ \left(d_{m-2}\cdot\ldots\cdot d_{m_0} \right)r_m$$
for any pair of $m,m_0>\deg{f}$ such that $m=m_0+2k$ for some integer $k$.
Since the $d_i\geq1$, we get the inequality
$$ \left(\prod_{i=0}^{k}d_{m-2i} \right)a_{m+1}\geq a_{m_0-1}+r_{m_0}+r_{m_0+2}+\ldots+r_m.$$
Again by the definition of $a_{m_0-1}$ and (\ref{MathItem_5.45}) we have that $a_{m_0}\geq \left(\sfrac{c_{m_0-1}}{d_{m_0-1}} \right)\left\lVert u_{m_0-1}\right\rVert _{L^2}^2$ and so
$$ \left(\prod_{i=0}^{k}d_{m-2i} \right)a_{m+1}\geq \left(\frac{c_{m_0-1}}{d_{m_0-1}}- \left(1+B \right) \right)\left\lVert u_{m_0-1}\right\rVert _{L^2}^2.$$
Since the $d_m\rightarrow1$ as $m\rightarrow\infty$, we can assume that we defined $m_0$ before so that $\sfrac{c_{m_0-1}}{d_{m_0-1}}\geq B+1$ as well.
Furthermore, $\prod_{m=1}^{\infty}d_m$ converges by the infinite product criteria since $\sum_{m=1}^{\infty} \left(d_m-1 \right)<\infty$ and so the coefficient on the left-hand side is bounded by some fixed constant.
Hence, as before, the theorem follows from the fact that $a_{m+1}\rightarrow0$ as $m\rightarrow\infty$.

\begin{flushright}$\blacksquare$\end{flushright}

Theorem \ref{MathItem_5.38} has one disadvantage.
Though it tells us that $f$ being of finite degree implies that the solution $u$ is also of finite degree, it gives no information about the degree of $u$ itself.
The following proposition remedies this by assuming an additional condition on the metric $g$ and connection $\nabla^\mathcal{E}$.

\begin{definition} \label{MathItem_5.54}  Suppose that $ \left(M,g \right)$ is a Riemannian manifold, $\mathcal{E}$ is a smooth complex vector bundle over $M$, and that $\nabla^\mathcal{E}$ is a smooth connection in $\mathcal{E}$.
Letting $\mathbb{X}_+ : \Omega_m\rightarrow\Omega_{m+1}$ be the operators defined in (\ref{MathItem_5.41}), for any index $m\geq0$ we call elements of $\ker{\left.\mathbb{X}_+\right|_{\Omega_m}}$ \textbf{twisted conformal Killing tensors (CKTs)} of degree $m$ of the connection $\nabla^\mathcal{E}$.
A \textbf{nontrivial twisted CKT} is a twisted CKT of degree $m\geq1$ that is not identically zero.
\end{definition}

We remark that the word ``twisted" comes up in the above definition because there is a vector bundle involved.

\begin{proposition} \label{MathItem_5.55}  Suppose that $ \left(M\subseteq\overline{M},g \right)$ is an asymptotically hyperbolic space, $ \left(\mathcal{E},\langle \cdot,\cdot\rangle _\mathcal{E} \right)$ is a smooth complex Hermitian vector bundle over $\overline{M}$, $\Phi\in C^\infty \left(\overline{M};{\End}_{\mathrm{sk}}{\mathcal{E}} \right)$, and that $\nabla^\mathcal{E}$ is a smooth unitary connection in $\mathcal{E}$.
Assume also that the sectional curvatures of $g$ are negative and that $\nabla^\mathcal{E}$ has no nontrivial twisted CKTs in $\mathcal{R}^3 \left(SM;\pi ^\ast\mathcal{E} \right)$.
If $u\in\mathcal{R}^3 \left(SM;\pi ^\ast\mathcal{E} \right)$ solves
$$\mathbb{X}u+\Phi u=f$$
for some $f\in C^\infty \left(SM;\pi ^\ast\mathcal{E} \right)$ of finite degree $m$, then $u$ is of degree $\max{\left\{m-1,0\right\}}$.
\end{proposition}

\noindent\textbf{Proof:} We know by Theorem \ref{MathItem_5.38} that $u$ has finite degree, call it $m^\prime\geq0$.
If $m^\prime=0$, then we are done.
So suppose that $m^\prime\geq1$.
We will prove this theorem by contradiction: suppose that $m^\prime\geq m$.
Then comparing the Fourier modes of order $m^\prime+1$ of both sides of the above equation gives
$$\mathbb{X}_+u_{m^\prime}=0.$$
Hence, $u_{m^\prime}$ is a twisted CKT that is in $\mathcal{R}^3 \left(SM;\pi ^\ast\mathcal{E} \right)$ by Proposition \ref{MathItem_5.36}.
Since we assumed that there are no such nontrivial CKTs, we conclude that $u_m$ is identically zero.
But this contradicts that the degree of $u$ is $m^\prime$, and hence proves the proposition.

\begin{flushright}$\blacksquare$\end{flushright}
\subsection{Regularity of Solutions to the Transport Equation} \label{Section_5.4} 

Before we prove the main result of our paper, we need to establish the regularity of solutions to transport equations of a specific form.
Here we use the material that we introduced in Section \ref{Section_2.5} up to (\ref{MathItem_2.16}) there and the two sentences after.
Below we consider sections $A\in C^\infty \left(\overline{M};\Hom{ \left(T\overline{M},\mathcal{E} \right)} \right)$ due to their role in Theorem \ref{MathItem_2.8} and its proof when we study $A={\widetilde{\nabla}}^\mathcal{E}-\nabla^\mathcal{E}$.
We point out that each $A\in C^\infty \left(\overline{M};\Hom{ \left(T\overline{M},\mathcal{E} \right)} \right)$ is canonically identified with a bundle homomorphism of the form $A : T\overline{M}\rightarrow\mathcal{E}$, whose restriction to the unit tangent bundle $A : SM\rightarrow\mathcal{E}$ we will also denote by the same letter.

The following is our regularity theorem:

\begin{proposition} \label{MathItem_5.56}  Suppose that $ \left(M\subseteq\overline{M},g \right)$ is a nontrapping asymptotically hyperbolic space, $\rho$ is a boundary defining function, $\mathcal{E}$ is a smooth complex vector bundle over $\overline{M}$, and that $\nabla^\mathcal{E}$ is a smooth connection in $\mathcal{E}$.
Suppose also that we have a $\Phi\in\rho C^\infty \left(\overline{M};\End{\mathcal{E}} \right)$, a $\phi\in\rho C^\infty \left(\overline{M};\mathcal{E} \right)$, and an $A\in C^\infty \left(\overline{M};\Hom{ \left(T\overline{M};\mathcal{E} \right)} \right)$.
Then for any given boundary data $h\in C^\infty \left(\partial _-{}^bS^\ast\overline{M};\left.\pi _b^\ast\mathcal{E}\right|_{\partial _-{}^bS^\ast\overline{M}} \right)$, there exists a unique solution $u\in C^\infty \left({}^bS^\ast\overline{M};\pi _b^\ast\mathcal{E} \right)$ to
\begin{equation} \label{MathItem_5.57} \mathbb{X}u+\Phi u=\phi+A\ \ \ \ \ \mathrm{on}  \ SM,\end{equation}
with $\left.u\right|_{\partial _-{}^bS^\ast\overline{M}}=h$.

Furthermore, for any $l\geq0$ there exists an $N\geq0$ dependent only on $ \left(M,g \right)$ and $l$ such that the following holds.
Suppose that $\Phi\in\rho^{N+1}C^\infty \left(\overline{M};\End{\mathcal{E}} \right)$, $\phi\in\rho^{N+1}C^\infty \left(\overline{M};\mathcal{E} \right)$, $A\in\rho^NC^\infty \left(\overline{M};\Hom{ \left(T\overline{M},\mathcal{E} \right)} \right)$, and that the connection symbols of $\nabla^\mathcal{E}$ are in $\rho^NC^\infty \left(\overline{M} \right)$ in any boundary coordinates and frame (in the sense of Definition \ref{MathItem_2.7}).
Then for all solutions $u$ as above that also satisfy $h\equiv0$ and $\left.u\right|_{\partial _+S^\ast\overline{M}}\equiv0$, it holds that $u\in\mathcal{R}^l \left(SM;\pi ^\ast\mathcal{E} \right)$.
\end{proposition}

\noindent\textbf{Proof:} Let $d=\rank{\mathcal{E}}$.
From Lemma 2.1 in \cite{Bibitem_15} we have that $X=\rho\overline{X}$ for some smooth vector field $\overline{X}$ over ${}^bS^\ast\overline{M}$ that is nonvanishing and transverse to $\partial {}^bS^\ast\overline{M}$.
Let $ \left(\rho,y^\mu \right)= \left(x^i \right)$ be asymptotic boundary normal coordinates of $\overline{M}$, $ \left(b_k \right)$ a frame for $\mathcal{E}$ over their domain, and consider the coordinates $\sfrac{v^i\partial }{\partial  x^i}\mapsto \left(x^i,v^i \right)$ of $TM$.
Let ${}^\mathcal{E}\Gamma_{ij}^k$ denote the connection symbols of $\nabla^\mathcal{E}$ with respect to $ \left(\sfrac{\partial }{\partial  x^i} \right)$ and $ \left(b_k \right)$.
Then observe that the components of (\ref{MathItem_5.57}) with respect to $ \left(\pi ^\ast b_k \right)$ and $ \left(\sfrac{\partial }{\partial  x^i} \right)$ are given by
\begin{equation} \label{MathItem_5.58} Xu^k+{}^\mathcal{E}\Gamma_{ij}^kv^iu^j+\Phi_j^ku^j=\phi^k+A_i^kv^i\end{equation}
where $k=1,\ldots,d$.
We point out that the $A_i^k$ are smooth on $\overline{M}$ over our coordinates' domain.
Dividing through by $\rho$ gives
\begin{equation} \label{MathItem_5.59} \overline{X}u^k+ \left(\rho^{-1} \right){}^\mathcal{E}\Gamma_{ij}^kv^iu^j+\rho^{-1}\Phi_j^ku^j=\rho^{-1}\phi^k+\rho^{-1}A_i^kv^i.\end{equation}
Take the coordinates $\eta_0\sfrac{d\rho}{\rho}+\eta_\lambda dy^\lambda\mapsto \left(x^i,\eta_i \right)$ of ${}^bT^\ast\overline{M}$ and observe that the canonical identification $\sharp \circ \left(F^\ast \right)^{-1} : \left.{}^bT^\ast\overline{M}\right|_M\rightarrow TM$ is given by
\begin{equation} \label{MathItem_5.60} v^i=g^{ii^\prime}\begin{cases}\sfrac{\eta_0}{\rho}&\mathrm{if}\ \ \ i^\prime=0\\\eta_{i^\prime}&\mathrm{if}\ \ \ i^\prime\geq1\\\end{cases}.\end{equation}
We denote the right-hand side by $g^{ii^\prime}\left\{\eta_{i^\prime}\right\}$.
Pulling the above equation (\ref{MathItem_5.59}) to ${}^bS^\ast\overline{M}$ (in the sense of Remark \ref{MathItem_4.10}) gives
\begin{equation} \label{MathItem_5.61} \overline{X}u^k+\rho^{-1}{}^\mathcal{E}\Gamma_{ij}^kg^{ii^\prime}\left\{\eta_{i^\prime}\right\}u^j+\rho^{-1}\Phi_j^ku^j=\rho^{-1}\phi^k+\rho^{-1}A_i^kg^{ii^\prime}\left\{\eta_{i^\prime}\right\}.\end{equation}
We remind the reader that each $g^{ij}$ is $\rho^2$ times something smooth on $\overline{M}$.
Hence it follows from our assumptions that all of the terms in this differential equation are smooth on ${}^bS^\ast\overline{M}$.

Since we assumed that $g$ is nontrapping, it follows from the proof of Corollary 2.5 in \cite{Bibitem_15} that any maximal integral curve\footnote{ i.e. integral curve whose interval domain cannot be extended.} $\sigma$ of $\overline{X}$ is of the form $\sigma : \left[a,b\right]\rightarrow{}^bS^\ast\overline{M}$ where $a$ and $b$ are finite with $\sigma \left(a \right)\in\partial _-{}^bS^\ast\overline{M}$ and $\sigma \left(b \right)\in\partial _+{}^bS^\ast\overline{M}$.
Hence, (\ref{MathItem_5.61}) can be viewed as ordinary differential equations (ODEs) along such curves $\sigma$.
Hence, since $\overline{X}$ is nonvanishing and transverse to $\partial {}^bS^\ast\overline{M}$, it follows from the theory of flows and the existence, uniqueness, and smooth dependence on initial condition of linear ODEs (see \cite{Bibitem_6} and \cite{Bibitem_25}) that indeed a unique smooth solution $u$ exists to (\ref{MathItem_5.61}) and hence (\ref{MathItem_5.57}) satisfying the given boundary data.

Next suppose that $h\equiv0$ and $\left.u\right|_{\partial _+S^\ast\overline{M}}\equiv0$.
If $l\geq1$, suppose that we also have smooth vector fields $V_1,...,V_l\in C^\infty \left({}^0S^\ast\overline{M};T{}^0S^\ast\overline{M} \right)$ over the 0-cosphere bundle.
Pick any point $x_0\in\partial  M$ contained in our coordinates $ \left(x^i \right)$.
We will show that for some compact neighborhood $K$ of $x_0$ (i.e. $x_0$ is in the interior of $K$)
$$\left|u \left(z \right)\right|\le\rho^{N_0} \left(z \right)L^\infty\left[\pi ^{-1}\left[K\right]\right]\ \ \ \mathrm{and} \ \ \left|V_j\ldots V_1u \left(z \right)\right|\le\rho^{N_j} \left(z \right)L^\infty\left[\pi ^{-1}\left[K\right]\right]$$
over $z\in\pi ^{-1}\left[K\right]$ where $j=1,\ldots,l$ and $N_j\geq\sfrac{ \left(n+1 \right)}{2}$ are constants.
From this and (\ref{MathItem_5.2}) it will follow that $u\in\mathcal{R}^l$.
Throughout the proof we will make $N\geq0$ is as big as we need whenever we need it.
At the end of the proof, we will discuss why there exists a maximum upper bound on the size of $N$ that we need that is dependent only on $ \left(M,g \right)$ and $l$.

We begin by establishing a few facts about the flow of $X$.
Let $\varphi : SM\times\mathbb{R}\rightarrow SM$ denote the flow of $X$.
For any point $\zeta=\eta_0\sfrac{d\rho}{\rho}+\eta_\lambda dy^\lambda\in{}^bS^\ast M$ we write its identified point on $SM$ as $z\in SM$.
In \cite{Bibitem_15} the authors explain that close enough to the boundary, the flow $\varphi$ moves away from the boundary when $\eta_0>0$ and towards the boundary when $\eta_0<0$.
More precisely, by Lemma 2.3 in \cite{Bibitem_15} and its proof there exists constants $C,\varepsilon>0$ such that if we take any $z\in SM\cap\left\{\rho<\varepsilon\right\}$ and write $\rho \left(t \right)=\rho\circ\varphi_z \left(t \right)$,

	\begin{enumerate} \item if $\eta_0\geq0$ then $\lim_{t\rightarrow-\infty}{\varphi_z \left(t \right)}\in\partial _-{}^bS^\ast\overline{M}$ and for $t\le0$ we have that $\rho \left(t \right)$ is both increasing and satisfies $\rho \left(t \right)\le Ce^t$,

	\item if $\eta_0<0$ then $\lim_{t\rightarrow+\infty}{\varphi_z \left(t \right)}\in\partial _+{}^bS^\ast\overline{M}$ and for $t\geq0$ we have that $\rho \left(t \right)$ is both decreasing and satisfies $\rho \left(t \right)\le Ce^{-t}$.\end{enumerate}

By the same lemma and its proof, it follows that there exist compact neighborhoods $K,K^\prime\subseteq\left\{\rho<\varepsilon\right\}$ of $x_0$ in $\overline{M}$ satisfying $K\subseteq K^\prime\subseteq\dom{ \left(x^i \right)}$ such that for any $z\in\pi ^{-1}\left[K\right]$, $\varphi_z$ will always be contained in $\pi ^{-1}\left[K^\prime\right]$ for $t\le0$ if $\eta_0\geq0$ and $t\geq0$ if $\eta_0<0$.

Fix some $N_0\geq\sfrac{ \left(n+1 \right)}{2}$ whose size will be increased later if needed.
Our first goal is to show that each $u^k\in\rho^{N_0}L^\infty\left[\pi ^{-1}\left[K\right]\right]$ where our approach will be to study the growth of the solution to (\ref{MathItem_5.58}) by writing that equation as an ODE along integral curves of $X$.
Fix any $\zeta\cong z\in\pi ^{-1}\left[K\right]$.
Suppose that $\eta_0\geq0$ since the proof below is essentially the same for the case $\eta_0<0$.
We set $B$ and $b$ to be the $d\times d$ matrix and $d\times1$ column vector respectively given by
\begin{equation} \label{MathItem_5.62} B_i^k= \left({}^\mathcal{E}\Gamma_{ij}^kv^i+\Phi_j^k \right)\ \ \ \ \ \ \ \ \mathrm{and} \ \ \ \ \ \ \ b^k=\phi^k+A_i^kv^i,\end{equation}
and we denote $B \left(t \right)=B\circ\varphi_z \left(t \right)$, $b \left(t \right)=b\circ\varphi_z \left(t \right)$, and $u \left(t \right)=u\circ\varphi_z \left(t \right)$.
This way, (\ref{MathItem_5.58}) along $\varphi_z$ is given by the ODE
\begin{equation} \label{MathItem_5.63} \frac{du^k}{dt}+B_i^ku^i=b^k.\end{equation}
Since $\left.u\right|_{\partial {}^bS^\ast\overline{M}}\equiv0$ and $\lim_{t\rightarrow-\infty}{\varphi_z \left(t \right)}\in\partial _-{}^bS^\ast\overline{M}$, we have that our solution $u$ satisfies the ``initial condition" $\lim_{t\rightarrow-\infty}{u \left(t \right)}=0$.
Now, for any matrix $M$ or column vector $w$, let $\left|M\right|$ and $\left|w\right|$ denote the norms $\sum_{ik}\left|M_i^k\right|$ and $\sum_{i}\left|w^i\right|$.
Fixing any time $t_0<0$, we get from (\ref{MathItem_5.63}), the fundamental theorem of calculus, and the triangle inequality that
\begin{equation} \label{MathItem_5.64} \left|u \left(t \right)\right|\le\left|u \left(t_0 \right)\right|+\int_{t_0}^{t}\left|B \left(s \right)\right|\left|u \left(s \right)\right|ds+\int_{t_0}^{t}\left|b \left(s \right)\right|ds.\end{equation}
We will use this to bound $\left|u \left(0 \right)\right|=\left|u \left(z \right)\right|$ in terms of a power of $\rho$.
We employ the standard technique in ODEs of defining the function $R :  \left(-\infty,0\right]\rightarrow\mathbb{R}$ given by the first integral on the right-hand side above.
Whenever $\left|B \left(t \right)\right|\neq0$,
$$\frac{1}{\left|B \left(t \right)\right|}R^\prime \left(t \right)\le\left|u \left(t_0 \right)\right|+R \left(t \right)+\int_{t_0}^{t}\left|b \left(s \right)\right|ds$$
and so
$$R^\prime \left(t \right)\le\left|B \left(t \right)\right|R \left(t \right)+\left|B \left(t \right)\right| \left(\left|u \left(t_0 \right)\right|+\int_{t_0}^{t}\left|b \left(s \right)\right|ds \right).$$
Observe that this inequality also holds when $\left|B \left(t \right)\right|=0$ and hence for all $t\le0$.
This resembles a separable ODE in $R$ and $t$ and hence can be treated similarly.
In particular, if we take $\left|B \left(t \right)\right|R \left(t \right)$ to the left-hand side, multiply through by $\exp{\left[-\int_{t_0}^{t}\left|B \left(s \right)\right|ds\right]}$, integrate from $t=t_0$ to $t=0$ (using that $R \left(t_0 \right)=0$), and finally divide through by $\exp{\left[-\int_{t_0}^{0}\left|B \left(s \right)\right|ds\right]}$, we will get
\begin{equation} \label{MathItem_5.65} R \left(0 \right)=\int_{t_0}^{0}\left|B \left(s \right)\right|\left|u \left(s \right)\right|ds\end{equation}
$$\le e^{\int_{t_0}^{0}\left|B \left(s \right)\right|ds}\int_{t_0}^{0}{e^{-\int_{t_0}^{t}\left|B \left(s \right)\right|ds}\left|B \left(t \right)\right| \left(\left|u \left(t_0 \right)\right|+\int_{t_0}^{t}\left|b \left(s \right)\right|ds \right)dt}.$$
Before we let $t_0\rightarrow-\infty$, let us discuss integrability.
Fix $N_0\geq\sfrac{ \left(n+1 \right)}{2}$ whose size will be increased later if needed.
Because we are working on the sphere bundle and $g=\rho^{-2}\overline{g}$ for a smooth metric $\overline{g}$, there exists a constant $C_2>0$ such that the magnitude of $v^i$ in the expressions for $B$ and $b$ in (\ref{MathItem_5.62}) are less than $C_2\rho$ on $\pi ^{-1}\left[K^\prime\right]$.
Hence by (\ref{MathItem_5.62}) we have that there exist constants $C_3,C_{N_0}>0$ such that $\left|B\right|\le C_3\rho$ and $\left|b\right|\le C_{N_0}\rho^{N_0+1}$ on $\pi ^{-1}\left[K^\prime\right]$ (i.e. $C_{N_0}$ depends on $N_0$).
In particular, since $\rho \left(t \right)\le Ce^t$ for $t\le0$ we have that all of the integrals on the right-hand side of (\ref{MathItem_5.65}) converge if we let $t_0\rightarrow-\infty$.
Hence letting $t_0\rightarrow-\infty$ in (\ref{MathItem_5.65}) gives that
\begin{equation} \label{MathItem_5.66} \int_{-\infty}^{0}\left|B \left(s \right)\right|\left|u \left(s \right)\right|ds\le e^{\int_{-\infty}^{0}\left|B \left(s \right)\right|ds}\int_{-\infty}^{0}{\left|B \left(t \right)\right|\int_{-\infty}^{t}\left|b \left(s \right)\right|dsdt}\end{equation}
where we have used that $\lim_{t\rightarrow-\infty}{u \left(t \right)}=0$ and that $\exp{\left[-\int_{t_0}^{t}\left|B \left(s \right)\right|ds\right]}\le1$.
To estimate the right-hand side further, observe that since $\rho \left(t \right)$ is increasing on $t\le0$ we have that $\left|b \left(t \right)\right|\le C_{N_0}\rho \left(t \right)\rho^{N_0} \left(0 \right)\le C_{N_0}Ce^t\rho^{N_0} \left(0 \right)$ and so
\begin{equation} \label{MathItem_5.67} \int_{-\infty}^{0}\left|b \left(s \right)\right|ds\le C_{N_0}C\rho^{N_0} \left(0 \right).\end{equation}
Doing a similar sort of thing for $\left|B \left(t \right)\right|$, we get from (\ref{MathItem_5.66}) that
\begin{equation} \label{MathItem_5.68} \int_{-\infty}^{0}\left|B \left(s \right)\right|\left|u \left(s \right)\right|ds\le e^{C_3C}C_3CC_{N_0}C\rho^{N_0} \left(0 \right).\end{equation}
Finally, letting $t_0\rightarrow-\infty$ in (\ref{MathItem_5.64}) and plugging (\ref{MathItem_5.67}) and (\ref{MathItem_5.68}) into there finally gives us that for some constant $C_0^\prime>0$,
\begin{equation} \label{MathItem_5.69} \left|u \left(z \right)\right|\le C_0^\prime\rho^{N_0} \left(z \right)\end{equation}
for all $z\in\pi ^{-1}\left[K\right]$, where we have used that $u \left(0 \right)=u \left(z \right)$.
Hence indeed each $u^k\in\rho^{N_0}L^\infty\left[\pi ^{-1}\left[K\right]\right]$.

Next, suppose $l\geq1$ and fix $N_1\geq\sfrac{ \left(n+1 \right)}{2}$ whose size will be increased later if needed.
We show that $V_1u^k\in\rho^{N_1}L^\infty\left[\pi ^{-1}\left[K\right]\right]$, after which it should be clear how the cases of the higher order derivatives are handled (i.e. $V_2V_1u^k$, ...) if $l\geq2$.
We do this by applying $V_1$ to the ODE (\ref{MathItem_5.63}) and studying the growth of the solution.
Applying $V_1$ to (\ref{MathItem_5.63}) and rearranging, we obtain
\begin{equation} \label{MathItem_5.70} \frac{d\left[V_1 \left(u^k\circ\varphi \right)\right]}{dt}+ \left(B_i^k\circ\varphi \right)V_1 \left(u^k\circ\varphi \right)=V_1 \left(b^k\circ\varphi \right)-V_1 \left(B_i^k\circ\varphi \right) \left(u^i\circ\varphi \right).\end{equation}
Considering that in form this is a similar ODE for $V_1 \left(u^k\circ\varphi \right)$ as (\ref{MathItem_5.63}) is for $u^k$, up to a few details we describe below, a similar proof as for (\ref{MathItem_5.69}) shows that for some constant $C_1^\prime>0$
$$\left|V_1 \left(u\circ\varphi \right) \left(z,0 \right)\right|=\left|V_1u \left(z \right)\right|\le C_1^\prime\rho^{N_1} \left(z \right)$$
for all $z\in\pi ^{-1}\left[K\right]$ and hence indeed each $V_1u^k\in\rho^{N_1}L^\infty\left[\pi ^{-1}\left[K\right]\right]$.
The only analogous steps that we are missing and need to show is that $V_1 \left(u^k\circ\varphi \right)$ satisfies the ``initial condition"
\begin{equation} \label{MathItem_5.71} \lim_{t\rightarrow-\infty}{V_1 \left(u^k\circ\varphi \right)}=0\end{equation}
and that there exists a constant $\widetilde{C}_{N_1}>0$ such that the right-hand side of (\ref{MathItem_5.70}) satisfies
\begin{equation} \label{MathItem_5.72} \left|V_1 \left(b^k\circ\varphi \right)-V_1 \left(B_i^k\circ\varphi \right) \left(u^i\circ\varphi \right)\right|\le \widetilde{C}_{N_1}\rho^{N_1+1}\circ\varphi\end{equation}
for all $z\in\pi ^{-1}\left[K^\prime\right]$.

By shrinking everything if needed, we assume that there is another compact subset $K^{\prime\prime}$ of our coordinates' domain such that $K^\prime\subseteq K^{\prime\prime}\subseteq\dom{ \left(x^i \right)}$ and such that for any $z\in\pi ^{-1}\left[K^\prime\right]$, $\varphi_z$ will always be contained in $\pi ^{-1}\left[K^{\prime\prime}\right]$ for $t\le0$ if $\eta_0\geq0$ and $t\geq0$ if $\eta_0<0$.
Throughout the rest of the proof, we again assume that we are working with points such that $\eta_0\geq0$ because the proof is similar for the case $\eta_0<0$.
We begin with showing the ``initial condition" (\ref{MathItem_5.71}).
Let $ \left(\varphi^i,\varphi_i \right)$ denote the components of $\varphi$ in the coordinates $ \left(x^i,\eta_i \right)$ of ${}^bT^\ast\overline{M}$.
Then over the interior of ${}^bS^\ast\overline{M}$, we have by the chain rule that
\begin{equation} \label{MathItem_5.73} V_1 \left(u^k\circ\varphi \right)= \left(\frac{\partial  u^k}{\partial  x^i}\circ\varphi \right)V_1\varphi^i+ \left(\frac{\partial  u^k}{\partial \eta_i}\circ\varphi \right)V_1\varphi_i.\end{equation}
Since $\left.u\right|_{\partial {}^bS^\ast\overline{M}}\equiv0$, and by (\ref{MathItem_5.69}) $u$ vanishes like $\rho^{N_0}$ at the boundary $\partial {}^bS^\ast\overline{M}$ where we can require $N_0\geq2$, all of $u^k$'s partials vanish at the boundary like $\rho^{N_0-1}$.
In particular, for fixed $z\in\pi ^{-1}\left[K\right]$,
\begin{equation} \label{MathItem_5.74}  \left(\frac{\partial  u^k}{\partial  x^i}\circ\varphi \right), \left(\frac{\partial  u^k}{\partial \eta_i}\circ\varphi \right)\in O \left(e^{ \left(N_0 -1 \right)t} \right)\ \ \ \ \ \mathrm{as}  \ t\rightarrow-\infty.\end{equation}
So let us take a look at size of the terms $V_1\varphi^i$ and $V_1\varphi_i$ in (\ref{MathItem_5.73}).
Consider the coordinates ${\overline{\eta}}_i\sfrac{dx^i}{\rho}\mapsto \left(x^i,{\overline{\eta}}_i \right)$ of ${}^0T^\ast\overline{M}$ (note the bars to distinguish these from our coordinates of ${}^bT^\ast\overline{M}$ above) and let $ \left({\widetilde{\varphi}}^i,{\widetilde{\varphi}}_i \right)$ denote the components of $\varphi$ with respect to these coordinates.
We note that $ \left({\widetilde{\varphi}}^i,{\widetilde{\varphi}}_i \right)$ extends smoothly to the boundary $\partial {}^0S^\ast\overline{M}$ since it is the flow of $X$ and the latter is smooth on ${}^0S^\ast\overline{M}$ by (\ref{MathItem_5.13}).
Since canonical identification is given by $x^i=x^i$, $\eta_0={\overline{\eta}}_0$ and $\eta_\lambda=\rho^{-1}{\overline{\eta}}_\lambda$ for $\lambda=1,\ldots,n$, we have that by identification $\varphi^i={\widetilde{\varphi}}^i$, $\varphi_0={\widetilde{\varphi}}_0$, and $\varphi_\lambda= \left({\widetilde{\varphi}}^0 \right)^{-1}{\widetilde{\varphi}}_\lambda$ for $\lambda=1,\ldots,n$ where we point out that ${\widetilde{\varphi}}^0=\rho\circ\widetilde{\varphi}$.
Hence
\begin{equation} \label{MathItem_5.75} V_1\varphi^i=V_1{\widetilde{\varphi}}^i,\ \ \ \ \ \ \ \ \ \ V_1\varphi_0=V_1{\widetilde{\varphi}}_0,\end{equation}
\begin{equation} \label{MathItem_5.76} V_1\varphi_\lambda=- \left(\rho\circ\varphi \right)^{-2}V_1 \left({\widetilde{\varphi}}^0 \right){\widetilde{\varphi}}_\lambda+ \left(\rho\circ\varphi \right)^{-1}V_1{\widetilde{\varphi}}_\lambda\ \ \ \ \ \mathrm{for} \ \ \lambda=1,\ldots,n.\end{equation}
We point out that ${\widetilde{\varphi}}_\lambda$ are bounded over $\pi ^{-1}\left[K^{\prime\prime}\right]$ because they are smooth over the compact $\pi _0^{-1}\left[K^{\prime\prime}\right]$.
Now, we showed in (\ref{MathItem_5.13}) that $X$ extends as a smooth vector field on ${}^0S^\ast\overline{M}$.
Hence, as is noted at the end of the proof of Lemma 3.13 in \cite{Bibitem_15}, $X$ has a Lipschitz constant uniformly bounded on the compact ${}^0S^\ast\overline{M}$.
Since $\widetilde{\varphi}$ is the flow of $X$ over ${}^0S^\ast\overline{M}$, as the authors explain there it follows by Gr\"onwall's inequality that there exist constants $C_4,c_4>0$ such that all of the partials
\begin{equation} \label{MathItem_5.77} \left|\frac{\partial {\widetilde{\varphi}}^i}{\partial  x^i}\right|,\left|\frac{\partial {\widetilde{\varphi}}^i}{\partial {\overline{\eta}}_i}\right|,\left|\frac{\partial {\widetilde{\varphi}}_i}{\partial  x^i}\right|,\left|\frac{\partial {\widetilde{\varphi}}_i}{\partial {\overline{\eta}}_i}\right|\le C_4e^{-c_4t}\end{equation}
on $ \left(z,t \right)$ for all $z\in\pi ^{-1}\left[K^\prime\right] : \eta_0\geq0$ and $t\le0$.
Next, if we write in coordinates $V_1= \left(V_1 \right)^i\sfrac{\partial }{\partial  x^i}+ \left(V_1 \right)_i\sfrac{\partial }{\partial {\overline{\eta}}_i}$ over ${}^0S^\ast\overline{M}$, we have that the components $ \left(V_1 \right)^i$ and $ \left(V_1 \right)_i$ are also bounded over $\pi ^{-1}\left[K^\prime\right]$ because they are smooth over the compact $\pi _0^{-1}\left[K^\prime\right]$.
Hence there exists a constant $C_5\geq0$ such that
\begin{equation} \label{MathItem_5.78} \left|V_1{\widetilde{\varphi}}^i\right|\le C_5e^{-c_4t}\ \ \ \mathrm{and} \ \ \left|V_1{\widetilde{\varphi}}_i\right|\le C_5e^{-c_4t}\end{equation}
on $ \left(z,t \right)$ for all $z\in\pi ^{-1}\left[K^\prime\right] : \eta_0\geq0$ and $t\le0$.
Plugging this into (\ref{MathItem_5.75}) and (\ref{MathItem_5.76}), and then into (\ref{MathItem_5.73}), we see that choosing $N_0$ big enough in (\ref{MathItem_5.74}) will make the ``initial conditions" (\ref{MathItem_5.71}) hold (in particular, any $N_0>c_4+1$ will work).

Finally, let us show that (\ref{MathItem_5.72}) holds.
By (\ref{MathItem_5.69}), we already have that $\left| \left(u^i\circ\varphi \right)\right|\le C_N^\prime\rho^{N_0}\circ\varphi$.
Next, the canonical identification $TM\cong{}^0T^\ast\overline{M}$ is given by $v^i=\rho^{-1}g^{ii^\prime}{\overline{\eta}}_{i^\prime}$ and so by (\ref{MathItem_5.62}) we have that over ${}^0S^\ast\overline{M}$
$$B_i^k= \left({}^\mathcal{E}\Gamma_{ij}^k\rho^{-1}g^{ii^\prime}{\overline{\eta}}_{i^\prime}+\Phi_j^k \right)\ \ \ \ \ \ \ \ \mathrm{and} \ \ \ \ \ \ \ b^k=\phi^k+A_i^k\rho^{-1}g^{ii^\prime}{\overline{\eta}}_{i^\prime}.$$
In particular, we see that these are $\rho^{N+1}$ times something smooth on ${}^0S^\ast\overline{M}$.
Hence, the partials $\sfrac{\partial  B_i^k}{\partial  x^l}$, $\sfrac{\partial  B_i^k}{\partial {\overline{\eta}}_l}$, $\sfrac{\partial  b^k}{\partial  x^l}$, and $\sfrac{\partial  b^k}{\partial {\overline{\eta}}_l}$ are all $\rho^N$ times something smooth on ${}^0S^\ast\overline{M}$.
In particular,
\begin{equation} \label{MathItem_5.79} \left|\frac{\partial  B_i^k}{\partial  x^l}\circ\varphi\right|,\left|\frac{\partial  B_i^k}{\partial {\overline{\eta}}_l}\circ\varphi\right|,\left|\frac{\partial  b^k}{\partial  x^l}\circ\varphi\right|,\left|\frac{\partial  b^k}{\partial {\overline{\eta}}_l}\circ\varphi\right|\le C_6\rho^N\circ\varphi\end{equation}
on $ \left(z,t \right)$ for all $z\in\pi ^{-1}\left[K^\prime\right] : \eta_0\geq0$ and $t\le0$.
Now, over ${}^0S^\ast\overline{M}$ we have that
\begin{equation} \label{MathItem_5.80} \begin{array}{r@{}l}V_1 \left(B_i^k\circ\varphi \right)= \left(\frac{\partial  B_i^k}{\partial  x^l}\circ\varphi \right)V_1{\widetilde{\varphi}}^l+ \left(\frac{\partial  B_i^k}{\partial {\overline{\eta}}_l}\circ\varphi \right)V_1{\widetilde{\varphi}}_l,\\V_1 \left(b^k\circ\varphi \right)= \left(\frac{\partial  b^k}{\partial  x^l}\circ\varphi \right)V_1{\widetilde{\varphi}}^l+ \left(\frac{\partial  b^k}{\partial {\overline{\eta}}_l}\circ\varphi \right)V_1{\widetilde{\varphi}}_l.\\\end{array}\end{equation}
Hence if in (\ref{MathItem_5.79}) we write $\rho^N\circ\varphi\le Ce^{c_4t} \left(\rho^{N-c_4}\circ\varphi \right)$ it follows from (\ref{MathItem_5.78}) and (\ref{MathItem_5.80}) that requiring $N\geq c_4+N_1+1$ will make (\ref{MathItem_5.72}) hold.
As discussed above, this completes the proof that $V_1u^k\in\rho^{N_1}L^\infty\left[\pi ^{-1}\left[K\right]\right]$.

As we mentioned above, if $l\geq2$ the proof that $V_2V_1u^k\in\rho^{N_3}L^\infty\left[\pi ^{-1}\left[K\right]\right]$ for some $N_3\geq\sfrac{ \left(n+1 \right)}{2}$ and similarly for the higher derivatives follows similarly, starting with applying $V_2$ to (\ref{MathItem_5.70}).
We end the proof with a discussion of how we can ensure that the $N\geq0$ that we used has an upper bound dependent only on $ \left(M,g \right)$ and $l$.
Throughout the proof we required that $N$ (which determines the size of $N_0$, $N_1$,...) is bigger than fixed numbers (e.g. 2), numbers dependent on the dimension (e.g. $\sfrac{ \left(n+1 \right)}{2}$), the constant $c_4$ in (\ref{MathItem_5.77}), and the analogs of $c_4$ in (\ref{MathItem_5.77}) when carrying out the proof for higher order derivatives (i.e. $V_2V_1u^k$, ...).
We explain why the latter two can be bounded.
Due to the compactness of the boundary $\partial \overline{M}$, we can cover $\partial \overline{M}$ with a finite collection of neighborhoods of the form $K$, $K^\prime$, and $K^{\prime\prime}$ as described above which in turn determine their own $c_4$'s in (\ref{MathItem_5.77}) and their higher derivative analogs using only the geometry of the geodesic flow $\varphi$.
Hence, taking the maximum of all such $c_4$'s and their higher derivative analogs will provide us with an upper bound for $N$ that works throughout the whole proof.
In particular, this upper bound for $N$ will be dependent only on $ \left(M,g \right)$ and $l$.

\begin{flushright}$\blacksquare$\end{flushright}
\subsection{Proof of Theorem \ref{MathItem_2.8}} \label{Section_5.5} 

In our proof, we will assume that the $N\geq0$ in the statement of the theorem is as big as we need.
It will be clear that the size of $N$ that we need is dependent only on $ \left(M,g \right)$.
The first step is to provide a formulation of (\ref{MathItem_2.4}) as a single transport equation of endomorphism fields over $SM$.
To begin, take coordinates $ \left(x^i \right)$ of $\overline{M}$, a frame $ \left(b_i \right)$ for $\mathcal{E}$ over their domain, and let ${}^\mathcal{E}\Gamma_{ij}^k$ denote the connection symbols of $\nabla^\mathcal{E}$ with respect to $ \left(\sfrac{\partial }{\partial  x^i} \right)$ and $ \left(b_k \right)$.
Recall from (\ref{MathItem_4.26}) that in these coordinates
$$\nabla_v^{\End{\mathcal{E}}}U=vU+ \left({}^\mathcal{E}\Gamma \right)U-U \left({}^\mathcal{E}\Gamma \right)$$
and that we denote $\mathbb{X}U=\nabla_X^{\pi ^\ast\End{\mathcal{E}}}U$.
It is easy to check from here that this connection is unitary with respect to $\langle \cdot,\cdot\rangle _{\End{\mathcal{E}}}$ (see (\ref{MathItem_4.27})) whose connection symbols are in $\rho^NC^\infty \left(\overline{M} \right)$ in any boundary coordinates and frame (in the sense of Definition \ref{MathItem_2.7}).

Next, by looking in local coordinates and an orthonormal frame for $\mathcal{E}$, it follows that ${\widetilde{\nabla}}^\mathcal{E}=\nabla^\mathcal{E}+A$ where $A\in\rho^NC^\infty \left(\overline{M};\Hom{ \left(T\overline{M},\mathcal{E} \right)} \right)$.
Consider solutions $U$ and $\widetilde{U}$ to the following transport equations on $SM$:
\begin{equation} \label{MathItem_5.81} \begin{cases}\mathbb{X}U+\Phi U=0,&\left.U\right|_{\partial _-{}^bS^\ast\overline{M}}=\mathrm{id} ,\\\mathbb{X}\widetilde{U}+A\widetilde{U}+\widetilde{\Phi}\widetilde{U}=0,&\left.\widetilde{U}\right|_{\partial _-{}^bS^\ast\overline{M}}=\mathrm{id} ,\\\end{cases}\end{equation}
whose existence we now justify.
It is not hard to check using (\ref{MathItem_4.9}) that in the second equation ``$\mathbb{X}\widetilde{U}+A\widetilde{U}$" is given by ``$ \left(\pi ^\ast\nabla \right)_X\widetilde{U}$" for the connection $\nabla=\nabla^{\End{\mathcal{E}}}+A$.
Thus, by thinking of $U\mapsto\Phi U$ and $\widetilde{U}\mapsto\widetilde{\Phi}\widetilde{U}$ as elements of $\rho^{N+1}C^\infty \left(\overline{M};\End{\End{\mathcal{E}}} \right)$, by plugging $\End{\mathcal{E}}$ into $\mathcal{E}$ in Proposition \ref{MathItem_5.56} we get that the solutions to the above two equations indeed exist and are unique.

We demonstrate the usefulness of $U$ and $\widetilde{U}$.
Suppose that $\gamma :  \left(-\infty,\infty \right)\rightarrow M$ is a complete unit-speed geodesic and that $u :  \left(-\infty,\infty \right)\rightarrow\mathcal{E}$ is the smooth solution along $\gamma$ to the initial value problem
\begin{equation} \label{MathItem_5.82} \nabla_{\dot{\gamma} \left(t \right)}^\mathcal{E}u \left(t \right)=0,\ \ \ \ \lim_{t\rightarrow-\infty}{u \left(\gamma \left(t \right) \right)}=e,\end{equation}
where $e$ is an element in $\mathcal{E}_{x_0}$ and $x_0\in\partial \overline{M}$ is the limit of $\gamma \left(t \right)$ as $t\rightarrow-\infty$.
We point out that such a $u$ exists by Lemma \ref{MathItem_2.6} part 1) with $\Phi=0$.

Let $\sigma :  \left(-\infty,\infty \right)\rightarrow SM$ be the integral curve of $X$ satisfying $\gamma=\pi \circ\sigma$ and let $\hat{u}$ denote $u$ lifted to $\sigma$ via the canonical identification of $\mathcal{E}$ and $\pi ^\ast\mathcal{E}$.
We claim that $ \left(U\circ\sigma \right)u \left(t \right)$ and $ \left(\widetilde{U}\circ\sigma \right)u \left(t \right)$ are solutions to (\ref{MathItem_2.4}) and (\ref{MathItem_2.4}) with $\Phi$ and $\nabla^\mathcal{E}$ replaced by $\widetilde{\Phi}$ and ${\widetilde{\nabla}}^\mathcal{E}$ respectively.
It follows from (\ref{MathItem_4.8}) that $\nabla_{\dot{\sigma}}^{\pi ^\ast\mathcal{E}}\hat{u}\cong\nabla_{\dot{\gamma}}^\mathcal{E}u=0$.
Thus
\begin{equation} \label{MathItem_5.83} \nabla_{\dot{\gamma}}^\mathcal{E} \left( \left(U\circ\sigma \right)u \right)+\Phi \left( \left(U\circ\sigma \right)u \right)=\left[\nabla_{\dot{\gamma}}^{\End{\mathcal{E}}} \left(U\circ\sigma \right)\right]u+ \left(U\circ\sigma \right)\nabla_{\dot{\gamma}}^\mathcal{E}u+\Phi \left( \left(U\circ\sigma \right)u \right)\end{equation}
$$\cong\left[\nabla_X^{\pi ^\ast\End{\mathcal{E}}}U\right]\hat{u}+U\nabla_{\dot{\sigma}}^{\pi ^\ast\mathcal{E}}\hat{u}+\Phi \left(U\hat{u} \right)=0,\ \ \ \mathrm{with}\ \ \ \lim_{t\rightarrow-\infty}{ \left(U\circ\sigma \right)u \left(t \right)}=e,$$
where $\nabla_{\dot{\gamma}}^{\End{\mathcal{E}}} \left(U\circ\sigma \right)\cong\nabla_X^{\pi ^\ast\End{\mathcal{E}}}U$ follows from a local coordinate calculation using (\ref{MathItem_4.9}).
Similarly
\begin{equation} \label{MathItem_5.84} {\widetilde{\nabla}}_{\dot{\gamma}}^\mathcal{E} \left( \left(\widetilde{U}\circ\sigma \right)u \right)+\widetilde{\Phi} \left( \left(\widetilde{U}\circ\sigma \right)u \right)\end{equation}
$$=\nabla_{\dot{\gamma}}^\mathcal{E} \left( \left(\widetilde{U}\circ\sigma \right)u \right)+ \left(A\circ\sigma \right) \left(\widetilde{U}\circ\sigma \right)u+\widetilde{\Phi} \left( \left(\widetilde{U}\circ\sigma \right)u \right)$$
$$=\left[\nabla_{\dot{\gamma}}^{\End{\mathcal{E}}} \left(\widetilde{U}\circ\sigma \right)\right]u+ \left(\widetilde{U}\circ\sigma \right)\nabla_{\dot{\gamma}}^\mathcal{E}u+ \left(A\circ\sigma \right) \left(\widetilde{U}\circ\sigma \right)u+\widetilde{\Phi} \left( \left(\widetilde{U}\circ\sigma \right)u \right)$$
$$\cong \left(\nabla_X^{\pi ^\ast\End{\mathcal{E}}}\widetilde{U} \right)\hat{u}+\widetilde{U}\nabla_{\dot{\sigma}}^{\pi ^\ast\mathcal{E}}\hat{u}+A \left(\widetilde{U}\hat{u} \right)+\widetilde{\Phi} \left(U\hat{u} \right)=0,\ \ \ \mathrm{with}\ \ \ \lim_{t\rightarrow-\infty}{ \left(\widetilde{U}\circ\sigma \right)u \left(t \right)}=e.$$
Hence our claim above is indeed true.
By our assumption the data (\ref{MathItem_2.5}) is the same for $ \left(\nabla^\mathcal{E},\Phi \right)$ and $({\widetilde{\nabla}}^\mathcal{E},\widetilde{\Phi})$ and so
$$\lim_{t\rightarrow\infty}{ \left(U\circ\sigma \right)u \left(t \right)}=\lim_{t\rightarrow\infty}{ \left(\widetilde{U}\circ\sigma \right)u \left(t \right)}.$$
By Lemma \ref{MathItem_2.6} part 2) parallel transport such as (\ref{MathItem_5.82}) above is an isomorphism between fibers.
Hence varying the ``$e$" in (\ref{MathItem_5.82}) implies here that $U=\widetilde{U}$ on $\partial _+{}^bS^\ast\overline{M}$ as well (i.e. they already agree on $\partial _-{}^bS^\ast\overline{M}$ by definition).

Intuitively speaking, we have demonstrated that knowing the parallel transport (\ref{MathItem_5.82}), the endomorphism fields $U$ and $\widetilde{U}$ encode the transform that takes all possible $ \left(\gamma,e \right)$ to the data (\ref{MathItem_2.5}) of $\nabla^\mathcal{E},\Phi$ and ${\widetilde{\nabla}}^\mathcal{E},\widetilde{\Phi}$ respectively.
Furthermore, the assumption that the two transforms are equal gives us that $U$ and $\widetilde{U}$ are equal on the boundary $\partial {}^bS^\ast\overline{M}$.
Hence we have reformulated our task to showing that $U=\widetilde{U}$ on $\partial {}^bS^\ast M$ implies the gauge equivalence stated in the theorem.

Guided by the observation (\ref{MathItem_2.13}) in the introduction, we next study the behavior of $U{\widetilde{U}}^{-1}$ over the interior $SM$.
We note that both $U$ and $\widetilde{U}$ are invertible because if one rewrites (\ref{MathItem_5.81}) as in (\ref{MathItem_5.59}), we get that they satisfy a matrix ordinary differential equations with an invertible initial value (e.g. see (1.8) in Chapter 3 of \cite{Bibitem_6}).
A quick computation using (\ref{MathItem_4.26}) and (\ref{MathItem_4.9}) shows that $\mathbb{X} \left(\mathrm{id} \right)\equiv0$ and that $\mathbb{X} \left(W_1W_2 \right)=\mathbb{X} \left(W_1 \right)W_2+W_1\mathbb{X} \left(W_2 \right)$ for all $W_1,W_2\in C^\infty \left(SM;\pi ^\ast\End{\mathcal{E}} \right)$.
Hence by (\ref{MathItem_5.81}),
$$\mathbb{X} \left({\widetilde{U}}^{-1} \right)=-{\widetilde{U}}^{-1}\mathbb{X} \left(\widetilde{U} \right){\widetilde{U}}^{-1}={\widetilde{U}}^{-1} \left(A\widetilde{U}+\widetilde{\Phi}\widetilde{U} \right){\widetilde{U}}^{-1}={\widetilde{U}}^{-1}A+{\widetilde{U}}^{-1}\widetilde{\Phi}.$$
Next, we have that
$$\mathbb{X} \left(U{\widetilde{U}}^{-1} \right)= \left(-\Phi U \right){\widetilde{U}}^{-1}+U \left({\widetilde{U}}^{-1}A+{\widetilde{U}}^{-1}\widetilde{\Phi} \right),$$
and so we finally arrive at that $Q=U{\widetilde{U}}^{-1}$ satisfies
\begin{equation} \label{MathItem_5.85} \mathbb{X}Q+\Phi Q-QA-Q\widetilde{\Phi}=0.\end{equation}
This is a transport equation over $SM$.
However, to apply our finite degree theorems from Section \ref{Section_5.3} above, we need our solution to vanish at ``infinity." Hence we instead consider $W=Q-\mathrm{id}$ which satisfies
\begin{equation} \label{MathItem_5.86} \mathbb{X}W+\Phi W-WA-W\widetilde{\Phi}=-\Phi+A+\widetilde{\Phi},\end{equation}
and vanishes on the boundary $\partial {}^bS^\ast\overline{M}$.
It is an elementary exercise to check that $W\mapsto\Phi W-W\widetilde{\Phi}$ is an element of $\rho^{N+1}C^\infty \left(\overline{M};{\End}_{\mathrm{sk}}{\End{\mathcal{E}}} \right)$ which we shall call $\Psi$ (note the ``sk" with respect to $\langle \cdot,\cdot\rangle _{\End{\mathcal{E}}}$).
Similarly, $\mathbb{X}W-WA$ is given by $ \left(\pi ^\ast\nabla^\prime \right)_XW$ for the unitary connection $\nabla^\prime U=\nabla^{\End{\mathcal{E}}}U-UA$ whose connection symbols are in $\rho^NC^\infty \left(\overline{M} \right)$ in any boundary coordinates and frame (in the sense of Definition \ref{MathItem_2.7}).
Thus (\ref{MathItem_5.86}) can be rewritten as
\begin{equation} \label{MathItem_5.87}  \left(\pi ^\ast\nabla^\prime \right)_XW+\Psi \left(W \right)=f\end{equation}
where $f$ is the right-hand side of (\ref{MathItem_5.86}) (note that $\Psi \left(W \right)$ is not a matrix multiplication but rather $\Psi$ applied to $U$).

Hence assuming that $N$ is large enough, by (\ref{MathItem_5.87}) and Proposition \ref{MathItem_5.56} and we have that $W$ is in $\mathcal{R}^3 \left(SM;\pi ^\ast\End{\mathcal{E}} \right)$.
Furthermore, by Remark \ref{MathItem_5.6} we have that $W$ extends continuously to ${}^0S^\ast\overline{M}$ and vanishes on the boundary $\partial {}^0S^\ast\overline{M}$.
Next, looking at (\ref{MathItem_5.87}) in coordinates and frame, above any fixed point $x\in M$ the entries of the right-hand side $f$ are restrictions of homogeneous polynomials of order zero and one in the variable $v$.
By the theory of spherical harmonics these are elements of (Fourier) degree zero and one respectively (see Section 2.H in \cite{Bibitem_13}).
Since we assumed in the theorem statement that $\nabla^\prime$ has no nontrivial twisted CKTs in $\mathcal{R}^3 \left(SM;\pi ^\ast\End{\mathcal{E}} \right)$, by Proposition \ref{MathItem_5.55} we have that $W$ and hence $Q$ are of degree zero.
Thus, we get that $Q\in C^\infty \left(M;\End{\mathcal{E}} \right)\cap C^0 \left(\overline{M};\End{\mathcal{E}} \right)$ (i.e. up to identification) such that $\left.Q\right|_{\partial \overline{M}}=\mathrm{id}$ and $ \left(Q-\mathrm{id} \right)\in\mathcal{R}^3 \left(SM;\pi ^\ast\mathcal{E} \right)$.

As the final step, let us show that this $Q$ is the gauge that we wanted.
In our coordinates above, we have that
$$\mathbb{X}Q=v \left(Q \right)+ \left({}^\mathcal{E}\Gamma \right)Q-Q \left({}^\mathcal{E}\Gamma \right),$$
whose entries are first order polynomials in $v$.
Plugging this into (\ref{MathItem_5.85}) and equating the zeroth and first Fourier modes gives
\begin{equation} \label{MathItem_5.88} \Phi Q-Q\widetilde{\Phi}=0\ \ \ \ \ \ \ \ \ \ \mathrm{and} \ \ \ \ \ \ \ \ \ \mathbb{X}Q-QA=0,\end{equation}
The first equation gives $\widetilde{\Phi}=Q^{-1}\Phi Q$ over $M$ and hence over $\overline{M}$ by continuity.
The second equation gives
$$A=Q^{-1}\mathbb{X}Q=Q^{-1}\left[v \left(Q \right)+ \left({}^\mathcal{E}\Gamma \right)Q-Q \left({}^\mathcal{E}\Gamma \right)\right],$$
$$\Longrightarrow\ \ \ \ \ \ \ \ A+ \left({}^\mathcal{E}\Gamma \right)=Q^{-1}v \left(Q \right)+Q^{-1} \left({}^\mathcal{E}\Gamma \right)Q,$$
Now, take any section $u\in C^\infty(\overline{M};\mathcal{E}$) and write it as a column vector with respect to the basis $ \left(b_i \right)$.
We have that (here $v \left(u \right)$ denotes applying $v$ to every entry of $u$)
\begin{align*}{\widetilde{\nabla}}_v^\mathcal{E}u=v \left(u \right)+\left[A+ \left({}^\mathcal{E}\Gamma \right)\right]u& &\mathrm{eq.\ for}\ {\widetilde{\nabla}}_v^\mathcal{E},\\=v \left(Q^{-1}Qu \right)+Q^{-1}v \left(Q \right)u+Q^{-1} \left({}^\mathcal{E}\Gamma \right)Qu& &QQ^{-1}=\mathrm{id} \ \mathrm{and\ plug\ in\ eq.\ above} ,\\=v \left(Q^{-1}Qu \right)-v \left(Q^{-1} \right)Qu+Q^{-1} \left({}^\mathcal{E}\Gamma \right)Qu& &Q^{-1}v \left(Q \right)=-v \left(Q^{-1} \right)Q\ \mathrm{(prod.\ rule)} ,\\=Q^{-1}v \left(Qu \right)+Q^{-1} \left({}^\mathcal{E}\Gamma \right)Qu& &\mathrm{prod.\ rule},\\=Q^{-1}\nabla_v^\mathcal{E} \left(Qu \right)& &\mathrm{eq.\ for}\ \nabla_v^\mathcal{E}.\end{align*}
Again this relation extends to $\partial \overline{M}$ by continuity and hence this proves (\ref{MathItem_2.9}).

Lastly, let us show that $Q$ is unitary.
This is trivial on $\partial \overline{M}$ since $Q=\mathrm{id}$ there.
Over $SM$, since $Q=U{\widetilde{U}}^{-1}$, this will follow if we prove that $U$ and $\widetilde{U}$ are unitary.
We start with $U$.
Choose any $e_0\in \left(\pi ^\ast\mathcal{E} \right)_{ \left(x_0,v_0 \right)}$ where $v_0\in S_{x_0}M$.
We will show that $\left|Ue_0\right|_\mathcal{E}=\left|e_0\right|_\mathcal{E}$.
Consider the $g$-geodesic $\gamma :  \left(-\infty,\infty \right)\rightarrow M$ satisfying $\gamma \left(0 \right)=x_0$ and $\gamma^\prime \left(0 \right)=v_0$ and its lift $\sigma :  \left(-\infty,\infty \right)\rightarrow SM$ to $SM$ which passes through $v_0$ at $t=0$ (i.e. $\gamma=\pi \circ\sigma$).
By Lemma \ref{MathItem_2.6} parts 1) and 2), there exists a parallel section $E :  \left(-\infty,\infty \right)\rightarrow\mathcal{E}$ (i.e. $\nabla_{\dot{\gamma}}^\mathcal{E}E\equiv0$) such that $\lim_{t\rightarrow-\infty}{E \left(t \right)}=e_{-\infty}$ exists and $E \left(0 \right)=e_0$.
Notice that since $\sfrac{d}{dt}\left|E\right|_\mathcal{E}^2=2\langle \nabla_{\dot{\gamma}}^\mathcal{E}E,E\rangle =0$, we have that $\left|e_{-\infty}\right|_\mathcal{E}^2=\left|e_0\right|_\mathcal{E}^2$.
Now, by the calculation in (\ref{MathItem_5.83}) we have that
$$\frac{d}{dt}\left| \left(U\circ\sigma \right)E\right|_\mathcal{E}^2=2\langle \nabla_{\dot{\gamma}}^\mathcal{E}\left[ \left(U\circ\sigma \right)E\right], \left(U\circ\sigma \right)E\rangle _\mathcal{E}=\langle -\Phi \left(U\circ\sigma \right)E, \left(U\circ\sigma \right)E\rangle _\mathcal{E}=0,$$
where in the last equality we have used that $\Phi$ is skew-Hermitian: $\langle -\Phi U E,UE\rangle _\mathcal{E}=\langle UE,\Phi U E\rangle _\mathcal{E}=\langle \Phi U E,UE\rangle _\mathcal{E}$ and hence is zero.
By (\ref{MathItem_5.81}) we have that $U$ is the identity on $\partial _-{}^bS^\ast\overline{M}$ and so $\left| \left(U\circ\sigma \right)E\right|_\mathcal{E}^2\rightarrow\left|e_{-\infty}\right|^2$ as $t\rightarrow-\infty$.
These two observations imply that
$$\left|Ue_0\right|_\mathcal{E}^2=\left| \left(U\circ\sigma \right)E \left(0 \right)\right|_\mathcal{E}^2=\left|e_{-\infty}\right|_\mathcal{E}^2=\left|e_0\right|_\mathcal{E}^2.$$
So $U$ is indeed unitary over $SM$.
A similar proof using (\ref{MathItem_5.84}) shows that $\widetilde{U}$ is unitary, and so as discussed above we get that $Q$ is indeed unitary.
This fact extends to $\partial \overline{M}$ by continuity.

\begin{flushright}$\blacksquare$\end{flushright}
\section{Proof of Theorem \ref{MathItem_2.15}} \label{Section_6} 

Suppose that $u\in\mathcal{R}^3 \left(SM;\pi ^\ast\mathcal{E} \right)$ is a nontrivial twisted CKT of degree $m\geq1$.
We will show that it is identically zero everywhere.
We already proved the existence of a $\kappa$ as in the statement of the lemma in the proof of Lemma \ref{MathItem_5.44} above.
Recalling that $\mathbb{X}_+u=0$ by definition, the idea here is to use our curvature bounds and Proposition \ref{MathItem_5.43} to conclude that $\mathbb{X}_-u=0$ as well.
To do this, observe that plugging (\ref{MathItem_5.47}) and (\ref{MathItem_5.48}) into the equation in Proposition \ref{MathItem_5.43} gives
$$ \left(2m+n \right)\left\lVert \mathbb{X}_+u\right\rVert _{L^2}^2\geq\left\lVert {\buildrel\mathrm{h}\over\nabla}{}^{\pi ^\ast\mathcal{E}}u\right\rVert _{L^2}^2+ \left(2m+n-2 \right)\left\lVert \mathbb{X}_-u\right\rVert _{L^2}^2+ \left(\kappa\lambda_m-\left\lVert F^\mathcal{E}\right\rVert _{L^\infty}\lambda_m^{\sfrac{1}{2}} \right)\left\lVert u\right\rVert _{L^2}^2.$$
Recalling the definition of $\lambda_m$, our assumption $\left\lVert F^\mathcal{E}\right\rVert _{L^\infty}\le\kappa\sqrt n$ implies that the coefficient of $\left\lVert u\right\rVert _{L^2}^2$ on the right-hand side is nonnegative.
Hence both sides of the inequality are nonnegative, and thus the left-hand side being zero implies that $\mathbb{X}_-u=0$.
Hence $\mathbb{X}u=\mathbb{X}_+u+\mathbb{X}_-u=0$ as well.

Now, choose any point $v_0\in S_xM$.
We will show that $u \left(v_0 \right)=0$.
Consider the $g$-geodesic $\gamma :  \left(-\infty,\infty \right)\rightarrow M$ satisfying $\gamma \left(0 \right)=x_0$ and $\gamma^\prime \left(0 \right)=v_0$ and its lift $\sigma :  \left(-\infty,\infty \right)\rightarrow SM$ to $SM$ which passes through $v_0$ at $t=0$ (i.e. $\gamma=\pi \circ\sigma$ and $\sigma$ is an integral curve of $X$).
As we explained in Section \ref{Section_2.3} above, since $ \left(M,g \right)$ is nontrapping the limit $\lim_{t\rightarrow-\infty}{\gamma \left(t \right)}\in\partial \overline{M}$ exists.
We have that $u$ vanishes on the compact boundary $\partial {}^0S^\ast\overline{M}$ by Remark \ref{MathItem_5.6} and hence it follows that that $u\circ\sigma \left(t \right)\rightarrow0$ as $t\rightarrow-\infty$.
This combined with the fact that $\mathbb{X}u\equiv0$ imply that $u\circ\sigma \left(t \right)\equiv0$ and so $u \left(v_0 \right)=0$.
Thus indeed $u$ vanishes everywhere.

We point out that the conclusion of the theorem follows even faster if one assumes $\left\lVert F^\mathcal{E}\right\rVert _{L^\infty}<\kappa\sqrt n$ since then the coefficient of $\left\lVert u\right\rVert _{L^2}^2$ in the above inequality is positive, and hence $\left\lVert u\right\rVert _{L^2}^2=0$ giving $u\equiv0$.

\begin{flushright}$\blacksquare$\end{flushright}
\section{Injectivity over Higgs Fields} \label{Section_7} 

In this section we prove Corollary \ref{MathItem_2.11}. 

\begin{lemma} \label{MathItem_7.1}  Suppose that $\overline{M}$ is a smooth manifold with smooth boundary, $\mathcal{E}$ is a smooth complex vector bundle over $\overline{M}$, and that $\nabla^\mathcal{E}$ is a smooth connection in $\mathcal{E}$.
Let $\nabla^\prime$ be as defined in Theorem \ref{MathItem_2.8} with ${\widetilde{\nabla}}^\mathcal{E}=\nabla^\mathcal{E}$.
If $\nabla^\mathcal{E}$ has zero curvature, then so does $\nabla^\prime$.
\end{lemma}

\noindent\textbf{Proof:} Since in this case $\nabla^\mathcal{E}={\widetilde{\nabla}}^\mathcal{E}$, the $A$ in Theorem \ref{MathItem_2.8} is zero and hence $\nabla^\prime=\nabla^{\End{\mathcal{E}}}$.
So we must simply show that $\nabla^{\End{\mathcal{E}}}$ has zero curvature.
Let $ \left(x^i \right)$ be coordinates of $\overline{M}$ and $ \left(b_i \right)$ a frame for $\mathcal{E}$ over these coordinates' domain.
Let ${}^\mathcal{E}\Gamma_{ij}^k$ denote the connection symbols of $\nabla^\mathcal{E}$ with respect to $ \left(\sfrac{\partial }{\partial  x^i} \right)$ and $ \left(b_i \right)$.
Because the curvature of $\nabla^\mathcal{E}$ is zero, by Proposition 1.2 of Appendix C in Volume II of \cite{Bibitem_51}, we may choose the $ \left(b_i \right)$ so that all of the ${}^\mathcal{E}\Gamma_{ij}^k\equiv0$.
Let us write sections of $\End{\mathcal{E}}$ as square matrices with respect to this frame $ \left(b_i \right)$ (i.e. write $U : \overline{M}\rightarrow\End{\mathcal{E}}$ as the matrix $\left[U_i^j\right]_{ij=1}^n$ where $U \left(\alpha^ib_i \right)=U_i^j\alpha^ib_j$).
Then, if we let $E_j^i$ denote the matrix with all zero entries except for a one in the $i^{\mathrm{th}}$ row and $j^{\mathrm{th}}$ column, $ \left(E_j^i \right)$ is a smooth frame of $\End{\mathcal{E}}$.
Observe that we may write any $U : \overline{M}\rightarrow\End{\mathcal{E}}$ locally as $U=U_i^jE_j^i$.

Now, from the explicit equation (\ref{MathItem_4.26}) for $\nabla^{\End{\mathcal{E}}}$ we have that
$$\nabla_v^{\End{\mathcal{E}}}U=v \left(U_i^j \right)E_j^i.$$
In other words, the connection symbols of $\nabla^{\End{\mathcal{E}}}$ with respect to $ \left(\sfrac{\partial }{\partial  x^i} \right)$ and $ \left(E_j^i \right)$ are zero and hence by (\ref{MathItem_4.23}) has curvature zero.

\begin{flushright}$\blacksquare$\end{flushright}

\noindent\textit{Proof of Corollary \ref{MathItem_2.11}:}

Recall that $d=\rank{\mathcal{E}}$.
By Lemma \ref{MathItem_7.1} we have that the curvature of $\nabla^{\End{\mathcal{E}}}$ is zero.
Hence, by plugging $\End{\mathcal{E}}$ into $\mathcal{E}$ in Theorem \ref{MathItem_2.15} we have that $\nabla^{\End{\mathcal{E}}}$ has no nontrivial twisted CKTs in $\mathcal{R}^3 \left(SM;\pi ^\ast\End{\mathcal{E}} \right)$.
Thus, by Theorem \ref{MathItem_2.8}, if $N\geq0$ is big enough there exists a unitary $Q\in C^0 \left(\overline{M};\End{\mathcal{E}} \right)\cap C^\infty \left(M;\End{\mathcal{E}} \right)$ such that $\left.Q\right|_{\partial \overline{M}}=\mathrm{id}$ and satisfies (\ref{MathItem_2.9}) with ${\widetilde{\nabla}}^\mathcal{E}=\nabla^\mathcal{E}$.
So, we will be done if we show that $Q\equiv\mathrm{id}$ everywhere.
Take coordinates $ \left(x^i \right)$ of $\overline{M}$, a frame $ \left(b_i \right)$ for $\mathcal{E}$ over these coordinates' domain, and consider the coordinates $\sfrac{v^i\partial }{\partial  x^i}\mapsto \left(x^i,v^i \right)$ of $T\overline{M}$.
Let ${}^\mathcal{E}\Gamma_{ij}^k$ denote the connection symbols of $\nabla^\mathcal{E}$ with respect to $ \left(\sfrac{\partial }{\partial  x^i} \right)$ and $ \left(b_k \right)$.
Since the curvature of $\nabla^\mathcal{E}$ is zero, by Proposition 1.2 of Appendix C in Volume II of \cite{Bibitem_51}, we may suppose that the $ \left(b_i \right)$ were chosen so that all of the ${}^\mathcal{E}\Gamma_{ij}^k\equiv0$.
Let us represent $Q$ as a $d\times d$ matrix in the basis $ \left(b_i \right)$.
Similarly we represent any section $u\in C^\infty \left(\overline{M};\mathcal{E} \right)$ as a $d\times1$ column vector in this basis.
Then for any section $u=u^jb_j\cong\left[u^1,\ldots,u^d\right]$ whose component functions $u^j$ are constant, we have that for any $v\in TM$ in our coordinates
$$\nabla_v^\mathcal{E}u=0,$$
$${\widetilde{\nabla}}_v^\mathcal{E}u=Q^{-1}\nabla_v^\mathcal{E} \left(Qu \right)=Q^{-1}v \left(Q \right)u.$$
where $v \left(Q \right)$ denotes applying $v$ to the entries of $Q$.
Since ${\widetilde{\nabla}}_v^\mathcal{E}u=\nabla_v^\mathcal{E}u$ by assumption and the above is true for all such $u$, we get that $v \left(Q \right)\equiv0$.
Hence $Q$ is locally constant.

Thus the sets $\left\{x\in\overline{M} : Q=\mathrm{id} \right\}$ and $\left\{x\in\overline{M} : Q\neq\mathrm{id} \right\}$ are both open.
Observe that because of nontrapping, for any point $x\in M$ and any geodesic $\gamma : \mathbb{R}\rightarrow M$ passing through $x$ has limits $\lim_{t\rightarrow\pm\infty}{\gamma \left(t \right)}\in\partial \overline{M}$.
From this it follows that the intersection of every connected component of $\overline{M}$ with $\partial  M$, and hence with $\left\{x\in\overline{M} : Q=\mathrm{id} \right\}$, is nonempty.
Thus indeed $Q$ is equal to ``id" everywhere.

\begin{flushright}$\blacksquare$\end{flushright}


\begin{thebibliography}{99} \bibitem{Bibitem_1} Yu. E. Anikonov and V. G. Romanov. On uniqueness of determination of a form of first degree by its integrals along geodesics. \textit{J. Inverse Ill-Posed Probl.}, 5(6):487-490, 1997.

	\bibitem{Bibitem_2} Carlos A. Berenstein and Enrico Casadio Tarabusi. Inversion formulas for the k-dimensional Radon transform in real hyperbolic spaces. \textit{Duke Math. J.}, 62(3):613-631, 1991.

	\bibitem{Bibitem_3} Jan Bohr. Stability of the non-abelian X-ray transform in dimension $\geq3$. \textit{J. Geom. Anal.}, 31(11):11226-11269, 2021.

	\bibitem{Bibitem_4} Jan Bohr and Gabriel P. Paternain. The transport Oka-Grauert principle for simple surfaces. \textit{J. Ec. polytech. Math.}, 10:727-769, 2023.

	\bibitem{Bibitem_5} Mihajlo Cekić and Thibault Lefeuvre. Generic dynamical properties of connections on vector bundles. \textit{Int. Math. Res. Not. IMRN}, (14):10649 - 10703, 2022.

	\bibitem{Bibitem_6} Earl A. Coddington and Norman Levinson. \textit{Theory of ordinary differential equations}. McGraw-Hill Book Co., Inc., New York-Toronto-London, 1955.

	\bibitem{Bibitem_7} Bartlomiej Czech, Lampros Lamprou, Samuel McCandlish, and James Sully. Integral geometry and holography. \textit{J. High Energy Phys.}, (10):175, front matter+40, 2015.

	\bibitem{Bibitem_8} N. S. Dairbekov and V. A. Sharafutdinov. On conformal Killing symmetric tensor fields on Riemannian manifolds. \textit{Mat. Tr.}, 13(1):85-145, 2010.

	\bibitem{Bibitem_9} Naeem M. Desai, William R. B. Lionheart, Morten Sales, Markus Strobl, and S{\o}ren Schmidt. Polarimetric neutron tomography of magnetic fields: uniqueness of solution and reconstruction. \textit{Inverse Problems}, 36(4):045001, 17, 2020.

	\bibitem{Bibitem_10} Nikolas Eptaminitakis. Stability estimates for the X-ray transform on simple asymptotically hyperbolic manifolds. \textit{Pure Appl. Anal.}, 4(3):487-516, 2022.

	\bibitem{Bibitem_11} Nikolas Eptaminitakis and C. Robin Graham. Local X-ray transform on asymptotically hyperbolic manifolds via projective compactification. \textit{New Zealand J. Math.}, 52:733-763, 2021 [2021-2022].

	\bibitem{Bibitem_12} David Finch and Gunther Uhlmann. The X-ray transform for a non-abelian connection in two dimensions. volume 17, pages 695-701. 2001. Special issue to celebrate Pierre Sabatier's 65th birthday (Montpellier, 2000).

	\bibitem{Bibitem_13} Gerald B. Folland. \textit{Introduction to partial differential equations}. Princeton University Press, Princeton, NJ, second edition, 1995.

	\bibitem{Bibitem_14} Gerald B. Folland. \textit{Real analysis}. Pure and Applied Mathematics (New York). John Wiley \& Sons, Inc., New York, second edition, 1999. Modern techniques and their applications, A Wiley-Interscience Publication.

	\bibitem{Bibitem_15} C. Robin Graham, Colin Guillarmou, Plamen Stefanov, and Gunther Uhlmann. X-ray transform and boundary rigidity for asymptotically hyperbolic manifolds. \textit{Ann. Inst. Fourier (Grenoble)}, 69(7):2857-2919, 2019.

	\bibitem{Bibitem_16} C. Robin Graham and John M. Lee. Einstein metrics with prescribed conformal infinity on the ball. \textit{Adv. Math.}, 87(2):186-225, 1991.

	\bibitem{Bibitem_17} Colin Guillarmou. Lens rigidity for manifolds with hyperbolic trapped sets. \textit{J. Amer. Math. Soc.}, 30(2):561-599, 2017.

	\bibitem{Bibitem_18} Colin Guillarmou, Gabriel P. Paternain, Mikko Salo, and Gunther Uhlmann. The X-ray transform for connections in negative curvature. \textit{Comm. Math. Phys.}, 343(1):83-127, 2016.

	\bibitem{Bibitem_19} Victor Guillemin and David Kazhdan. Some inverse spectral results for negatively curved n-manifolds. In \textit{Geometry of the Laplace operator (Proc. Sympos. Pure Math., Univ. Hawaii, Honolulu, Hawaii, 1979)}, volume XXXVI of \textit{Proc. Sympos. Pure Math.}, pages 153-180. Amer. Math. Soc., Providence, RI, 1980.

	\bibitem{Bibitem_20} Sigurdur Helgason. The totally-geodesic Radon transform on constant curvature spaces. In \textit{Integral geometry and tomography (Arcata, CA, 1989)}, volume 113 of \textit{Contemp. Math.}, pages 141-149. Amer. Math. Soc., Providence, RI, 1990.

	\bibitem{Bibitem_21} Sigurdur Helgason. \textit{Geometric analysis on symmetric spaces}, volume 39 of \textit{Mathematical Surveys and Monographs}. American Mathematical Society, Providence, RI, 1994.

	\bibitem{Bibitem_22} A. Hilger, I. Manke, N. Kardjilov, M. Osenberg, H. Mark\"otter, and J. Banhart. Tensorial neutron tomography of three-dimensional magnetic vector fields in bulk materials. \textit{Nat Commun}, 9(4023), 2018.

	\bibitem{Bibitem_23} Sean Holman and Gunther Uhlmann. On the microlocal analysis of the geodesic X-ray transform with conjugate points. \textit{J. Differential Geom.}, 108(3):459-494, 2018.

	\bibitem{Bibitem_24} Joonas Ilmavirta. Coherent quantum tomography. \textit{SIAM J. Math. Anal.}, 48(5):3039-3064, 2016.

	\bibitem{Bibitem_25} John M. Lee. \textit{Introduction to smooth manifolds}, volume 218 of \textit{Graduate Texts in Mathematics}. Springer, New York, second edition, 2013.

	\bibitem{Bibitem_26} John M. Lee. \textit{Introduction to Riemannian manifolds}, volume 176 of \textit{Graduate Texts in Mathematics}. Springer, Cham, second edition, 2018.

	\bibitem{Bibitem_27} Jere Lehtonen. \textit{The geodesic ray transform on two-dimensional Cartan Hadamard manifolds}. PhD thesis, University of Jyv\"askyl\"a, 2016.

	\bibitem{Bibitem_28} Jere Lehtonen, Jesse Railo, and Mikko Salo. Tensor tomography on Cartan-Hadamard manifolds. \textit{Inverse Problems}, 34(4):044004, 27, 2018.

	\bibitem{Bibitem_29} Juan Maldacena. The large N limit of superconformal field theories and supergravity [MR1633016 (99e:81204a)]. In \textit{Trends in theoretical physics, II (Buenos Aires, 1998)}, volume 484 of \textit{AIP Conf. Proc.}, pages 51-63. Amer. Inst. Phys., Woodbury, NY, 1999.

	\bibitem{Bibitem_30} S. V. Manakov and V. E. Zakharov. Three-dimensional model of relativistic-invariant field theory, integrable by the inverse scattering transform. \textit{Lett. Math. Phys.}, 5(3):247-253, 1981.

	\bibitem{Bibitem_31} Rafe Mazzeo. \textit{Hodge Cohomology of Negatively Curved Manifolds}. PhD thesis, Massachusetts Institute of Technology, 1986.

	\bibitem{Bibitem_32} Richard B. Melrose. \textit{The Atiyah-Patodi-Singer index theorem}, volume 4 of \textit{Research Notes in Mathematics}. A K Peters, Ltd., Wellesley, MA, 1993.

	\bibitem{Bibitem_33} Francois Monard, Richard Nickl, and Gabriel P. Paternain. Efficient nonparametric Bayesian inference for X-ray transforms. \textit{The Annals of Statistics}, 47(2):1113 - 1147, 2019.

	\bibitem{Bibitem_34} Francois Monard, Plamen Stefanov, and Gunther Uhlmann. The geodesic ray transform on Riemannian surfaces with conjugate points. \textit{Comm. Math. Phys.}, 337(3):1491-1513, 2015.

	\bibitem{Bibitem_35} R. G. Muhometov. The problem of recovery of a two-dimensional Riemannian metric and integral geometry. Doklady Akademii Nauk, 232(1):32-35, 1977.

	\bibitem{Bibitem_36} R. G. Novikov. Non-Abelian Radon transform and its applications. In \textit{The Radon transform—the first 100 years and beyond}, volume 22 of \textit{Radon Ser. Comput. Appl. Math.}, pages 115-127. Walter de Gruyter, Berlin, Copyright 2019.

	\bibitem{Bibitem_37} R. G. Novikov. On determination of a gauge field on $\mathbb{R}^d$ from its non-abelian Radon transform along oriented straight lines. \textit{J. Inst. Math. Jussieu}, 1(4):559-629, 2002.

	\bibitem{Bibitem_38} Gabriel P. Paternain, Mikko Salo, and Gunther Uhlmann. The attenuated ray transform for connections and Higgs fields. \textit{Geom. Funct. Anal.}, 22(5):1460-1489, 2012.

	\bibitem{Bibitem_39} Gabriel P. Paternain, Mikko Salo, and Gunther Uhlmann. Tensor tomography on surfaces. \textit{Invent. Math.}, 193(1):229-247, 2013.

	\bibitem{Bibitem_40} Gabriel P. Paternain, Mikko Salo, and Gunther Uhlmann. Tensor tomography: progress and challenges. \textit{Chinese Ann. Math. Ser. B}, 35(3):399-428, 2014.

	\bibitem{Bibitem_41} Gabriel P. Paternain, Mikko Salo, and Gunther Uhlmann. Invariant distributions, Beurling transforms and tensor tomography in higher dimensions. \textit{Math. Ann.}, 363(1-2):305-362, 2015.

	\bibitem{Bibitem_42} Gabriel P. Paternain, Mikko Salo, and Gunther Uhlmann. \textit{Geometric inverse problems - with emphasis on two dimensions}, volume 204 of \textit{Cambridge Studies in Advanced Mathematics}. Cambridge University Press, Cambridge, 2023. With a foreword by Andr\'as Vasy.

	\bibitem{Bibitem_43} L. N. Pestov and V. A. Sharafutdinov. Integral geometry of tensor fields on a manifold of negative curvature. \textit{Sibirsk. Mat. Zh.}, 29(3):114-130, 221, 1988.

	\bibitem{Bibitem_44} V. A. Sharafutdinov. \textit{Integral geometry of tensor fields}. Inverse and Ill-posed Problems Series. VSP, Utrecht, 1994.

	\bibitem{Bibitem_45} V. A. Sharafutdinov. On the inverse problem of determining a connection on a vector bundle. \textit{J. Inverse Ill-Posed Probl.}, 8(1):51-88, 2000.

	\bibitem{Bibitem_46} V. A. Sharafutdinov. Variations of Dirichlet-to-Neumann map and deformation boundary rigidity of simple 2-manifolds. \textit{J. Geom. Anal.}, 17(1):147-187, 2007.

	\bibitem{Bibitem_47} Plamen Stefanov and Gunther Uhlmann. Boundary rigidity and stability for generic simple metrics. \textit{J. Amer. Math. Soc.}, 18(4):975-1003, 2005.

	\bibitem{Bibitem_48} Plamen Stefanov and Gunther Uhlmann. Local lens rigidity with incomplete data for a class of non-simple Riemannian manifolds. \textit{J. Differential Geom.}, 82(2):383-409, 2009.

	\bibitem{Bibitem_49} Plamen Stefanov and Gunther Uhlmann. The geodesic X-ray transform with fold caustics. \textit{Anal. PDE}, 5(2):219-260, 2012.

	\bibitem{Bibitem_50} Plamen Stefanov, Gunther Uhlmann, and Andr\'as Vasy. Inverting the local geodesic X-ray transform on tensors. \textit{J. Anal. Math.}, 136(1):151-208, 2018.

	\bibitem{Bibitem_51} Michael E. Taylor. \textit{Partial differential equations}, volume 115 of \textit{Applied Mathematical Sciences}. Springer, New York, second edition, 2011.

	\bibitem{Bibitem_52} Gunther Uhlmann and Andr\'as Vasy. The inverse problem for the local geodesic ray transform. \textit{Invent. Math.}, 205(1):83-120, 2016.

	\bibitem{Bibitem_53} R. S. Ward. Soliton solutions in an integrable chiral model in 2 + 1 dimensions. \textit{J. Math. Phys.}, 29(2):386-389, 1988.\end{thebibliography}
\end{document}